\documentclass[aps, pre,twocolumn, superscriptaddress,longbibliography]{revtex4-1}

\pdfoutput=1   %This is for the ArXiv submission only

\usepackage{bm,graphicx,epsfig,subfigure,color,amsmath}
\usepackage{amssymb}
\usepackage{verbatim}

\usepackage[colorlinks=true, urlcolor=blue, anchorcolor=blue, citecolor=blue,filecolor=blue,linkcolor=blue,menucolor=blue]{hyperref}

\usepackage{subfigure}

\DeclareGraphicsExtensions{.jpg,.pdf, .mps, .png, .eps, .ps, .EPS,.gif}
\graphicspath{{./_figs_deterministic/}}

\begin{document}

\def\be{\begin{equation}}
\def\ee{\end{equation}}

\def\bc{\begin{center}}
\def\ec{\end{center}}
\def\bea{\begin{eqnarray}}
\def\eea{\end{eqnarray}}
\newcommand{\avg}[1]{\langle{#1}\rangle}
\newcommand{\Avg}[1]{\left\langle{#1}\right\rangle}
\newcommand{\Bavg}[1]{\Bigl\langle{#1}\Bigr\rangle}
\newcommand{\eq}[1]{(\ref{#1})}

\def\ie{\textit{i.e.}}
\def\etal{\textit{et al.}}
\def\m{\vec{m}}
\def\G{\mathcal{G}}
\def\RRR{\hbox{\msytw R}}
\def\CCC{\hbox{\msytw C}}
\def\NNN{\hbox{\msytw N}} 
\def\ZZZ{\hbox{\msytw Z}}
\def\TTT{\hbox{\msytw T}}

\title
%%%%%%%%%%%%%%[Deterministic simplicial complexes]
{Deterministic simplicial complexes}

\author{S. N. Dorogovtsev}
\email{sdorogov@ua.pt}
%%%%%%%%%%%%%%%%%%%%%%%%%%%%%%%%\address
\affiliation{Departamento de F\'{i}sica da Universidade de Aveiro \& I3N, Campus Universit\'{a}rio de Santiago, Aveiro, 3810-193, Portugal}
%%%%%%%%%%%%%%%%%%%%%%%%%%%%%%%%\address
\affiliation{A. F. Ioffe Physico-Technical Institute, St. Petersburg, 194021, Russia}

\author{P. L. Krapivsky}
\email{pkrapivsky@gmail.com}
%%%%%%%%%%%%%%%%%%%%%%%%%%%%%%%%\address
\affiliation{Department of Physics, Boston University, Boston, Massachusetts 02215, USA}
%%%%%%%%%%%%%%%%%%%%%%%%%%%%%%%%\address
\affiliation{Santa Fe Institute, Santa Fe, New Mexico 87501, USA}

%%%%%%%%%%%%%%%%%%%%%%%%%%%%%%\eads{\mailto{sdorogov@ua.pt}, ~\mailto{pkrapivsky@gmail.com}}

\vspace{10pt}

\begin{abstract}
We investigate simplicial complexes deterministically growing from a single vertex. In the first step, a vertex and an edge connecting it to the primordial vertex are added. The resulting simplicial complex has a 1-dimensional simplex and two 0-dimensional faces (the vertices). The process continues recursively: On the $n$-th step, every existing $d-$dimensional simplex ($d\leq n-1$) joins a new vertex forming a $(d+1)-$dimensional simplex; all $2^{d+1}-2$ new faces are also added so that the resulting object remains a simplicial complex. The emerging simplicial complex has intriguing local and global characteristics. The number of simplices grows faster than $n!$, and the upper-degree distributions follow a power law. Here, the upper degree (or $d$-degree) of a $d$-simplex refers to the number of $(d{+}1)$-simplices that share it as a face. Interestingly, the $d$-degree distributions evolve quite differently for different values of $d$. We compute the Hodge Laplacian spectra of simplicial complexes and show that the spectral and Hausdorff dimensions are infinite. We also explore a constrained version where the dimension of the added simplices is fixed to a finite value $m$. In the constrained model, the number of simplices grows exponentially. In particular, for $m=1$, the spectral dimension is $2$. For $m=2$, the spectral dimension is finite, and the degree distribution follows a power law, while the $1$-degree distribution decays exponentially. 
\end{abstract}
 
\maketitle 

%%\vspace{2pc}

%%\vspace{2pc}

\section{Introduction}
\label{s1}

Simplicial complexes are getting a lot of attention not only because they are central to algebraic topology \cite{Hatcher}, but also because they offer a natural way to represent higher-order interactions in complex systems.  Recall that a $0$-dimensional simplex is a vertex; a $1$-dimensional simplex is an edge with two vertices; a $2$-dimensional simplex is a triangle; a $3$-dimensional simplex is a tetrahedron; etc. A $d$-dimensional simplex ($d$-simplex) has $2^d-1$ non-empty sub-simplices, faces, in short: $\binom{d}{1}$ vertices, $\binom{d}{2}$ edges, $\binom{d}{3}$ triangles, etc. 
A simplicial complex $\mathcal{K}$ is a set of simplices satisfying two conditions: (i) every face of a simplex from $\mathcal{K}$ is also in $\mathcal{K}$ and (ii) non-empty intersection of any two simplices from $\mathcal{K}$ is a face of both simplices. A $d$-dimensional simplicial complex contains simplices of maximal dimension $d$. 

Large random simplicial complexes are gaining popularity, see, e.g., \cite{Pippenger,Linial16,Farber16} and \cite{Kahle,Petri20,Ginestra21, bobrowski2022random} for review. Models of randomly growing simplicial complexes based on preferential attachment were studied e.g. in Refs.~\cite{bianconi2015complex, bianconi2017emergent,da2018complex,xie2023combinatorial}. Deterministic simplicial complexes have the advantage of being easier to analyze and, when possible, to describe explicitly compared to their stochastically growing counterparts. Deterministic graphs (a graph is a $1$-dimensional simplicial complex, where, by definition, the dimension of a simplicial complex is the maximum dimension of its simplices) are already well explored \cite{barabasi2001deterministic,Sergey02,andrade2005apollonian,hwang2010spectral,hasegawa2013hierarchical,dorogovtsev2022the}, and their structural characteristics and spectra were found to be remarkably close to their random counterparts. A detailed comparison between a deterministic graph and its stochastically growing counterpart---previously described in \cite{dorogovtsev2001size}---can be found in \cite{Sergey02}. 

Deterministic simplicial complexes generated by progressive adding maximal simplices of a given fixed dimension with all their faces were studied in \cite{bianconi2018topological,bianconi2020spectral,reitz2020higher}. [A maximal simplex (facet) is a simplex not contained in a larger simplex.] The $1$-skeletons of these simplicial complexes are the Apollonian deterministic graphs \cite{andrade2005apollonian}, the pseudo-fractals \cite{Sergey02}, and their immediate generalizations. In this article, we introduce a distinct class of more elegant deterministic simplicial complexes, where, at each growth step, new maximal simplices of various dimensions are added simultaneously. In the first, parameter-free model we introduce, the range of dimensions of new maximal simplices grows without bound as the complex evolves, though each step always includes edges ($1$-dimensional simplices). We refer to this type of growth as `unconstrained' and denote the model by the acronym DSC, standing for deterministic simplicial complex. In more constrained versions, we fix the maximum dimension of new maximal simplices to a finite value $m$ at every step---the `constrained growth', which makes the analysis significantly easier. We denote these models by the acronym DSC$(m)$. Importantly, the simplicial complexes generated by all these models are compact in the sense that the Hausdorff dimension of their $1$-skeletons is infinite.   

For these models of deterministic simplicial complexes, we derive key structural characteristics. First, we determine the exact number $N_d(n)$ of $d$-simplices at any generation $n$, along with the total number of simplices $N(n) \equiv \sum_d N_d(n)$. For the unconstrained growth, $N_d(n)$ increases as $n!$ or faster. For the constrained growth, $N_d(n)$ increases exponentially with $n$. We also probe the degree sequences of the simplicial complexes, more precisely the upper degrees where the upper degree $k^{(d)}$ of a $d$-simplex is the number of $(d{+}1)$-simplices sharing this $d$-simplex as a face. For brevity, we refer to the upper degree of a $d$-simplex as its $d$-degree. In particular, $k^{(0)}$ is the ordinary degree of a vertex ($0$-simplex). The distributions of $d$-degrees in our simplicial complexes turn out to be markedly different for each $d$. Finally, we describe the Laplacian spectra of these simplicial complexes. When it comes to Laplacians in the context of simplicial complexes,  the Hodge Laplacians are particularly useful \cite{hodge1934dirichlet,Dodziuk-hodge,lim2020hodge, schaub2020random}, where the $0$-th Hodge Laplacian of a simplicial complex is the Laplacian of the $1$-skeleton of this simplicial complex. While spectral properties of graphs are a major focus in graph theory \cite{brouwer2011spectra}, the spectra of Hodge Laplacians on simplicial complexes are much less explored. We obtain the spectra of all $d$-th Hodge Laplacians of our simplicial complexes, including the spectral dimensions associated with these matrices. 

The paper is organized as follows. Section~\ref{sec:null} introduces the basic, parameter-free DSC model of deterministic simplicial complexes. We investigate the evolution of the number of $d$-simplices, the upper-degree distributions, and present the spectra of the Hodge Laplacians. In section~\ref{s3}, we study the DSC$(m)$ models, focusing on the cases $m=1$ and $m=2$.  We describe the growth of the numbers of $d$-simplices, the upper-degree distributions, and the Hodge Laplacian spectra for these models. We also obtain the adjacency matrix spectrum for the DSC$(1)$ model, which generates trees. Section~\ref{s-conclusions} summarizes our main findings. Appendices include lengthy expressions for upper-degree sequences and characteristic polynomials, which support our derivations and help verify the results presented in the main text.

\section{Unconstrained growth}
\label{sec:null}

\subsection{The DSC model}

\begin{figure}[t]
\begin{center}
\includegraphics[width=0.02745\textwidth]{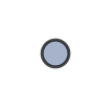}
\\[30pt]
\includegraphics[width=0.297\textwidth]{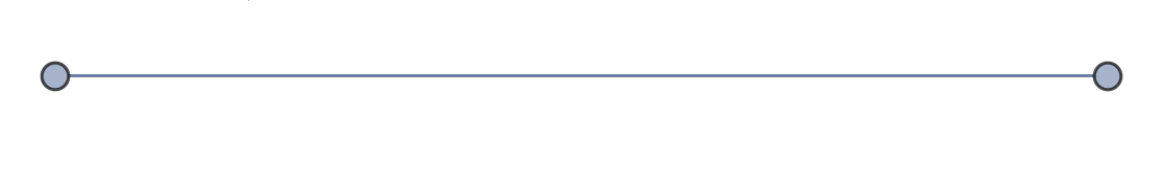}
\\[14pt]
\includegraphics[width=0.3555\textwidth]{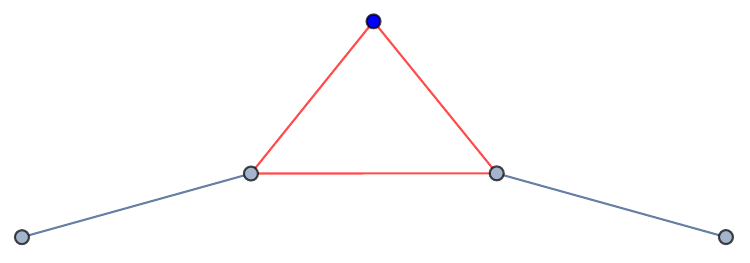}
\\[10pt]
\includegraphics[width=0.432\textwidth]{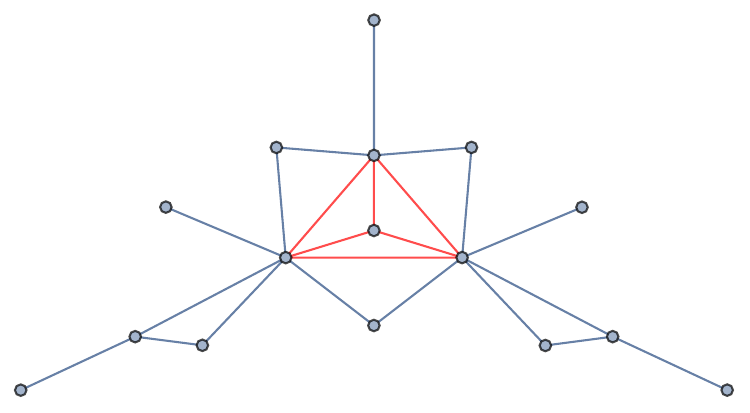}
\caption{The first four simplicial complexes generated by the DSC model, shown from top to bottom: $\mathcal{K}(0)$, $\mathcal{K}(1)$, $\mathcal{K}(2)$, $\mathcal{K}(3)$. 
} 
\label{fig:2}
\end{center}
\end{figure}

Here we analyze the basic model of deterministic simplicial complexes (DSC) explained in figure~\ref{fig:2}. The process begins with a primordial simplicial complex that is a single vertex: $\mathcal{K}(0)=[\{v_0\}]$. On the first step, a new vertex $\{v_1\}$  and an edge $\{v_0, v_1\}$  are added giving the simplicial complex  $\mathcal{K}(1)=[S_0(1), S_1(1)]$ with two $0$-dimensional faces 
\begin{equation}
S_0(1) = \{v_0\}, \{v_1\}
\label{10}
\end{equation}
and the $1$-dimensional simplex 
\begin{equation}
S_1(1) = \{v_0, v_1\}
.
\label{20}
\end{equation}
On the second step, every simplex in the simplicial complex $\mathcal{K}(1)$ is joined with a new vertex generating the simplicial complex $\mathcal{K}(2)=[S_0(1), S_1(1), S_2(1)]$ with five $0$-dimensional faces (vertices) 
\begin{equation}
S_0(2)=\{v_0\}, \{v_1\}, \{v_2\},  \{v_3\}, \{v_4\}
\label{30}
\end{equation}
and five $1$-dimensional faces (edges) 
\begin{equation}
S_1(2)=\{v_0, v_1\}, \{v_0, v_2\},  \{v_1, v_2\},  \{v_0, v_3\},  \{v_1, v_4\}
\label{40}
\end{equation}
together with the $2$-dimensional simplex (triangle)
\begin{equation}
S_2(2) = \{v_0, v_1,v_2\}
.
\label{50}
\end{equation}
On the third step, we obtain the simplicial complex $\mathcal{K}(3)$. This simplicial complex has a single $3$-dimensional simplex, the tetrahedron, and $48$ faces of smaller dimensions:  $16$ vertices, $23$ edges, and $9$ triangles, see figure~\ref{fig:2}. 

%%%%%%%%%%%%%%%%%%%%%%%%%%%%%%%%%%%%%%%%%%%%%%%
%%%%%%%%%%%%%%%%%%%%%%%%%%%%%%%%%%%%%%%%%%%%%%%

%%%%%%%%%%%%%%%%%%%%%%%%%%%%%%%%%%%%%%%%%%%%%%%
%%%%%%%%%%%%%%%%%%%%%%%%%%%%%%%%%%%%%%%%%%%%%%%

A $d$-dimensional simplex $\sigma_d$ is specified by its vertices: $\sigma_d=\langle w_0,w_1,\ldots,w_d\rangle$. Adding a new vertex $w_{d+1}$ transforms $\sigma_d$ to the $(d+1)$-dimensional simplex 
\begin{equation*}
\sigma_d=\langle w_0,w_1,\ldots,w_d\rangle \to \sigma_{d+1}=\langle w_0,w_1,\ldots,w_d,w_{d+1}\rangle
.
\label{60}
\end{equation*}

The simplicial complex $\mathcal{K}(n)$ is $n$-dimensional, so it admits the representation $\mathcal{K}(n)=[S_0(n), S_1(n),\ldots, S_n(n)]$. Denote by $N_d(n)=|S_d(n)|$ the number of simplices of dimension $d$. We intend to determine the vector 
\begin{equation}
\label{decomp}
\mathcal{N}(n)=[N_0(n),N_1(n),\ldots, N_n(n)]
.
\end{equation}
For a $d$-dimensional simplicial complex $\mathcal{K}$, a similar vector is known as the $\mathbf{f}$-vector and traditionally written as $(f_1,\ldots,f_{d+1})$ where $f_i$ is the number of $(i-1)$-dimensional faces. A complete characterization of all possible $\mathbf{f}$-vectors of simplicial complexes is provided by the Kruskal-Katona theorem \cite{Frankl}. 

Notably, the DSC model generates compact simplicial complexes whose Hausdorff dimension is infinite. Indeed, the diameter of the $1$-skeleton of the simplicial complex generated at $n$-th step is 
\be
\mathcal{D}(n) = 2n - 1
\label{80}
,
\ee
while the number of vertices grows with $n$ faster than any power of $n$ (see the next subsection).

\subsection{Number of simplices}

The recursive nature of the DSC model leads to the recurrences 
\bea
\label{N0}
&&N_0(n+1)=N_0(n)+\sum_{j=0}^n N_j(n)  
,
\\[3pt]
\label{N1}
&&N_1(n+1)=N_1(n)+\sum_{j=0}^n (j+1)N_j(n)  
,
\\[3pt]
\label{N2}
&&N_2(n+1)=N_2(n)+\sum_{j=1}^n {j+1 \choose 2}N_j(n)
\eea
for the numbers of vertices, edges, and triangles. The general recurrence  reads 
\begin{equation}
\label{Nk}
N_d(n+1)=N_d(n)+\sum_{j=d-1}^n {j+1 \choose d}N_j(n)
.
\end{equation} 

Using these recurrences and performing pedestrian computations we find 
\be 
\label{N1-6}
\begin{split}
&\mathcal{N}(0)=[1] , \\
&\mathcal{N}(1)=[2,1], \\
&\mathcal{N}(2)=[5,5,1], \\
&\mathcal{N}(3)=[16,23,9,1], \\
&\mathcal{N}(4)=[65,116,65,14,1], \\
&\mathcal{N}(5)=[326, 669, 470,145,20,1], \\
&\mathcal{N}(6)=[1957, 4429, 3634, 1415, 280, 27, 1].
\end{split}
\ee
The arrays for the next three generations are presented in Appendix~\ref{sa1}. 

Interpreting $\{N_d(n)| 0\leq d,n<\infty\}$ as an infinite matrix and recalling that $N_d(n)=0$ if $d>n$, we conclude that we get an infinite lower triangular matrix.  The diagonal elements are $N_n(n)=1$ since $\mathcal{K}(n)$ has a single $n-$dimensional simplex. The elements $N_{n-1}(n)$ on the next co-diagonal constitute a quadratic polynomial in $n$ which is convenient to write in the form 
\begin{equation}
\label{N-diag-1}
N_{n-1}(n)=\binom{n}{1}\, \frac{n+3}{2}
\end{equation}
which will prove useful later. The elements $N_{n-2}(n)$ on the following co-diagonal constitute a quartic polynomial which we write again in the form 
\begin{equation}
\label{N-diag-2}
N_{n-2}(n) = \binom{n}{2}\,\frac{(n+1)(3n+14)}{12}
\end{equation}
involving a binomial coefficient. The elements $N_{n-3}(n)$ constitute a sextic polynomial 
\begin{equation}
\label{N-diag-3}
N_{n-3}(n) = {n \choose 3}\,\frac{n^3+8n^2+11n-4}{8}
.
\end{equation}

We guessed \eq{N-diag-1}--\eq{N-diag-2} using the data in \eqref{N1-6} and longer arrays presented in Appendix~\ref{sa1}, equation\eqref{N7-9}, and comparing with sequences in the on-line encyclopedia of integer sequences (OIES) \cite{Sloane}. The OIES does not contain a sequence appearing on the $(n-3,n)$ diagonal. Since $N_{n-p}(n)$ are polynomials of degree $2p$ when $p=0,1,2$, we guessed that the same holds for $p=3$ and obtained \eq{N-diag-3} by seeking a sextic polynomial reproducing the data in \eqref{N1-6} and \eqref{N7-9}. 

The sequence $N_0(n)$ giving the number of vertices appears in numerous combinatorial problems, see \cite{Flajolet}. For instance, it gives the number of permuatations of all subsets of a set with $n+1$ elements. The quantities $N_0(n)$ can be computed iterating the recurrence $N_0(n)=nN_0(n-1)+1$ starting from $N_0(0)=1$. Using the recurrence one derives an explicit formula 
\begin{equation}
\label{N0-sol}
N_0(n)=n!\left[e-\sum_{k>n}\frac{1}{k!}\right]
,
\end{equation}
which leads to a very accurate large $n$ asymptotic
\begin{equation}
\label{N0-asymp}
N_0(n)=e\, n! -\frac{1}{n}+\frac{1}{n^3}-\frac{1}{n^4}-\frac{2}{n^5}+\ldots 
.
\end{equation}

The sequence $N_1(n)$ appears in 
\begin{equation}
e^x\,\frac{\log(1-x)}{x-1}=\sum_{n\geq 1}\frac{N_1(n)}{n!}\,x^n
,
\label{N1n}
\end{equation}
from which one deduces the large $n$ asymptotics 
\begin{equation}
\label{N1-asymp}
N_1(n) \simeq e\, n!(\log n + \gamma)
,
\end{equation}
where $\gamma=0.5772\ldots$ is the Euler constant. 

The total number of simplices 
\begin{equation}
N(n)=\sum_{d=0}^n N_d(n) 
\label{210}
\end{equation}
can be expressed via the number of vertices: 
\begin{equation}
\label{N-all}
N(n) = N_0(n+1) - N_0(n) 
.
\end{equation}
This follows from the definition of the process and agrees with \eq{N0}. Combining \eq{N-all} and \eq{N0-sol} we obtain an explicit solution for the total number of simplices. The large $n$ asymptotic is
\begin{equation}
\label{N-asymp}
N(n)=e\, n\cdot n! +\frac{1}{n^2}-\frac{1}{n^3}-\frac{2}{n^4}+\ldots 
.
\end{equation}
The exponential generating function of the series $N(n)$ has a neat form 
\be
\sum_{n\geq 0} N(n)\,\frac{x^n}{n!} = \frac{e^x}{(1-x)^2}
.
\label{235}
\ee

The exponential generating functions associated with $N_d(n)$ are also simple for all $d\geq 0$. One extracts
\begin{equation}
\label{Ndn-0}
N_d(n) =\left. \frac{(-1)^d}{d!} \frac{d^n}{d x^n} \left\{e^x \frac{[\ln(1 - x)]^d}{1 - x}\right\} \right |_{x=0} 
\end{equation}
from these generating functions, or equivalently the representation of $N_d(n)$ via contour integrals
\begin{equation}
\label{Ndn-int}
N_d(n) = e n! \,
\frac{1}{2\pi {\rm i}} 
\frac{(-1)^d}{d!} \oint\limits_C \frac{dx}{x^{n+1}} \, \frac{[\ln(1 - x)]^d}{(1 - x)e^{1-x}}\,.
\end{equation}
Here $C$ is a positively oriented simple closed contour around the origin which does not enclose the singularities of the function $[\ln(1 - x)]^d/[(1 - x)e^{1-x}]$; for instance, $C$ can be a circle around the origin of radius $R<1$. 

The representation \eq{Ndn-0}, \eq{Ndn-int} is valid for all $0\leq d\leq n$. Specializing equation~\eq{Ndn-0} to $d=0$ and  $d=1$ we recover equations~\eq{N0-sol} and \eq{N1n}. We also used equation~\eq{Ndn-0} to confirm all numbers appearing in equations \eqref{N1-6} and \eqref{N7-9}. When $d=O(1)$ is fixed and $n\gg 1$, we deduce the leading asymptotic behavior
\begin{equation}
N_d(n) \cong e n! \frac{(\ln n)^d}{d!}
\label{Ndn-asymp}
\end{equation}
from equation~\eq{Ndn-int}, see \cite{evgrafov2020asymptotic}.  

Using equation~\eq{Ndn-0}, one can also derive simple exact formulae  for the number of simplices of small co-dimension: $m=n-d=O(1)$. When $m=1,2,3$, i.e., $d=n-1, n-2, n-3$,  we recover equations~\eq{N-diag-1}--\eq{N-diag-3} which we previously guessed. One can derive similar formulae for larger $m$ relying on the pattern evident from equations~\eq{N-diag-1}--\eq{N-diag-3} and generally correct, namely
\begin{equation}
\label{N-m}
N_{n-m}(n) = \binom{n}{m}\,P_m(n)
,
\end{equation}
where $P_m(n)$ are the polynomials of $n$ of degree $m$. We already know $P_m(n)$ with $m=0,1,2,3$: $P_0(n)=1$ and $P_{1,2,3}(n)$ are given by equations~\eq{N-diag-1}, \eq{N-diag-2}, \eq{N-diag-3}, respectively. The next two polynomials are
\bea
&&
\hspace{-18pt}
P_4(n) = \frac{15n^4 + 150n^3+245n^2-378n-248}{240}\, 
, 
\\[3pt]
&&
\hspace{-18pt}
P_5(n) = \frac{3 n^5 + 35 n^4 + 55 n^3  - 243 n^2 - 202 n + 256}{96}\,
. 
\label{P4-5}
\eea
The leading asymptotic behavior is
\begin{equation}
\label{N-m-asymp}
N_{n-m}(n) \simeq  \frac{n^{2m}}{2^m\,m!}
\end{equation}
when $m=O(1)$ is fixed and $n\gg 1$.

\subsection{Upper degrees}

The quantities $N_d(n)$ describe the global properties of emerging simplicial complexes. We haven't yet examined local properties. The basic local characteristic is the upper degree $k^{(d)}_\sigma$ of a $d$-simplex $\sigma$. The mean value $\avg{k^{(d)}}$ of the upper degree of a $d$-simplex in a simplicial complex can be expressed in terms of $N_d$ and $N_{d+1}$,  
\be
\avg{k^{(d)}} = \frac{(d+2)N_{d+1}}{N_d}
.
\label{310}
\ee
Let us first introduce the set of lists $K^{(d)}$ of distinct $d$-degrees in a simplicial complex, where each $K^{(d)}$ contains the distinct $d$-degrees ordered in ascending order and listed without repetitions, regardless of how many simplices 
%%share 
have the same $d$-degree, 
\be
K^{(d)} \equiv \left[k^{(d)}_1, k^{(d)}_2, \ldots,k^{(d)}_\alpha, \ldots, k^{(d)}_{\alpha^{(d)}_{\text{max}}-1}, k^{(d)}_{\alpha^{(d)}_{\text{max}}}\right]
, 
\label{320}
\ee
where $\alpha^{(d)}_{\text{max}}$ is the number of distinct $d$-degrees in a simplicial complex. 

To uncover the rules governing the evolution of $K^{(d)}(n)$, we employ a simple script to generate the DSC simplicial complexes up to $n=8$. We measured the upper degrees of their simplices at each step, and inferred the patterns that describe the evolution of the distinct upper degrees. Here we show $K^{(d)}(n)$ for $n \leq 5$,
\bea
&&
K^{(0)}(1) = [1]
,
\nonumber
\\[3pt]
&&
K^{(0)}(2) = [1,2,3]
,
\nonumber
\\[3pt]
&&
K^{(1)}(2) = [0,1]
,
\nonumber
\\[3pt]
&&
K^{(0)}(3) = [1,2,3,6,8]
,
\nonumber
\\[3pt]
&&
K^{(1)}(3) = [0,1,2,3]
,
\nonumber
\\[3pt]
&&
K^{(2)}(3) = [0,1]
,
\nonumber
\\[5pt]
&&
K^{(0)}(4) = [1,2,3,4,6,8,11,19,24]
,
\nonumber
\\[3pt]
&&
K^{(1)}(4) = [0,1,2,3,6,8]
,
\nonumber
\\[3pt]
&&
K^{(2)}(4) = [0,1,2,3]
,
\nonumber
\\[3pt]
&&
K^{(3)}(4) = [0,1]
,
\nonumber
\\[5pt]
&&
K^{(0)}(5) = [1,2,3,4,5,6,8,11,19,20,24,46,73,89]
,
\nonumber
\\[3pt]
&&
K^{(1)}(5) = [0,1,2,3,4,6,8,11,19,24]
,
\nonumber
\\[3pt]
&&
K^{(2)}(5) = [0,1,2,3,6,8]
,
\nonumber
\\[3pt]
&&
K^{(3)}(5) = [0,1,2,3]
,
\nonumber
\\[3pt]
&&
K^{(4)}(5) = [0,1].  
\label{330}
\eea

Numerous coinciding sequences in \eqref{330} hint that the quantities $K^{(n-p)}(n)$ with fixed $p$ are stationary, viz. independent on $n$ when $n\geq p+1$. One can extract such stationary sequences with $p\leq 4$ from equations \eqref{330}. We list such sequences for $p\leq 8$ in Appendix~\ref{sa1}. The knowledge of such infinite sequences allows one to fix all $K^{(d)}(n)$ for $n \leq 8$.  Another useful general relation  
\begin{equation}
\label{K-stat}
K^{(d)}(n) = 0 \cup K^{(0)}(n-d)
\end{equation}
is conjecturally valid for all $n \geq d + 1$ and $d \geq 1$. 

Let us look more carefully at the first and second largest $d$-degrees. We shortly write $M^{(d)}(n)\equiv k^{(d)}_{\alpha^{(d)}_{\text{max}}}(n)$ and $\widetilde{M}^{(d)}(n) \equiv k^{(d)}_{\alpha^{(d)}_{\text{max}}-1}(n)$. We have 
\be
M^{(d)}(n) = M^{(0)}(n-d)
,
\label{350}
\ee
for $n \geq d + 1$ and $d \geq 1$, and 
\be
\widetilde{M}^{(d)}(n) = \widetilde{M}^{(0)}(n-d)
\label{360}
\ee
for $n \geq d + 2$ and $d \geq 1$. 

Clearly, 
\be
N_0(n) = M^{(0)}(n+1) - M^{(0)}(n) 
%%.
\label{365}
\ee
and hence 
\be
N(n) = M^{(0)}(n+2) - 2M^{(0)}(n+1) + M^{(0)}(n)
, 
\label{370}
\ee
where we used \eq{N-all}. Equations~\eq{N1-6}, \eqref{N7-9} and \eq{365} provide the following series $M^{(0)}(n)$: 
\begin{equation}
\label{Mn-small}
1,  3,  8, 24, 89, 415, 2372, 16072, \ldots
.
\end{equation}
This series can be observed in equations~\eq{330} and \eq{a10}.  
The numbers $M^{(0)}(n)$ are known as logarithmic numbers since they appear in the expansion 
\begin{equation}
\label{Mn-exp}
e^x \ln(1 - x) = -\sum_{n\geq 1} \frac{M^{(0)}(n)}{n!}\,x^n
.
\end{equation}
One can express $M^{(0)}(n)$ via contour integrals, 
\begin{equation}
\label{Mn-int}
M^{(0)}(n) = - e n! \,\frac{1}{2\pi \text{i}}  \oint\limits_C \frac{dx}{x^{n+1}} \, e^{-(1-x)} \ln(1 - x) 
,
\end{equation} 
where $C$ is a positively oriented simple closed contour around the origin which does not enclose the singularities of the function $e^{-(1-x)}\ln(1 - x)$. 
Using \eq{Mn-int} we deduce the large $n$ asymptotics  
\begin{equation}
M^{(0)}(n) \cong e n! /n
.
\label{Mn-asymp}
\end{equation} 
Furthermore, inspecting the lists $K^{(0)}(n)$ for $n \leq 8$, equations~\eq{330} and \eq{a10}, we notice that the difference between the largest degree $M^{(0)}(n)$ and the second largest degree $\widetilde{M}^{(0)}(n)$ satisfies 
\be
M^{(0)}(n) - \widetilde{M}^{(0)}(n) = \left. \frac{d^{n-2}}{dx^{n-2}} \left( \frac{e^x}{1-x} \right) \right|_{x=0}
.
\label{375}
\ee
We guess that equation~\eq{375} is valid for an arbitrary $n$, which leads to the asymptotics: 
\be
M^{(0)}(n) - \widetilde{M}^{(0)}(n) \cong \frac{e\, n!}{(n-1)n}
.
\label{380}
\ee
Once we derived equations~\eq{Mn-asymp} and \eq{380}, the large-$n$ asymptotics of $M^{(d)}(n)$ and $M^{(d)}(n) - \widetilde{M}^{(d)}(n)$ became immediately apparent from equations~\eq{350} and \eq{360}, 
\bea
\hspace{54.5pt} 
M^{(d)}(n) &\cong& e\, (n-d-1)! 
,
\label{382}
\\[3pt] 
M^{(d)}(n) - \widetilde{M}^{(d)}(n) &\cong& e\, (n-d-2)!
.
\label{383}
\eea

The complete upper-degree sequences account for the number of occurrences of distinct upper degrees in a simplicial complex, that is, the degeneracies of upper degrees. Let $D^{(d)}(k^{(d)},n)$ be the number of $d$-simplices of $d$-degree $k^{(d)}$ in the simplicial complex of generation $n$. It is convenient to introduce the arrays of upper degrees along with their degeneracies: 
\bea
&&\hspace{-20pt}
Q^{(d)}(n) \equiv \left[\hspace{1.5pt}[k_1^{(d)}, D^{(d)}(k_1^{(d)},n)] , \ldots , 
\right. 
\nonumber 
\\[3pt]
&&\hspace{55pt}
\left. 
\left[k^{(d)}_{\alpha^{(d)}_{\text{max}}(n)}, D^{(d)}\left(k^{(d)}_{\alpha^{(d)}_{\text{max}}(n)},n\right)\right]
\hspace{1.5pt}\right]
.
\label{385}
\eea
Figure~\ref{fig:2} provides the arrays $Q^{(d)}(n)$ for $n=1,2,3$:    
\bea
&&
\hspace{-72pt}
Q^{(0)}(1) = [\hspace{1.5pt}[1,2]\hspace{1.5pt}]
,
\nonumber
\\[3pt]
&&
\hspace{-72pt}
Q^{(0)}(2) = [\hspace{1.5pt}[1,2],[2,1],[3,2]\hspace{1.5pt}]
,
\nonumber
\\[3pt]
&&
\hspace{-72pt}
Q^{(1)}(2) = [\hspace{1.5pt}[0,2],[1,3]\hspace{1.5pt}]
,
\nonumber
\\[3pt]
&&
\hspace{-72pt}
Q^{(0)}(3) = [\hspace{1.5pt}[1,5],[2,5],[3,3],[6,1],[8,2]\hspace{1.5pt}]
,
\nonumber
\\[3pt]
&&
\hspace{-72pt}
Q^{(1)}(3) = [\hspace{1.5pt}[0,5],[1,12],[2,3],[3,3]\hspace{1.5pt}]
,
\nonumber
\\[3pt]
&&
\hspace{-72pt}
Q^{(2)}(3) = [\hspace{1.5pt}[0,5],[1,4]\hspace{1.5pt}]
,
\label{390}
\eea
We used Mathematica to determine the arrays $Q^{(d)}(n)$ with $4\leq n \leq 8$, see Appendix~\ref{sa1}. 

The number $D^{(0)}(1,n)$ of vertices of degree $1$ is 
\begin{equation}
D^{(0)}(1,n) = D^{(1)}(0,n) = N_0(n-1)
%%.
\label{392}
\end{equation}
and the number $D^{(0)}(2,n)$ of vertices of degree $2$ is 
\begin{equation}
D^{(0)}(2,n) = N_1(n-1)
.
\label{394}
\end{equation}
These results are valid for $n \geq 2$. The sum rules
\bea
&&
\sum_{j=1}^{M^{(0)}(n)} D^{(0)}(j,n)=N_0(n)
, 
\label{D-sums}
\\[3pt]
&&
\sum_{j=1}^{M^{(0)}(n)} j D^{(0)}(j,n)=2N_1(n) 
,
\label{3945}
\eea 
lead to 
\begin{equation}
\sum_{j=2}^{M^{(0)}(n)} (j-1) D^{(0)}(j,n)=2N_1(n)-N_0(n)
,
\label{395}
\end{equation}
which in conjunction with $D^{(0)}_{M^{(0)}(n)}(n)=2$ (see figure~\ref{fig:2}) and the asymptotic behaviors \eq{N0-asymp} and \eq{N1-asymp} leads to the bound $M^{(0)}(n) < e\, n!\,\ln n$ asserting that the maximal degree is smaller than the total number of edges $N_1(n)$. This upper bond is not tight as is evident from the asymptotics of $M^{(0)}(n)$, equation~\eq{Mn-asymp}. 

In general, for $n \geq 2$ and $1 \leq d \leq n-1$, the number $D^{(d)}(0,n)$ of $d$-simplices of $d$-degree $0$ is 
\be
D^{(d)}(0,n) = N_{d-1}(n-1) 
.
\label{396}
\ee
For $n \geq 2$, the numbers $D^{(n-1)}(0,n)$ and $D^{(n-1)}(1,n)$  of $(d=n-1)$-simplices of $d$-degrees $0$ and $1$ are  
\begin{equation}
\label{Dn1n}
\begin{split}
D^{(n-1)}(0,n) &= \frac{(n - 1)(n + 2)}{2}\,, \\
D^{(n-1)}(1,n) &= n + 1.
\end{split}
\end{equation}
For $n \geq 3$, the numbers $D^{(n-2)}(0,n)$,  $D^{(n-2)}(1,n)$,  $D^{(n-2)}(2,n)$, and $D^{(n-2)}(3,n)$  of $(d=n-2)$-simplices of $d$-degrees $0$, $1$, $2$, and $3$, respectively, are  
\bea
D^{(n-2)}(0,n) &=& \frac{(n - 2)(n - 1)n(3n + 11)}{24}\,, 
\nonumber
\\[3pt]
D^{(n-2)}(1,n) &=& \frac{n(n^2 + n - 4)}{2}\,, 
\nonumber
\\[3pt]
D^{(n-2)}(2,n) &=& \frac{(n - 1)n}{2}\,, 
\nonumber
\\[3pt]
D^{(n-2)}(3,n) &=& n.
\label{Dn2n}
\eea
Furthermore, 
\begin{equation}
D^{(d)}\!\left(M^{(d)}(n), n\right) = d + 2
\end{equation}
for $n\geq 1$ and
\begin{equation}
\label{399}
D^{(d)}\!\left(\widetilde{M}^{(d)}(n)
,n\right) = \frac{(d+1)(d+2)}{2}
\end{equation}
for $n\geq d+2$. 
We guessed \eq{Dn1n}--\eq{399} from $Q^{(d)}(n)$ with $n \leq 8$ obtained with the help of Mathematica, see \eqref{Qn1}--\eq{a20}; we conjecture that \eq{Dn1n}--\eq{399} hold for all $n$. 

The fraction of $d$-simplices with a given $d$-degree is 
\begin{equation}
p^{(d)}(k^{(d)},n) = \frac{D^{(d)}(k^{(d)},n)}{N_d(n)}
.
\label{400}
\end{equation}
We shortly refer to $p^{(d)}(k^{(d)})$ as the upper-degree distribution, or simply the $d$-degree distribution. Specifically, $p^{(0)}(k^{(0)})$ is the degree distribution of the $1$-skeleton of a simplicial complex. In particular, using \eq{392}, we compute the fraction of vertices of the smallest degree 
\begin{equation}
p^{(0)}(1,n)=\frac{N_0(n-1)}{N_0(n)}
.
\label{401}
\end{equation} 
Hence $p^{(0)}(1,n) \simeq n^{-1}$ as $n\gg 1$, i.e., it asymptotically vanishes. 

The fraction $p^{(d)}(k^{(d)})$ cannot be directly compared to the degree distributions of random simplicial complexes, as deterministic complexes often exhibit non-uniform gaps between degree values. In our simplicial complex, these gaps become significantly wider as the degrees increase. To relate the exponent of this discrete upper-degree distribution to the standard $\gamma$ exponent of a continuous degree distribution for random scale-free networks, we use the cumulative distribution  
\begin{equation}
p^{(d)}_{\text{cum}}(k^{(d)},n) \equiv \sum_{q \geq k} p^{(d)}(q^{(d)},n)
, 
\label{402}
\end{equation}
for which exponent $\gamma$ is defined as $p_{\text{cum}}(k) \sim k^{-\gamma + 1}$. 

%%%%%%%%%%%%%%%%%%%%%%%%%%%%%%%%%%%%%%%%%%%%%%%
%%%%%%%%%%%%%%%%%%%%%%%%%%%%%%%%%%%%%%%%%%%%%%%

\begin{figure}[t]
\begin{center}
\includegraphics[width=0.47\textwidth]{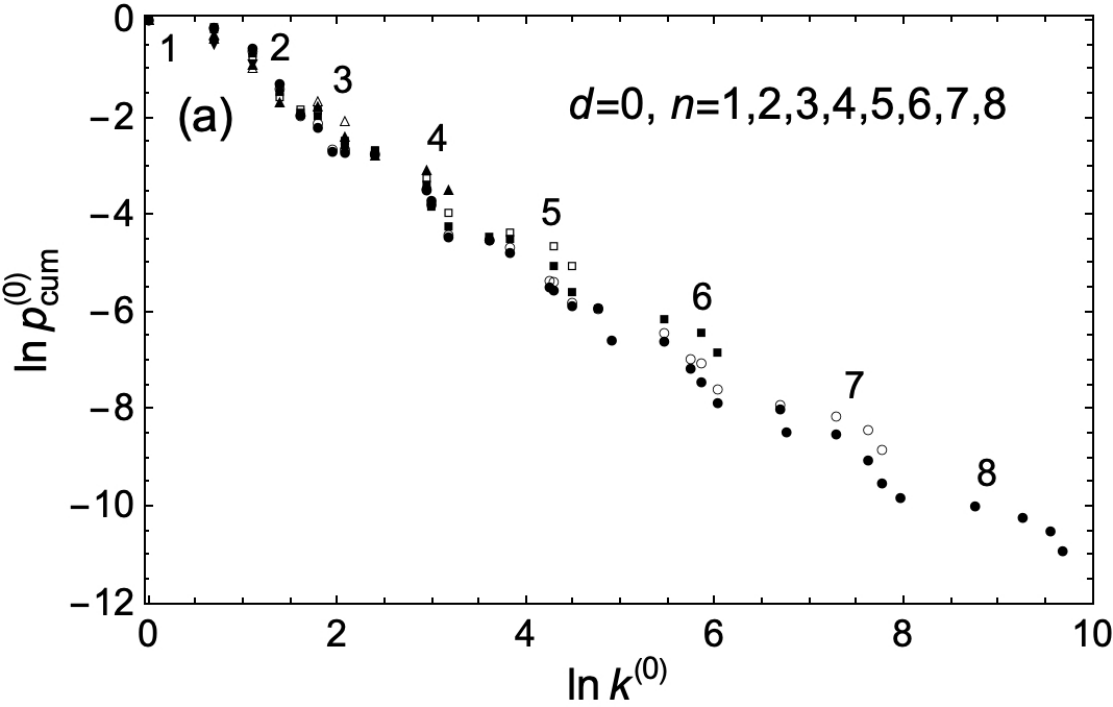}
\\[3pt]
\includegraphics[width=0.47\textwidth]{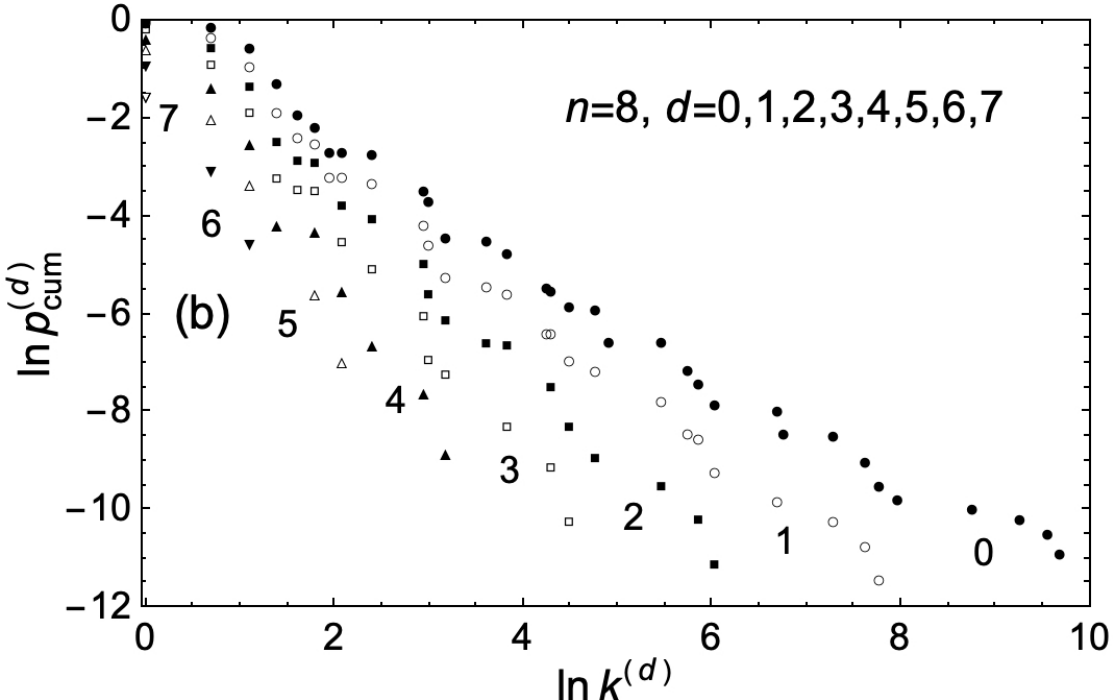}
\caption{The log--log plots of the cumulative upper-degree distributions $p^{(d)}_{\text{cum}}(k^{(d)},n)$ for simplicial complexes generated by the DSC model at various values of $n$ and $d$. 
(a)~The log--log plot of the cumulative degree distributions for the $1$-skeletons of the DSC simplicial complexes at different generations $n$. 
Markers represent different values of $n$: 
{\Large \raisebox{-1pt}{$\bullet$}}, {\Large \raisebox{-1pt}{$\circ$}}, 
%%$\Circle[f]$, $\Circle$, 
%%%%%%%%%%%%%%%%%%%%%%%%%%%%%%%%%%%%%%$\fullsquare$, $\opensquare$, 
{\scriptsize $\blacksquare$}, {\scriptsize $\square$}, 
$\blacktriangle$, $\protect\rotatebox[origin=c]{90}{$\rhd$}$, $\blacktriangledown$, and $\triangledown$ 
correspond to $n=8$ to $1$, respectively. 
(b)~The log--log plots of the cumulative $d$-degree distributions of the DSC simplicial complexes at generation $n=8$ for various $d$. 
Markers represent different values of $d$: 
{\Large \raisebox{-1pt}{$\bullet$}}, {\Large \raisebox{-1pt}{$\circ$}}, 
%%$\Circle[f]$, $\Circle$, 
%%%%%%%%%%%%%%%%%%%%%%%%%%%%%%%%%%%%%%$\fullsquare$,  $\opensquare$, 
{\scriptsize $\blacksquare$}, {\scriptsize $\square$}, 
$\blacktriangle$, $\protect\rotatebox[origin=c]{90}{$\rhd$}$, $\blacktriangledown$, and $\triangledown$ 
correspond to $d=0$ to  
$7$, respectively. 
} 
\label{f2}
\end{center}
\end{figure}

%%%%%%%%%%%%%%%%%%%%%%%%%%%%%%%%%%%%%%%%%%%%%%%
%%%%%%%%%%%%%%%%%%%%%%%%%%%%%%%%%%%%%%%%%%%%%%%

\begin{figure}[t]
\begin{center}
\includegraphics[width=0.47\textwidth]{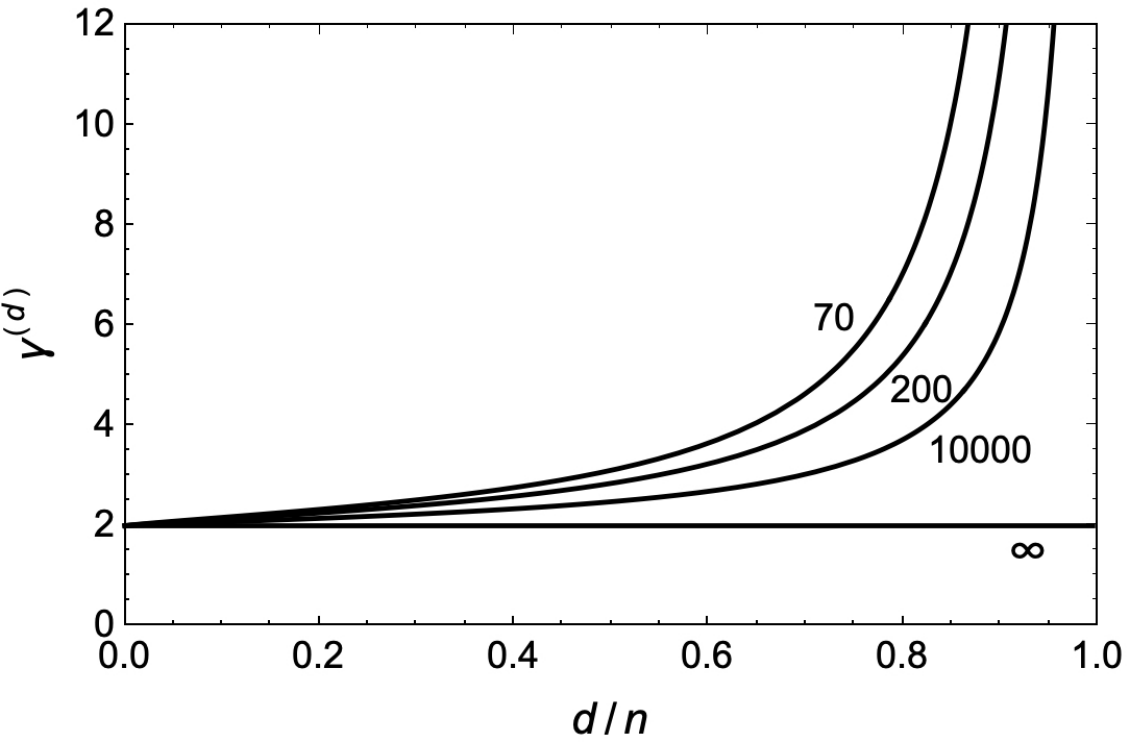}
\caption{Estimated values of the exponent $\gamma^{(d)}$ using the expression in first line of equation~\protect\eq{405}, shown for generations $n=70$, $200$, $10000$, and in the limit $n \to \infty$.
} 
\label{f3}
\end{center}
\end{figure}

We obtained numerically the cumulative upper-degree distributions $p^{(d)}_{\text{cum}}(k^{(d)},n)$ in the simplicial complexes generated by the DSC model for $n \leq 8$, see figure~\ref{f2}. Figure~\ref{f2}(a) shows that $p^{(0)}_{\text{cum}}(k^{(0)},n)$ rapidly converges to a stationary cumulative distribution as $n$ increases. Figure~\ref{f2}(b) suggests that the cumulative $d$-degree distributions $p^{(d)}_{\text{cum}}(k^{(d)},n)$ follow a power-law behavior. Specifically, the slope of the curve for $d=0$ and $n=8$ is close to $1$, indicating exponent $\gamma^{(0)}$ near $2$. 
To estimate the convergence rate of the $d$-degree distributions and their exponents $\gamma^{(d)}$, we compare the logarithms of $p^{(d)}_{\text{cum}}(k^{(d)},n)$ at the endpoints $k^{(d)}=1$ and $k^{(d)}=M^{(d)(n)}$,    
\be
\gamma^{(d)}(n) \approx 1 - \frac{\ln p^{(d)}_{\text{cum}}(M^{(d)}(n),n) - \ln p^{(d)}_{\text{cum}}(1,n)}{\ln M^{(d)}(n) - \ln 1 }
.
\label{403}
\ee
For $d=0$ and sufficiently large $n$, this gives  
\be
\gamma^{(0)}(n) \approx 1 + \frac{\ln(en!)}{\ln(en!) - \ln n} \cong 2 +  \frac{1}{n} + \frac{1}{n\,\ln n}\,
,                                                                      
\label{404}
\ee
where we used equations~\eq{N0-asymp} and \eq{Mn-asymp}. The $1/n$ correction indicates that that the cumulative distribution quickly approaches its limiting stationary form. For $d \geq 1$, the situation is notably different. For instance, if $n\to\infty$ and $d\to\infty$ with ratio $\delta\equiv d/n$ kept constant, we use equations \eq{Ndn-asymp}, \eq{382}, and \eq{396} to get  
\bea  
\gamma^{(d)}(n) &\approx& 1 + \frac{\ln n + \delta(1 - \delta)^{-1} [\ln\ln n - \ln \delta] - 1}{\ln n + \ln(1 - \delta) - 1} 
\nonumber
\\[3pt]
&\cong& 2 +  \frac{\delta}{1-\delta}\frac{\ln\ln n}{\ln n} 
.                                                                     
\label{405}
\eea
The correction term $\ln\ln n/\ln n$ implies that the cumulative distribution converges to $\gamma^{(d)}(\infty)=2$ at an anomalously slow pace. Figure~\ref{f3} clearly illustrates how slow this convergence is.  For $n=70$, $200$, and $10000$, used in the figure, the number of vertices in this simplicial complex reaches staggering values---on the order of $10^{100}$, $10^{375}$, and $10^{35\,660}$, respectively. For instance, these numbers greatly exceed the number of atoms in the observable Universe estimated to be between $10^{78}$ or $10^{82}$.

%%%%%%%%%%%%%%%%%%%%%%%%%%%%%%%%%%%%%%%%%%%%%%%
%%%%%%%%%%%%%%%%%%%%%%%%%%%%%%%%%%%%%%%%%%%%%%%

%%%%%%%%%%%%%%%%%%%%%%%%%%%%%%%%%%%%%%%%%%%%%%%
%%%%%%%%%%%%%%%%%%%%%%%%%%%%%%%%%%%%%%%%%%%%%%%

\subsection{Hodge Laplacian spectra in the DSC model}

We begin by introducing necessary definitions. Let $B^{(d)}$ be the $N_{d-1} \times N_d$ incidence matrix of an oriented simplicial complex. The element $(B^{(d)})_{\sigma_{d-1} \sigma_d}$ of the incidence matrix is non-zero $(\pm 1)$ if simplex $\sigma_{d-1}$ is a face of simplex $\sigma_d$. The sign of this element is positive when the orientations of $\sigma_d$ and $\sigma_{d-1}$ fit to each other and negative when the orientations are opposite. Thus
\be
(B^{(d)})_{\sigma_{d-1} \sigma_d} = (-1)^p 
\label{410}
\ee
for the simplices 
\bea
\sigma_d &=& [i_0, i_1, \ldots, i_d] 
,
\nonumber
\\[3pt]
\sigma_{d-1} &=& [i_0, i_1, \ldots, i_{p-1}, i_{p+1}, \ldots, i_d] 
;
\label{415}
\eea
otherwise, $(B^{(d)})_{\sigma_{d-1} \sigma_d} = 0$. 

The $d$-th Hodge Laplacian is defined as 
\be
L^{(d)} = [B^{(d)}]^T B^{(d)} + B^{(d+1)} [B^{(d+1)}]^T
.
\label{420}
\ee
It is an $N_d \times N_d$ matrix. The $d$-th Hodge Laplacian is associated with diffusion and flows in the array of $d$-simplices in a simplicial complex, connected through strong overlaps. By definition, two $d$-simplices strongly overlap when they share a $(d {-} 1)$-face. In particular, the $0$-th Hodge Laplacian is the ordinary Laplacian $L$ of the $1$-skeleton graph of a given simplicial complex,
\be
L \equiv L^{(0)} = B^{(1)} [B^{(1)}]^T
.
\label{425}
\ee
This Laplacian is associated with diffusion and flows on the $1$-skeleton graph. See \cite{Ginestra21,lim2020hodge, schaub2020random} for reviews of Hodge Laplacians and diffusion on simplicial complexes. Instructive examples of Hodge Laplacians are given in \cite{dorogovtsev2022the}, section~$13.2.2$. 

It is convenient to introduce the cumulative spectral density of the $d$-th Hodge Laplacian, 
\be
g_{\text{cum}}^{(d)}(\lambda) \equiv \frac{\mathfrak{N}_{\leq\phantom{\!\!\!\!\! y_y\!\!}}^{(d)}(\lambda)}{N_d}\,,
\label{430}
\ee
where $\mathfrak{N}_\leq^{(d)}(\lambda)$ is the number of eigenvalues smaller or equal than $\lambda$, since we focus on the small-$\lambda$ region of these spectra. Similarly to the spectral dimension $d_{\text{s}}$ of the Laplacian of a graph, defined by the power-law asymptotics 
$g_{\text{cum}}(\lambda) \sim \lambda^{d_s/2}$ (if it exists), one can introduce the higher-order spectral dimensions $d_{s}^{(d)}$ defined by the small-$\lambda$ asymptotics 
\be
g_{\text{cum}}^{(d)}(\lambda) \sim \lambda^{d_s^{(d)}\!/2}
.  
\label{435}
\ee
In particular, $d_{s}^{(0)} = d_{s}$. 

%%%%%%%%%%%%%%%%%%%%%%%%%%%%%%%%%%%%%%%%%%%%%%%
%%%%%%%%%%%%%%%%%%%%%%%%%%%%%%%%%%%%%%%%%%%%%%%

\begin{figure}[t]
\begin{center}
\includegraphics[width=0.462\textwidth]{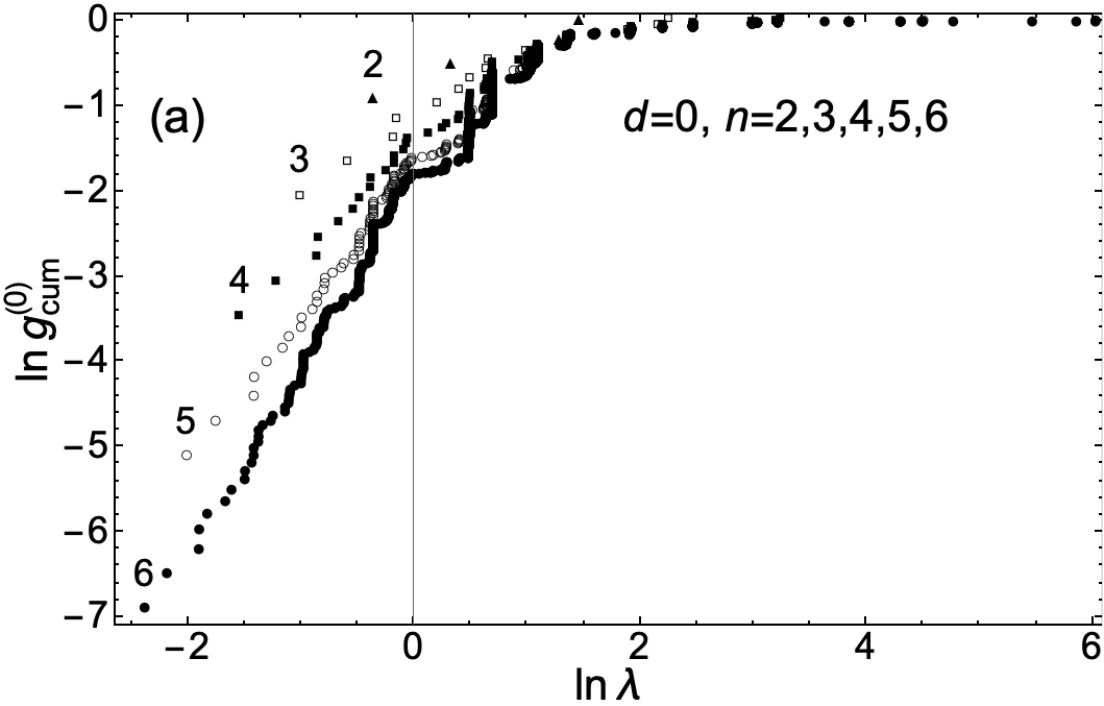}
\\[3pt]
\includegraphics[width=0.488\textwidth]{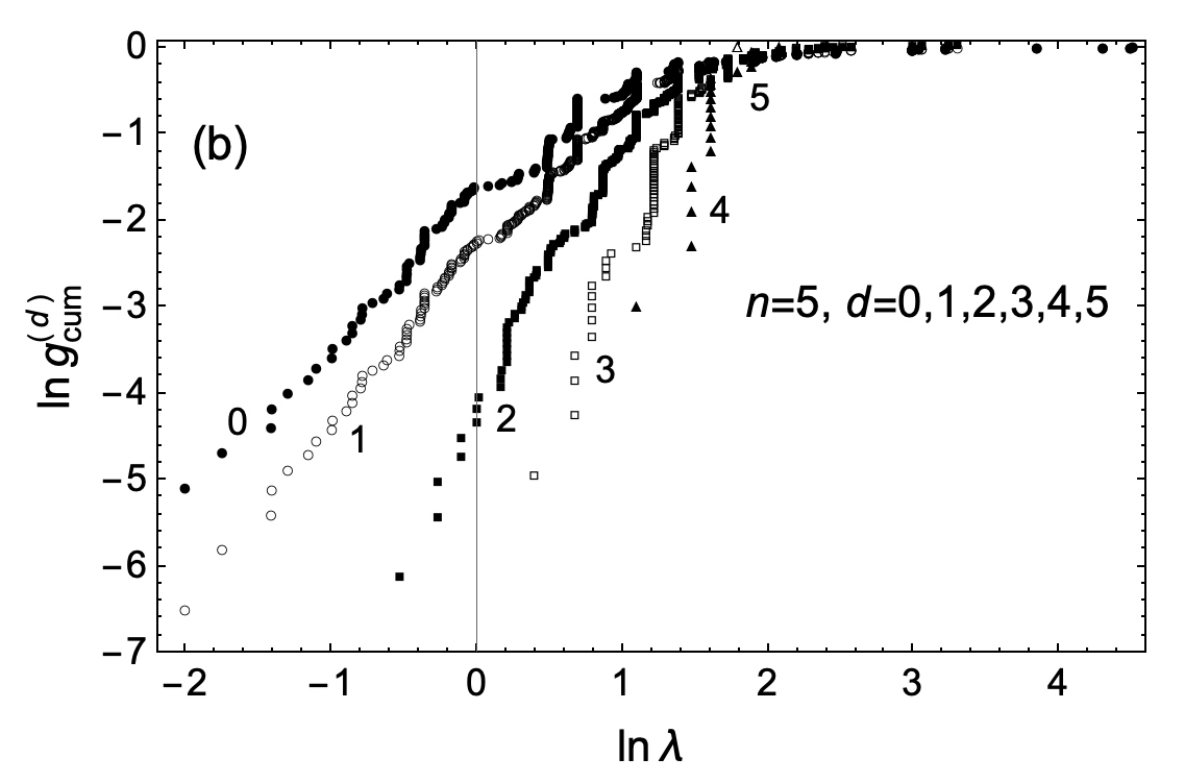}
\caption{The log--log plots of the cumulative spectral densities $g^{(d)}_{\text{cum}}(\lambda,n)$ of Hodge Laplacians for simplicial complexes generated by the DSC model at various values of $n$ and $d$.  (a)~The log--log plot of the cumulative spectral densities $g^{(0)}_{\text{cum}}(\lambda,n)$ of Laplacians for the $1$-skeletons of the DSC simplicial complexes at different generations $n$. Markers represent different values of $n$: 
{\Large \raisebox{-1pt}{$\bullet$}}, {\Large \raisebox{-1pt}{$\circ$}},  
%%$\Circle[f]$, $\Circle$, 
%%%%%%%%%%%%%%%%%%%%%%%%%%%%%%%%%%%$\fullsquare$,  $\opensquare$, 
{\scriptsize $\blacksquare$}, {\scriptsize $\square$}, 
and $\blacktriangle$ correspond to $n=6$ to $2$, respectively. The labels on the plot mark the smallest nonzero eigenvalues $\lambda_2(n)$ from the Laplacian spectra. (b)~The log--log plots of the cumulative spectral densities $g^{(d)}_{\text{cum}}(\lambda,n)$ of Hodge Laplacians for the DSC simplicial complexes at generation $n=5$ for various $d$. Markers represent different values of $d$: 
{\Large \raisebox{-1pt}{$\bullet$}}, {\Large \raisebox{-1pt}{$\circ$}}, 
%%$\Circle[f]$, $\Circle$, 
%%%%%%%%%%%%%%%%%%%%%%%%%%%%%%%%%%%$\fullsquare$,  $\opensquare$, 
{\scriptsize $\blacksquare$}, {\scriptsize $\square$}, 
$\blacktriangle$, and $\protect\rotatebox[origin=c]{90}{$\rhd$}$, correspond to $d=0$ to  
$5$, respectively. 
} 
\label{f4}
\end{center}
\end{figure}

%%%%%%%%%%%%%%%%%%%%%%%%%%%%%%%%%%%%%%%%%%%%%%%
%%%%%%%%%%%%%%%%%%%%%%%%%%%%%%%%%%%%%%%%%%%%%%%

Using equations~\eq{410}--\eq{435}, we numerically computed the Hodge Laplacian spectra of the DSC model for the first few generations, see figure~\ref{f4}. Note that the matching eigenvalues in the spectra of adjacent Hodge Laplacians $L{(d)}$ and $L{(d+1)}$ in figure~\ref{f4}(b) are a direct consequence of the Hodge decomposition theorem \cite{hodge1934dirichlet}. In contrast, the coincident eigenvalues in the spectra of $L^{(d)}$, $L^{(d+2)}$, etc., and the degenerate eigenvalues---visible as jumps in the cumulative spectral densities---arise from the deterministic nature of the process generating these simplicial complexes. A single zero eigenvalue appears only in the spectrum of $L = L^{(0)}$, reflecting the simple topology of our simplicial complex: it has no $d$-dimensional holes ($d$-holes) for $d > 1$, and it consists of a single connected component---$0$-hole. (The degeneracy of a zero eigenvalue in the spectrum of the $d$-th Hodge Laplacian is equal to the Betti number $b_d$, which counts the number of $d$-holes in the simplicial complex.)

\begin{figure}[t]
\begin{center}
\includegraphics[width=0.46\textwidth]{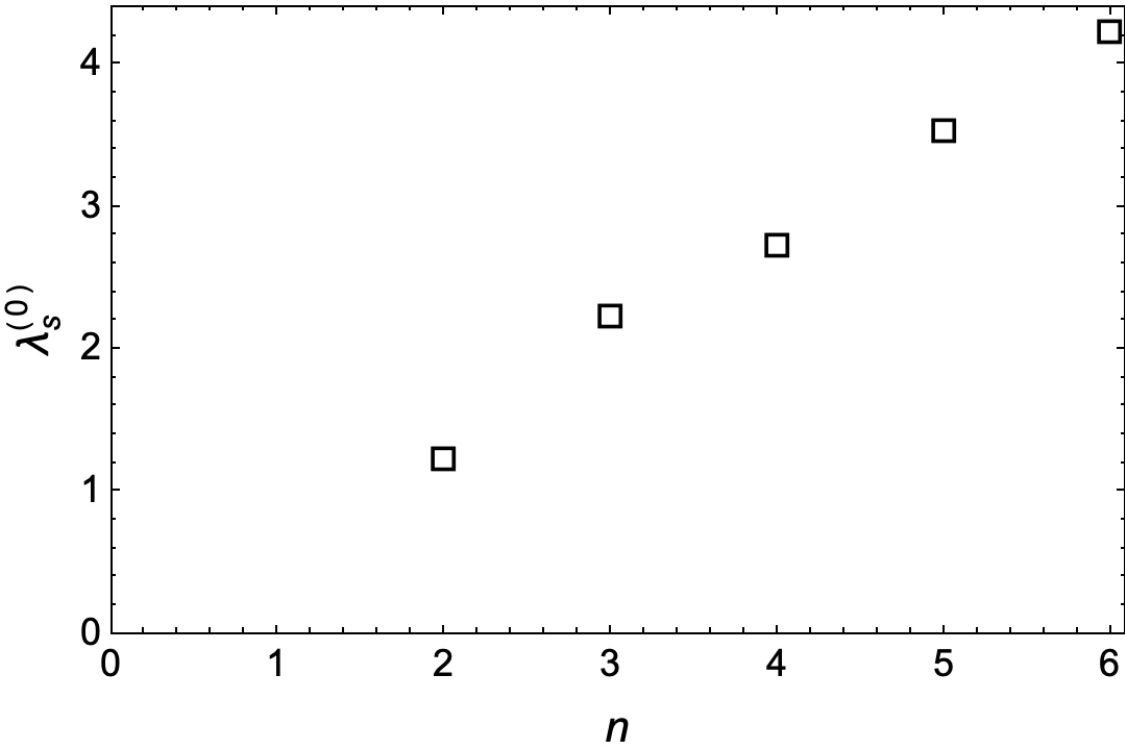}
\caption{Spectral dimension $d_{s}^{(0)}(n)$ in the DSC model estimated from the slopes of the curves in figure~\protect\ref{f4}(a) at low $\lambda$. 
} 
\label{f5}
\end{center}
\end{figure}

\begin{figure}[t]
\begin{center}
\includegraphics[width=0.46\textwidth]{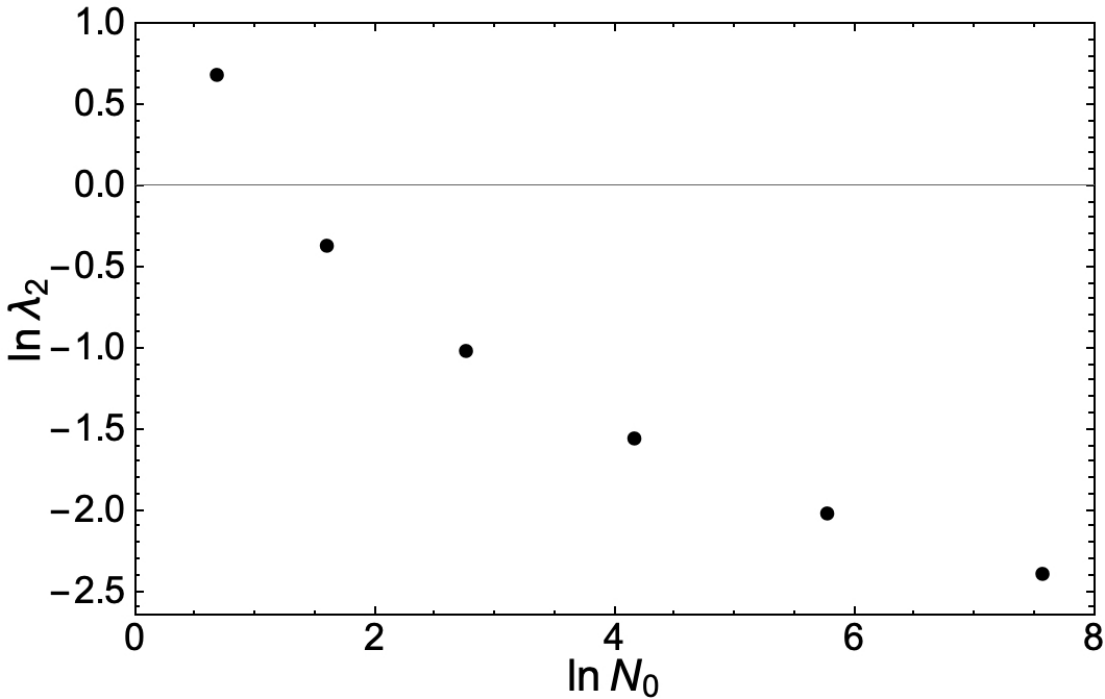}
\caption{Log--log plot of the smallest nonzero Laplacian eigenvalue $\lambda_2$ versus the number of vertices $N_0$ in simplicial complexes generated by the DSC model, where $\lambda_2(n)$ corresponds to the labelled points in figure~\protect\ref{f4}(a) and $N_0(n)$ is given by equation~\eq{N0-sol}. 
} 
\label{f5p}
\end{center}
\end{figure}

Figure~\ref{f4}(a), showing the Laplacian spectra up $n$ to $6$, suggests that the spectral dimension $d_{s}^{(0)}(n)$ grows with $n$ and should diverge (see also figure~\ref{f5}). Note that, strictly speaking, the spectral dimension is only well-defined in the limit $n\to\infty$. When we refer to the ``spectral dimension'' for finite $n$,  we mean the double effective slope of the curve in the log--log plot of the cumulative spectral density for small values of $\lambda$. This suggests that, in the infinite DSC simplicial complex, the cumulative spectral density of $L^{(0)}$ vanishes faster than any power-law function as $\lambda$ approaches zero. Figure~\ref{f5p} shows that $\lambda_2$ decreases with the increasing number of vertices $N_0$ more slowly than any power law. This behavior suggests that the spectral dimension of the DSC model is infinite. Figure~\ref{f4}(b) demonstrates that for a fixed generation $n$, the cumulative spectral density of $L^{(d)}$ decays more rapidly as $d$ increases. As a result, for any $d$,  the cumulative spectral density of $L^{(d)}$ falls off faster than any power-law function of $\lambda$.

%%%%%%%%%%%%%%%%%%%%%%%%%%%%%%%%%%%%%%%%%%%%%%%
%%%%%%%%%%%%%%%%%%%%%%%%%%%%%%%%%%%%%%%%%%%%%%%

%%%%%%%%%%%%%%%%%%%%%%%%%%%%%%%%%%%%%%%%%%%%%%%
%%%%%%%%%%%%%%%%%%%%%%%%%%%%%%%%%%%%%%%%%%%%%%%

%%%%%%%%%%%%%%%%%%%%%%%%%%%%%%%%%%%%%%%%%%%%%%%
%%%%%%%%%%%%%%%%%%%%%%%%%%%%%%%%%%%%%%%%%%%%%%%

%%%%%%%%%%%%%%%%%%%%%%%%%%%%%%%%%%%%%%%%%%%%%%%
%%%%%%%%%%%%%%%%%%%%%%%%%%%%%%%%%%%%%%%%%%%%%%%

\section{Constrained growth}
\label{s3}

\subsection{The DSC$(m)$ model} 

The definition of the DSC process admits various interesting deformations. We define the DSC$(m)$ model similarly to the DSC model but with the constraint that the dimension of the growing simplicial complex never exceeds $m$. In more precise terms, during each evolution step, every simplex of dimension $d<m$ gets a new vertex attached to it through new edges. We thus have a class of the DSC$(m)$ processes parametrized by positive integers $m=1,2,3,\ldots$. 

The $\mathbf{f}$-vector \eq{decomp} becomes 
\begin{equation}
\label{decomp-m}
\mathcal{N}(n)=[N_0(n),N_1(n),\ldots, N_m(n)]
\end{equation}
for the DSC$(m)$ process. The number of (non-zero) components of the $\mathbf{f}$-vector remains finite, and one can determine it for any $m$. 

\begin{figure}[t]
\begin{center}
\includegraphics[width=0.47\textwidth]{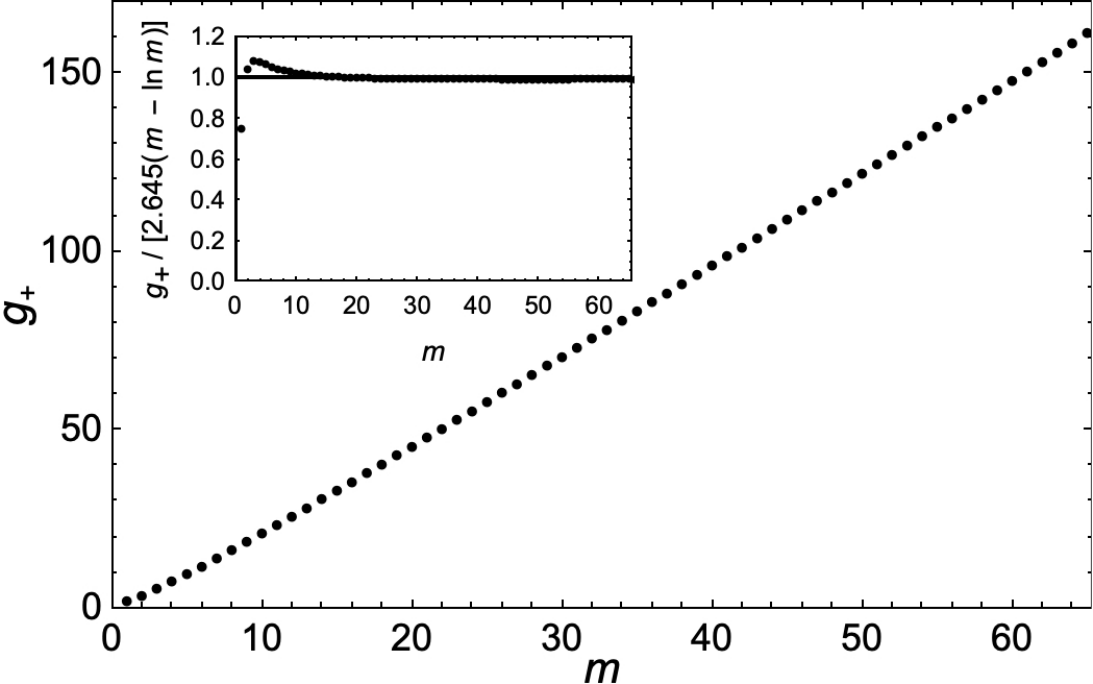}
\caption{The largest eigenvalue $g_+$ of recursion relation~\protect\eq{436} in the DSC$(m)$ model as a function of $m$.  
Inset: Ratio $g_+(m)/[2.645(m - \ln m)]$ vs. $m$, illustrating the quality of the fit~\protect\eq{438}.
} 
\label{f6}
\end{center}
\end{figure}

For the DSC$(m)$ model, the analogs of the DSC recursion relations, equations~\eq{N0}--\eq{N2}, are given by  
\be
N_d(n+1) = N_d(n) + \sum_{j=d-1}^{m-1} {j+1  \choose d} N_j(n) 
\label{436}
\ee
for $0 \leq d \leq m$. On the left-hand side, $N_m(n)$ is present only in the recursion $N_m(n+1) = N_m(n) + N_{m-1}(n)$. The recursions for $N_d(n)$, $0 \leq d \leq m-1$ are governed by a square matrix evident from equation~\eq{436}. The largest eigenvalue $g_+(m) \equiv \lambda_{\text{max}}(m) > 1$ dominates the leading asymptotic behavior: 
\be
N_d(n,m) \cong C_d(m) [g_+(m)]^n \,. 
\label{437}
\ee
We emphasize that $g_+(m)$ depnds only on $m$ and the amplitude $C_d(m)$ depends only of $m$ and $d$; the generation $n$ appears only as the exponent. Figure~\ref{f6} shows how $g_+$ depends on $m$. Numerical results indictae that dependence is nearly linear for large $m$. We tried various sub-leading terms and found that a logarithmic correction (see the inset in figure~\ref{f6})
\be
g_+ \approx 2.645(m - \ln m)
. 
\label{438}
\ee
provides a particularly good fit. Remarkably, the linear and logarithmic terms in the fit \eqref{438} share the same coefficient. In the $m \to \infty$ limit, the DSC$(m)$ model converges to the unconstrained DSC model. This convergence is illustrated in figure~\ref{f6} and equation~\eq{438}. 

%%%%%%%%%%%%%%%%%%%%%%%%%%%%%%%%%%%%%%%%%%%%%%%
%%%%%%%%%%%%%%%%%%%%%%%%%%%%%%%%%%%%%%%%%%%%%%%

%%%%%%%%%%%%%%%%%%%%%%%%%%%%%%%%%%%%%%%%%%%%%%%
%%%%%%%%%%%%%%%%%%%%%%%%%%%%%%%%%%%%%%%%%%%%%%%

For the DSC$(1)$ and DSC$(2)$ models, one can deduce exact formulas for the largest eigenvalue $g_+(m)$ and the amplitudes $C_d(m)$ in \eqref{437}. Several structural properties of the DSC$(1)$ model were studied in \cite{jung2002geometric,dorogovtsev2006degree,qi2009structural}. It also represents a special case of the $d$-pseudofractal deterministic simplicial complex introduced in \cite{bianconi2020spectral}. The simplicial complexes generated by the DSC$(1)$ process are trees. Thus $N_1=N_0-1$. More precisely,  
\begin{equation}
\mathcal{N}(n)=\big[2^n, 2^n-1\big]
.
\label{439}
\end{equation}
The total number of simplices in this model is 
\be
N = N_0 + N_1 = 2N_0 - 1 = 2N_1 + 1
. 
\label{4393}
\ee

The DSC$(2)$ model generates two-dimensional simplicial complexes. For the DSC$(2)$ model, the recursion relations~\eq{436} for $\mathcal{N}(n)=[N_0(n),N_1(n),N_2(n)]$
read
\bea
&&N_0(n+1) = 2N_0(n)+N_1(n)  
,
\nonumber
\\[3pt]
&&N_1(n+1) = N_0(n)+3N_1(n)  
,
\nonumber
\\[3pt]
&&N_2(n+1)=N_2(n)+N_1(n)
.
%%\end{split}
\label{N012-2}
\eea
The initial conditions are $N_0(0)=1, N_1(0)=N_2(1)=0$. The solution reads 
\bea
\label{N0-2}
&&N_0(n) = \frac{5-\sqrt{5}}{10}\, g_+^n + \frac{5+\sqrt{5}}{10}\, g_-^n
,
\\[3pt]
\label{N1-2}
&&N_1(n) = \frac{1}{\sqrt{5}}\, g_+^n - \frac{1}{\sqrt{5}}\,  g_-^n 
,
\\[3pt]
\label{N2-2}
&&N_2(n) = 1+ \frac{3\sqrt{5}-5}{10}\, g_+^n  - \frac{3\sqrt{5}+5}{10}\,  g_-^n 
,
\eea
where 
\begin{equation*}
\begin{split}
&g_+ \equiv g_+(2) = \frac{5 + \sqrt{5}}{2} = 3.618\ldots\\
&g_- \equiv g_-(2) = \frac{5 - \sqrt{5}}{2} = 1.381\ldots
\end{split}
\end{equation*}
Notice that $N_1(n) = N_0(n) + N_2(n) -1$ and hence the total number of simplices is 
\be
N(n) = N_1(n) + N_0(n) + N_2(n) = 2N_1(n) + 1
. 
\label{4396}
\ee
Interestingly, this expression for $N(n)$ in terms of $N_1(n)$ matches that of the DSC$(1)$ model, equation~\eq{4393}.

\subsection{Upper-degree distributions in the DSC$(1)$, DSC$(2)$, and DSC$(3)$ models}

The upper-degree distributions for the DSC$(1)$ and DSC$(2)$ models can be easily obtained explicitly. The degree sequence in the DSC$(1)$ model is $Q^{(0)}(1) = [\hspace{1.5pt}[1,2]\hspace{1.5pt}]$, $Q^{(0)}(2) = [\hspace{1.5pt}[1,2], [2,2]\hspace{1.5pt}]$, and 
\bea
&&\hspace{-20pt}
Q^{(0)}(n) = \left[\hspace{1.5pt}[1,2^{n-1}], [2,2^{n-2},\ldots, 
\right.
\nonumber
\\[3pt]
&&\hspace{85pt}
\left.
[n-2,4],[n-1,2],[n,2]\hspace{1.5pt}\right]
%%D(n)=\big[2^{n-1},2^{n-2},\ldots,4,2,2\big]
\label{440}
\eea 
for $n\geq 3$, using the notations introduced in equation~\eq{385}. The largest degree is $M^{(0)}(n) = n$ for the DSC$(1)$ model. 

The degree distribution 
\begin{equation}
p^{(0)}(k,n) = \frac{D^{(0)}(k,n)}{N_0(n)}
\label{450}
\end{equation} 
is $p^{(0)}(k,n)=2^{-k}$ when $k<n$, where we have used the shorthand notation $k \equiv k^{(0)}$. The degree distribution is remarkably stationary, independent on $n$ for $k<n$; only the fraction of vertices of the highest degree is non-stationary: $p^{(0)}(n,n) = p^{(0)}(n-1,n) =2^{-(n-1)}$. Thus 
\begin{equation}
p^{(0)}(k,n) = 
\left\{
\begin{array}{ll}
2^{-k}  ,         & \ k<n 
,
\\[3pt]
2^{-(n-1)}  ,    & \ k=n 
.
\label{460}
\end{array}
\right.
\end{equation} 
The same exponential degree distribution $p^{(0)}(k) = 2^{-k}$ characterizes recursive random trees, a paradigmatic parameter-free model of growing random trees that exhibits amusing behaviors, see \cite{Pittel94,KR01,KR02,Janson05,Drmota,Frieze,Janson15,Janson19}. The deterministic DSC$(1)$ model is simpler than the recursive random tree model and therefore more tractable. We briefly show how to analyze the DSC$(1)$ model, thereby demonstrating the methods applicable to the DSC$(m)$ models with $m\geq 2$, which are much more challenging for exact analyses. 

The sequence of $1$-degrees in the DSC$(1)$ model evolves as follows: For $n=1$, 
\be
Q^{(0)}(1) = [\hspace{1.5pt}[1,2]\hspace{1.5pt}]
, 
\label{461}
\ee
which corresponds to 
\be
p^{(0)}(k,1) = 2 \delta_{k,1}
;
\label{462}
\ee
for $n=2$, 
\be
Q^{(0)}(2) = [\hspace{1.5pt}[1,2],[2,1],[3,2]\hspace{1.5pt}] 
, 
\label{463}
\ee
which corresponds to 
\be
p^{(0)}(k,2) = 2 \delta_{k,1} + 1\delta_{k,2} + 2\delta_{k,3}
;
\label{464}
\ee
for $n=3$, 
\be
Q^{(0)}(3) = [\hspace{1.5pt}[1,5],[2,5],[3,2],[5,1],[7,2]\hspace{1.5pt}] 
, 
\label{465}
\ee
which corresponds to 
\be
p^{(0)}(k,3) = 5\delta_{k,1} + 5\delta_{k,2} + 2\delta_{k,3} + 1\delta_{k,5} + 2\delta_{k,7}
, 
\label{466}
\ee
and so on. 
In general, for $n\geq2$, we have the recursion relation 
\bea
p^{(0)}(k,n) &=& N_0(n-1)\delta_{k,1} + N_1(n-1)\delta_{k,2} 
\nonumber
\\[3pt]
&+& p^{(0)}((k-1)/2,n-1)
, 
\label{467}
\eea
where $k \geq 3$ in the third term on the right-hand side. 
For the spectrum of distinct degrees, which ignores their multiplicities, this gives
\bea
K^{(0)}(1) &=& [1] 
, 
\nonumber 
\\[3pt] 
K^{(0)}(2) &=& [1,2,3] 
, 
\nonumber 
\\[3pt] 
K^{(0)}(3) &=& [1,2,3,5,7] 
, 
\nonumber 
\\[3pt] 
K^{(0)}(4) &=& [1,2,3,5,7,11,15] 
, 
\nonumber 
\\[3pt] 
K^{(0)}(5) &=& [1,2,3,5,7,11,15,23,31] 
, 
\nonumber 
\\[3pt] 
K^{(0)}(6) &=& [1,2,3,5,7,11,15,23,31,47,63] 
,
\label{468}
\eea
etc. The distinct degrees $k_i$ with labels $i = 1,\ldots, 2n-1$ are given by 
\be
k_i = 
\left\{
\begin{array}{ll}
2^{(i+1)/2} - 1 ,  & \   \text{odd} \ i 
,
\\[3pt]
3\cdot 2^{(i/2) - 1} - 1  ,   & \  \text{even} \ i 
. 
\label{469}
\end{array}
\right.
\ee
Therefore, asymptotically, $p^{(0)}(2k,n)/p^{(0)}(k,n) \cong g_+$, where $g_+ = (5 + \sqrt{5})/2$ is obtained from equations~\eq{N0-2}--\eq{N2-2} and \eq{467}. Hence the degree distribution has a power-law shape, and the exponent $\gamma^{(0)}$ equals 
\be
\gamma^{(0)} = \frac{\ln(5 + \sqrt{5})}{\ln 2} = 2.855 \ldots
,
\label{4691}
\ee
where we used a cumulative degree distribution, since spaces between distinct degrees grow with increasing $k$. 

The sequence of $1$-degrees in the DSC$(2)$ model evolves as follows:  
For $n=2$,  
\be
Q^{(1)}(2) = [\hspace{1.5pt}[0,2],[1,3]\hspace{1.5pt}]
, 
\label{4692}
\ee
which corresponds to 
\be
p^{(1)}(k,2) = 2 \delta_{k,0} + 3 \delta_{k,1}
,  
\label{4693}
\ee
where we use the notation $k \equiv k^{(1)}$ for $1$-degrees for the sake of brevity;  for $n=3$,  
\be
Q^{(1)}(3) = [\hspace{1.5pt}[0,5],[1,12],[2,3]\hspace{1.5pt}]
, 
\label{4694}
\ee
which corresponds to 
\be
p^{(1)}(k,3) = 5 \delta_{k,0} + 12 \delta_{k,1} + 3 \delta_{k,2}
;  
\label{4695}
\ee
for $n=4$,  
\be
Q^{(1)}(4) = [\hspace{1.5pt}[0,15],[1,45],[2,12],[3,3]\hspace{1.5pt}]
, 
\label{4696}
\ee
which corresponds to 
\be
p^{(1)}(k,4) = 15 \delta_{k,0} + 45 \delta_{k,1} + 12 \delta_{k,2} + 3 \delta_{k,3}
,
\label{4697}
\ee
and so on. Generally for $n\geq 2$, 
\begin{equation}
\label{p1-0n}
p^{(1)}(0,n) = N_0(n-1), \quad  p^{(1)}(1,n) = \frac{3}{5}\,N_1(n) 
\end{equation}
and 
\be
p^{(1)}(k,n) = p^{(1)}(k-1,n-1), \ \ \ \ \ k \geq 2 
.
\label{p1-kn}
\ee
Combining \eqref{p1-0n}--\eqref{p1-kn} we obtain
\begin{eqnarray}
p^{(1)}(k,n) &=&\frac{3}{5}\,N_1(n-k+1)  \nonumber \\
                   &\cong& 3\,\frac{\sqrt{5}+1}{10} \left( \frac{5+\sqrt{5}}{2} \right)^{\!n-k}
\label{46995}
\end{eqnarray}
where in the last step we used equation \eqref{N1-2}. This exponential decay of the $1$-degree distribution is in sharp contrast with the power-law behavior of the $0$-degree distribution for this model. 

The number $(5+\sqrt{5})/2$ in equation~\eq{46995} is the largest eigenvalue $g_+(2)$ of the linear recursion relation \eq{N012-2} for the vector $\mathcal{N}(n)=[N_0(n),N_1(n),N_2(n)]$ in the DSC$(2)$ model. This eigenvalue determines the leading asymptotic behavior of $\mathcal{N}(n)$, equation~\eq{4396}.  

Finally, let us briefly describe the upper-degree distributions in the DSC$(3)$ model. We observe that the cumulative $0$- and $1$-degree distributions $p_{\text{cum}}^{(0)}(k^{(0)})$ and $p_{\text{cum}}^{(1)}(k^{(1)})$ follow power-law behavior, and exponents $\gamma^{(0)}$ and $\gamma^{(1)}$ approximately equal $2.35$ and $3.9$, respectively. The $2$-degree distribution $p^{(2)}(k^{(2)}, n) = D^{(2)}(k^{(2)}, n)/N_2(n)$ in the DSC$(3)$ model decays exponentially, as can be easily explained. 
Indeed, notice that
\be
D^{(2)} (k) = 0, n) = N_1 (n - 1)
%%,   
%%if   n \geq 2
\label{46996}
\ee
if $n \geq 2$, and
\be
D^{(2)} (k) \geq 1, n) = 4 \left[N_ 2 (n - k) - N_ 2 (n - k - 1)\right]
%%,    if   n \geq 3   and   1 \leq k^{(2)} \leq n-2
\label{46997}
\ee
if  $n \geq 3$  and $1 \leq k \leq n-2$, and we readily get the asymptotics 
\be
p^{(2)} (k) \cong 3.271 \times 5.491^{-k}
\label{46998}
\ee
as $n \to \infty$. 
The number $5.491=g_+(3)$ determines the leading asymptotic behavior of $N_d(n)$ in the DSC$(3)$ model, see equation~\eq{437}. In general, we suggest that the $(m{-}1)$-degree distribution in the DSC$(m)$ model decays exponentially as $p^{(m-1)} (k) \sim [g_+(m)]^{-k}$.

\subsection{Adjacency matrix spectrum of the DSC$(1)$ model}

We begin by examining the spectrum of the adjacency matrix. This problem proves to be more straightforward than that for the Laplacian spectrum. 

Let $\Upsilon_n(\lambda)$ be the characteristic polynomial of the adjacency matrix of the DSC$(1)$ tree of generation $n$. 
The characteristic polynomial of a tree, $\Upsilon_{\cal T}(\lambda)$, coincides with the matching polynomial \cite{godsil1981theory,Diestel}: $\Upsilon_{\cal T}(\lambda) = \sum_{k=0}^{\lfloor N_0/2\rfloor} (-1)^k {\cal M}_k({\cal T}) \lambda^{N_0 - 2k}$, where ${\cal M}_k({\cal T})$ is the number of matchings of size $k$. For our tree of $N_0 = 2^n$ vertices, the characteristic polynomial contains only even powers of $\lambda$ and it is an even function, $\Upsilon_n(\lambda) = \Upsilon_n(-\lambda)$. Clearly, $\Upsilon_1(0)=-1$ and $\Upsilon_{n\geq2}(0)=1$. Furthermore, since an DSC$(1)$ tree is composed of two identical branches connected by a central edge, $\Upsilon_n(\lambda)$ can be represented as the following product of two polynomials:  
\be
\Upsilon_n(\lambda) = \Omega_n(\lambda) \Omega_n(-\lambda)
,
\label{475}
\ee
where $\Omega_n(0)=1$ for $n\geq2$. The characteristic polynomial $\Upsilon_n(\lambda)$ is a polynomial of $\lambda$ of power $2^n$, while $\Omega_n(\lambda)$ is a polynomial of $\lambda$ of power $2^{n-1}$. 

To figure out the characteristic polynomials, we use Mathematica to compute $\Upsilon_n(\lambda)$ for $n \leq 9$, spot a recurrence relation, and guess that this pattern continues for larger $n$.    
For $n=1$, we have   
\be
\Upsilon_1 = (-1+\lambda) (1+\lambda)
. 
\label{U1}
\ee
Then, for $n=2$,   
\bea
\Upsilon_2 &=& (-1-\lambda+\lambda^2) (-1+\lambda+\lambda^2)  
\nonumber
\\[3pt]
%%\Upsilon_2 
&=& \Upsilon_1^2 - \lambda^2.
\label{495}
\eea
For $n=3$,   
\bea
\Upsilon_3 &=& (1+\lambda-3 \lambda^2-\lambda^3+\lambda^4) (1-\lambda-3 \lambda^2+\lambda^3+\lambda^4)  
\nonumber
\\[3pt]
&=& \Upsilon_2^2 - \lambda^2 \,\Upsilon_1^2
.
\label{505}
\eea
For $n=4$:  
\bea
\Upsilon_4 &=& (1+\lambda-7 \lambda^2-4 \lambda^3+13 \lambda^4+4 \lambda^5-7 \lambda^6-\lambda^7+\lambda^8) 
\nonumber
\\[3pt]
&\times& (1-\lambda-7 \lambda^2+4 \lambda^3+13 \lambda^4-4 \lambda^5-7 \lambda^6+\lambda^7+\lambda^8)  
\nonumber
\\[3pt]
&=& \Upsilon_3^2 - \lambda^2 \,\Upsilon_1^2 \,\Upsilon_2^2
,
\label{515}
\eea
and so on for the first $9$ generations. Hereinafter we often shortly write $\Upsilon_n$ instead of $\Upsilon_n(\lambda)$. We similarly suppress $\lambda$ in other polynomials when it does not lead to confusion.

Thus, to construct the characteristic polynomials we begin with $\Upsilon_1$ given by \eqref{U1} and iterate 
\begin{equation}
\label{Un}
\Upsilon_{n+1}= \left( \Upsilon_{n} - \lambda \prod_{i=1}^{n-1} \Upsilon_i \right)  \left( \Upsilon_{n} + \lambda \prod_{i=1}^{n-1} \Upsilon_i \right)
\end{equation}
for $n\geq 1$. We use the convention $\prod_{i=1}^{0} \equiv 1$. 

For the resulting polynomials $\Omega_n(\lambda)$ and $\Upsilon_n(\lambda)$, a few terms containing lowest and highest powers of $\lambda$ can be written explicitly.  One finds
\begin{eqnarray}
\label{Omega-n}
\Omega_n &=& 1 + \lambda^{2^{n-1}} -\left(\lambda - \lambda^{2^{n-1}-1}\right) \nonumber \\
&-&(2^{n-1} - 1)\left( \lambda^2 + \lambda^{2^{n-1}-2}\right) \nonumber \\
&+&(2^{n-1} - n)\left( \lambda^3 - \lambda^{2^{n-1}-3}\right) \nonumber \\
&+&(2^{2n-3} - 7\cdot 2^{n-2} + 2n + 1)\left( \lambda^4 + \lambda^{2^{n-1}-4}\right) \nonumber \\
&+& \ldots
\end{eqnarray}
for $n \geq 5$ and 
\begin{eqnarray}
\label{Y-n}
\Upsilon_n &=& 1 + \lambda^{2^{n}} -(2^{n} - 1)\left( \lambda^2 + \lambda^{2^{n}-2}\right) \nonumber \\
&+&(2^{2n-1} - 7\cdot 2^{n-1} + 2n + 3)\left( \lambda^4 + \lambda^{2^{n}-4}\right) \nonumber \\
&+& \ldots
\end{eqnarray}
 for $n \geq 4$. 
 
%%%%%%%%%%%%%%%%%%%%%%%%%%%%%%%%%%%%%%%%%%%%%%%
%%%%%%%%%%%%%%%%%%%%%%%%%%%%%%%%%%%%%%%%%%%%%%%

\begin{figure}[t]
\begin{center}
%%%%%%%%%%%%%%%%%%%%%%%%%%\hspace{8pt} \includegraphics[width=0.58\textwidth]{D-f7a.pdf}
\hfill \includegraphics[width=0.46\textwidth]{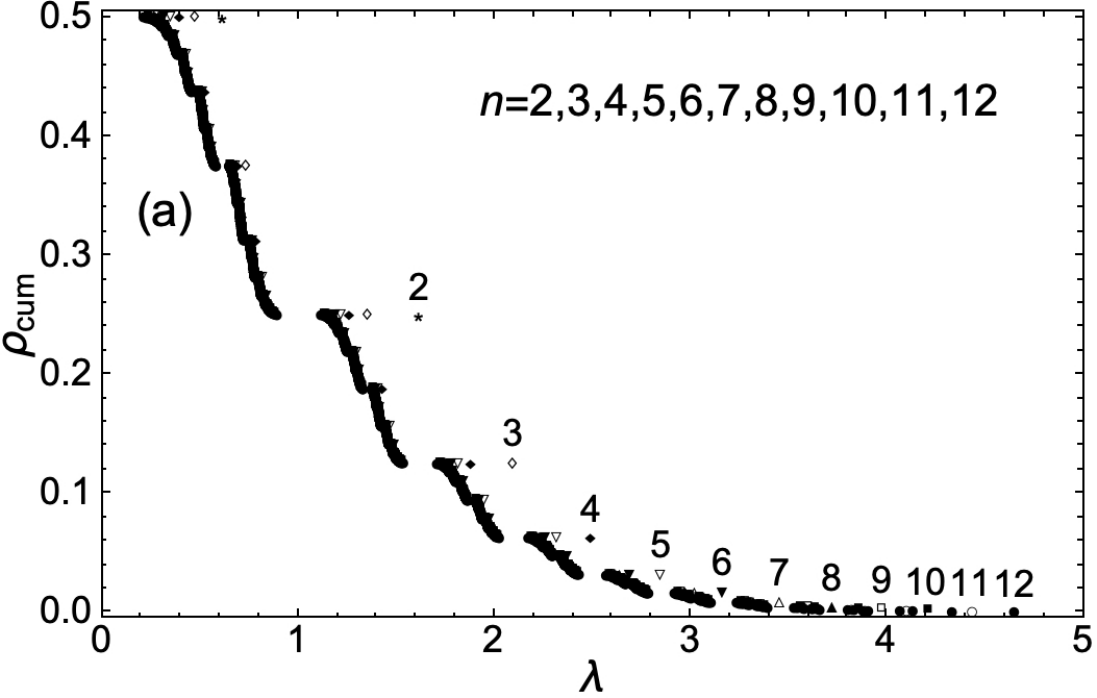} 
\\[3pt]
\hfill \includegraphics[width=0.48\textwidth]{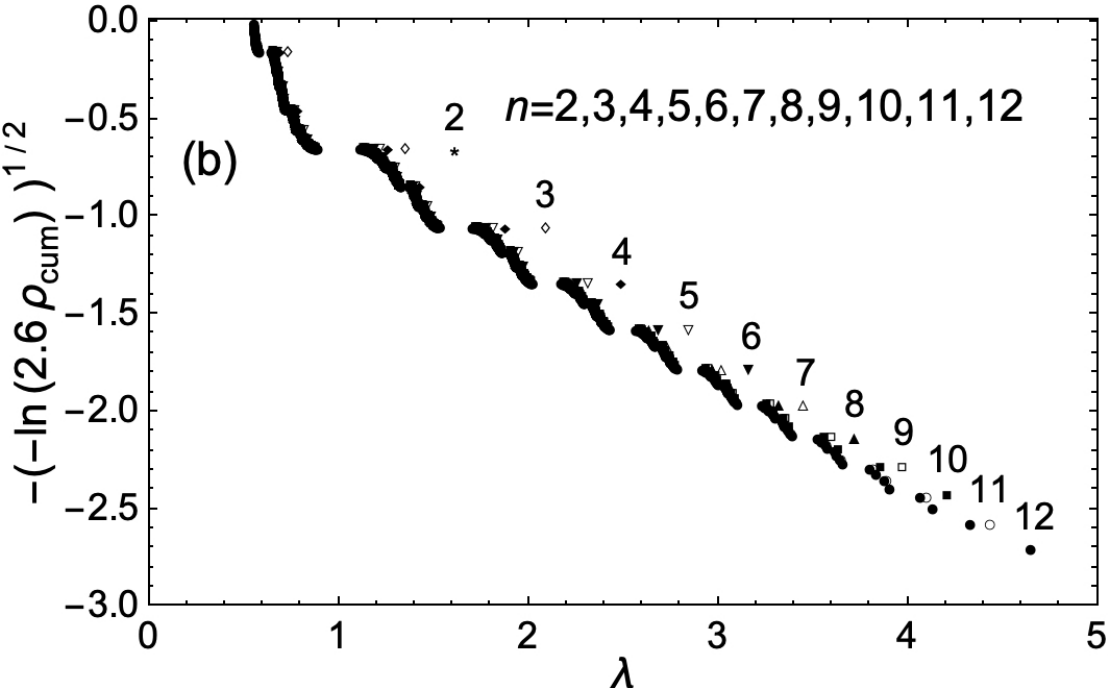}
\caption{(a) Cumulative spectral density $\rho_{\text{cum}}(\lambda)$ of the adjacency matrix for the DSC$(1)$ model, shown for several generations $n$. 
(b) Dependence of $-\sqrt{-\ln (2.6 \rho_{\text{cum}})}$ on $\lambda$, highlighting the large-$\lambda$ asymptotics of the cumulative spectral density $\rho_{\text{cum}} \cong \text{const}\, e^{-0.35 \lambda^2}$. 
Markers represent different values of $n$: 
{\Large \raisebox{-1pt}{$\bullet$}}, {\Large \raisebox{-1pt}{$\circ$}}, 
%%$\Circle[f]$, $\Circle$, 
%%%%%%%%%%%%%%%%%%%%%%%%%%%%%%%%%%%%%%%%%$\fullsquare$,  ${\scriptsize \opensquare}$, 
{\scriptsize $\blacksquare$}, {\scriptsize $\square$}, 
$\blacktriangle$,  
$\protect\rotatebox[origin=c]{90}{$\rhd$}$, $\blacktriangledown$, $\triangledown$, 
%%%%%%%%%%%%%%%%%%%%%%%%%%%%%$\protect\rotatebox[origin=c]{45}{$\fullsquare$}$, $\protect\rotatebox[origin=c]{45}{${\scriptsize \opensquare}$}$, 
$\protect\rotatebox[origin=c]{45}{{\scriptsize $\blacksquare$}}$, $\protect\rotatebox[origin=c]{45}{{\scriptsize $\square$}}$, 
and {\Large $\star$}  
correspond to $n=12$ to $2$, respectively. 
The labels on the plot mark the largest eigenvalues in the spectra. 
} 
\label{f7}
\end{center}
\end{figure}

%%%%%%%%%%%%%%%%%%%%%%%%%%%%%%%%%%%%%%%%%%%%%%%
%%%%%%%%%%%%%%%%%%%%%%%%%%%%%%%%%%%%%%%%%%%%%%%

Using the recursion relations \eq{Un}, we obtained the characteristic polynomials and computed the adjacency matrix spectra for $n \leq 12$. In figure~\ref{f7}(a) we plot the cumulative spectral density $\rho_{\text{cum}}(\lambda)$ at various $n$, 
\be
\rho_{\text{cum}}(\lambda) \equiv \frac{\mathfrak{N}_{\geq\phantom{\!\!\!\!\! y_y\!\!}}(\lambda)}{N_0}
,
\label{535}
\ee
where $\mathfrak{N}_{\geq\phantom{\!\!\!\!\! y_y\!\!}}(\lambda)$ is the number of eigenvalues of the adjacency matrix larger or equal to $\lambda$. 
Unlike the Hodge Laplacian spectra, the cumulative spectral density for the adjacency matrix is defined by summing from above, since we are primarily interested in the behavior at large $\lambda$. 
Figure~\ref{f7}(b) suggests 
\be
\rho_{\text{cum}}(\lambda) \cong \text{const}\, e^{-B \lambda^2}
,
\label{545}
\ee
where the coefficient in the exponential is $B = 0.35$. 
Notably, it is the tail of the cumulative spectral density---not the spectral density itself---that aligns with the Gaussian distribution. 
Note also that as $n \to \infty$, the cumulative spectral density approaches a stationary limit. 

\begin{figure}[t]
\begin{center}
\includegraphics[width=0.47\textwidth]{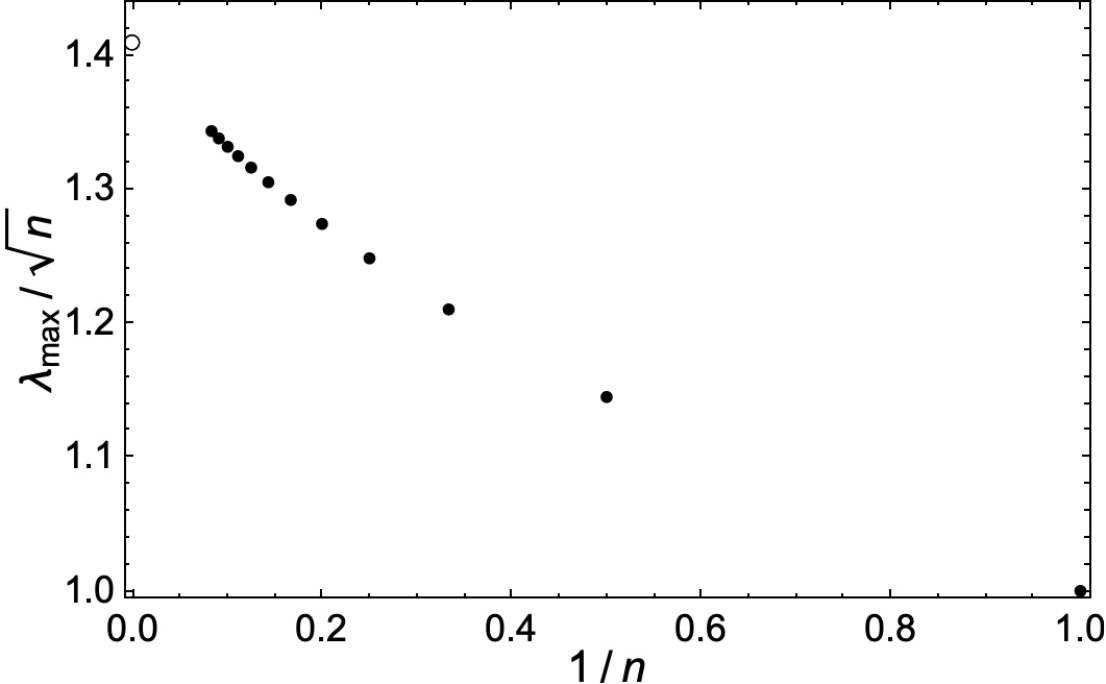}
\caption{The plot of $\lambda_{\text{max}}/\sqrt{n}$ vs. $1/n$ for the DSC$(1)$ model, illustrating the large-$n$ asymptotic behavior $\lambda_{\text{max}} \cong 1.41 \sqrt{n}$. The marker
{\Large \raisebox{-1pt}{$\circ$}}  
%%$\Circle$ 
indicates the extrapolated value $1.41$ as $n \to \infty$.  
} 
\label{f8}
\end{center}
\end{figure}

Inspecting the largest eigenvalues $\lambda_{\text{max}}(n)$ in the range $1 \leq n \leq 12$, we conclude that
\be
\lambda_{\text{max}}(n) \cong 1.41(1) \sqrt{n}
%%.
\label{555}
\ee
for large $n$, see figure~\ref{f8}. Equation~\eq{555} agrees with equation \eq{545}. Indeed, $\lambda_{\text{max}}$ can be estimated from the extreme statistics argument 
\be
2^n 
%%\int_{\lambda_{\text{max}}}^\infty \!\! d\lambda\, \rho_{\text{cum}}(\lambda) 
\rho_{\text{cum}}(\lambda_{\text{max}})
\sim 1
, 
\label{565}
\ee
taking into account that $\rho_{\text{cum}}(\lambda)$ is stationary in the limit $n \to \infty$. Remarkably, substituting $\rho_{\text{cum}}(\lambda)$ from \eq{545} into \eqref{565} not only reveals the correct functional form of $\lambda_{\text{max}}(n)$ in equation~\eq{555}, but also yields the coefficient in the asymptotic expression---albeit with lower precision, estimated as $1.40(5)$. 

We focused on the large-$\lambda$ part of the adjacency matrix spectrum, but also observed that the smallest positive eigenvalue follows a power-law decay: 
\be
\lambda_{\text{min}} \approx 0.75\, n^{-0.50}
. 
\label{567}
\ee
%%

%%%%%%%%%%%%%%%%%%%%%%%%%%%%%%%%%%%%%%%%%%%%%%%
%%%%%%%%%%%%%%%%%%%%%%%%%%%%%%%%%%%%%%%%%%%%%%%

%%%%%%%%%%%%%%%%%%%%%%%%%%%%%%%%%%%%%%%%%%%%%%%
%%%%%%%%%%%%%%%%%%%%%%%%%%%%%%%%%%%%%%%%%%%%%%%

\subsection{Laplacian spectrum of the DSC$(1)$ model}

Let $\Pi_n(\lambda)$ be  the characteristic polynomial of the Laplacian matrix of the DSC$(1)$ tree of generation $n$. As with the adjacency matrix spectrum, to figure out the form of $\Pi_n(\lambda)$ for any finite $n$, we examine $\Pi_n(\lambda)$ using results for $n \leq 9$ directly computed by Mathematica for trees generated by the DSC$(1)$ model, see appendix~\ref{sa2}. We then identify a recurrence relation and conjecture that this pattern continues for larger $n$.  We find
\begin{equation}
\label{Pi-factor}
\Pi_n = -\lambda\prod_{i=1}^n \pi_i
\end{equation}
where $\pi_i(\lambda)$ are polynomials of $\lambda$ of degree $2^{i-1}$. For $1 \leq i \leq 4$, these polynomials are generated by 
\bea
\pi_1 &=& 2 - \lambda
,
\nonumber
\\[3pt]
\pi_2 &=& \pi_1^2 - 2
,
\nonumber
\\[3pt]  
\theta_3 &=& 1 - \lambda
,
\nonumber
\\[3pt]
\pi_3 &=& \pi_2^2 - 2 \theta_3^2
, 
\nonumber
\\[3pt]
\theta_4 &=& \theta_3^2 - \lambda
,
\nonumber
\\[3pt] 
\pi_4 &=& \pi_3^2 - 2 \theta_3^2 \hspace{1pt} \theta_4^2 
.
\label{585}
\eea
For $i\geq 5$, the polynomials $\pi_i$ are also generated recurrenly. In addition to the auxiliary polynomials $\theta_i(\lambda)$ appearing already in \eqref{585}, we also use the auxiliary polynomials $\sigma_i(\lambda)$ and $\widetilde{\sigma}_i(\lambda)$. Namely, one should iterate
\bea
\sigma_n &=& \sigma_{n-1}\widetilde{\sigma}_{n-1} - \lambda \prod_{i=5}^{n-2} \sigma_i\widetilde{\sigma}_i
,
\nonumber
\\[3pt]
\widetilde{\sigma}_n &=& \widetilde{\sigma}_{n-1}^2 - \lambda \prod_{i=5}^{n-2} \sigma_i\widetilde{\sigma}_i
,
\nonumber
\\[3pt]  
 \theta_n &=& \theta_{n-1}^2 - \lambda \prod_{i=5}^{n} \sigma_{i}\widetilde{\sigma}_{i}
, 
\nonumber
\\[3pt]  
\pi_n &=& \pi_{n-1}^2 - 2 \prod_{i=3}^{n} \theta_i^2 
, 
\label{pi-rec}
\eea
and use $\sigma_{4} = 2$ and $\widetilde{\sigma}_{4} = 1$ together with $\theta_3$, $\theta_4$, and $\pi_4$ provided by equation~\eq{585} as initial conditions.

%%Here we present $\pi_i(\lambda)$ for $i\leq 5$: ********
%%%%
%%\begin{equation}
%%\label{pi:1-5}
%%\begin{split}
%%&\pi_1 = 2 - \lambda
%%\\
%%&\pi_2 = 2 - 4 \lambda + \lambda^2\\
%%&\pi_3 = 2 -12 \lambda + 18 \lambda^2 - 8 \lambda^3 + \lambda^4\\
%%&\pi_4 = 2 - 32 \lambda + 168 \lambda^2 - 396 \lambda^3 + 472 \lambda^4 - 296 \lambda^5 \\
%%           & + 98 \lambda^6 - 16 \lambda^7 +  \lambda^8 \\
%%&\pi_5 = 2 - 80 \lambda + 1208 \lambda^2 - 9552 \lambda^3 + 45476 \lambda^4\\
%%            &  - 140600 \lambda^5 + 295460 \lambda^6 - 434172 \lambda^7  + 453962 \lambda^8 \\
%%            &  - 340944 \lambda^9+ 184428 \lambda^{10} - 71516 \lambda^{11} + 19594 \lambda^{12} \\
%%            &  - 3684 \lambda^{13} + 450 \lambda^{14} - 32 \lambda^{15} +  \lambda^{16}
%%\end{split}
%%\end{equation}
%%
%%In Appendix~\ref{sa2}, we present $\pi_6$ and $\pi_7$. 

A few terms of the polynomials $\pi_n(\lambda)$ and $\Pi_n(\lambda)$ with lowest and highest powers of $\lambda$ can be written explicitly.  For $n \geq 3$,  
\be
\pi_n =  2 -  n\hspace{1pt} 2^{n - 1} \lambda - \ldots  -  2^{n} \lambda^{2^{n-1}-1} +  \lambda^{2^{n-1}} 
,
\label{600}
\ee
and 
for $n \geq 2$, 
\bea
&&\hspace{-25pt} 
\Pi_n = \lambda \left[ -2^n +  2^{n - 1} \left(1 + (n - 1) 2^n \right) \lambda - \ldots  
\right.
\nonumber
\\[3pt]
&&\hspace{62pt}
\left. 
-  \left(2^{n+1} - 2\right)\lambda^{2^n-2} +  \lambda^{2^n-1} \right]
.
\label{610}
\eea
%%

%%%%%%%%%%%%%%%%%%%%%%%%%%%%%%%%%%%%%%%%%%%%%%%
%%%%%%%%%%%%%%%%%%%%%%%%%%%%%%%%%%%%%%%%%%%%%%%

\begin{figure}[t]
\begin{center}
\includegraphics[width=0.47\textwidth]{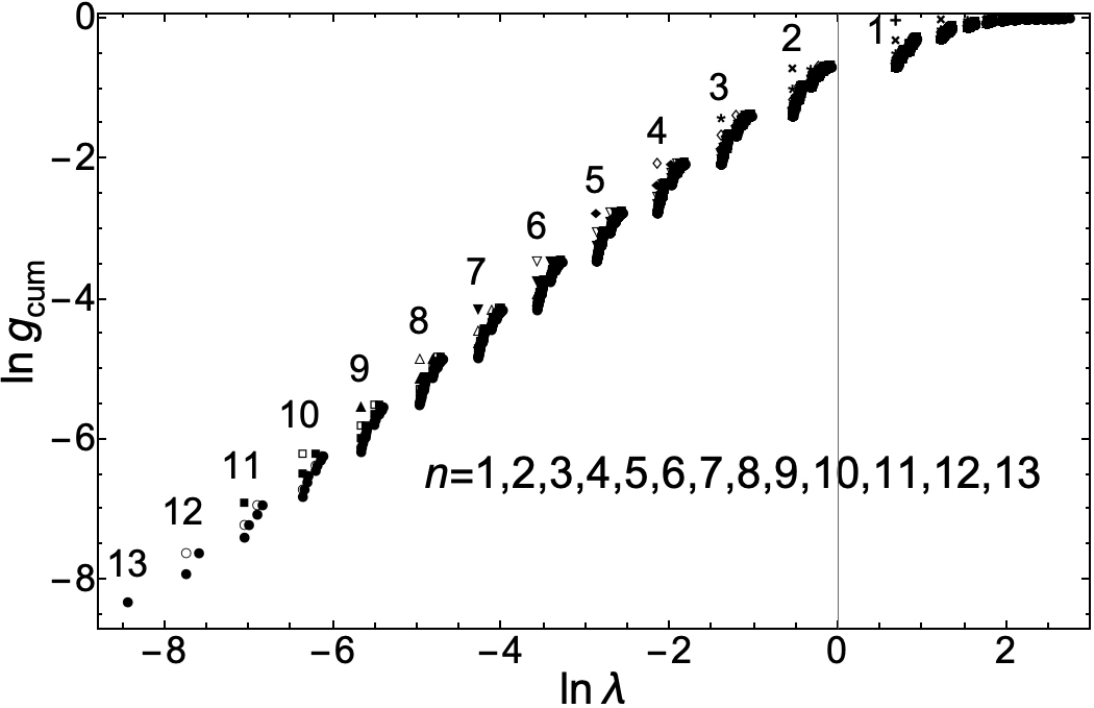}
\caption{Log--log plots of the cumulative Laplacian spectral density $g_{\text{cum}}(\lambda)$ 
%%of a Laplacian matrix 
for trees 
%%simplicial complexes 
generated by the DSC$(1)$ model at various values of $n$. 
Markers represent different values of $n$: 
{\Large \raisebox{-1pt}{$\bullet$}}, {\Large \raisebox{-1pt}{$\circ$}}, 
%%$\Circle[f]$, $\Circle$, 
%%%%%%%%%%%%%%%%%%%%%%%%%%%%%%%%%%%%$\fullsquare$,  $\opensquare$, 
{\scriptsize $\blacksquare$}, {\scriptsize $\square$}, 
$\blacktriangle$, $\protect\rotatebox[origin=c]{90}{$\rhd$}$, $\blacktriangledown$, $\triangledown$, 
%%%%%%%%%%%%%%%%%%%%%%%%%%%%%%%%%%%%$\protect\rotatebox[origin=c]{45}{$\fullsquare$}$, $\protect\rotatebox[origin=c]{45}{${\scriptsize \opensquare}$}$, 
$\protect\rotatebox[origin=c]{45}{{\scriptsize $\blacksquare$}}$, $\protect\rotatebox[origin=c]{45}{{\scriptsize $\square$}}$, 
{\Large $\star$}, $\boldsymbol{\times}$, and $\boldsymbol{+}$ correspond to $n=13$ to $1$, respectively. 
The labels on the plot mark the smallest nonzero eigenvalues $\lambda_2(n)$  
from the Laplacian spectra. 
} 
\label{f9}
\end{center}
\end{figure}

%%%%%%%%%%%%%%%%%%%%%%%%%%%%%%%%%%%%%%%%%%%%%%%
%%%%%%%%%%%%%%%%%%%%%%%%%%%%%%%%%%%%%%%%%%%%%%%

\begin{figure}[t]
\begin{center}
\includegraphics[width=0.47\textwidth]{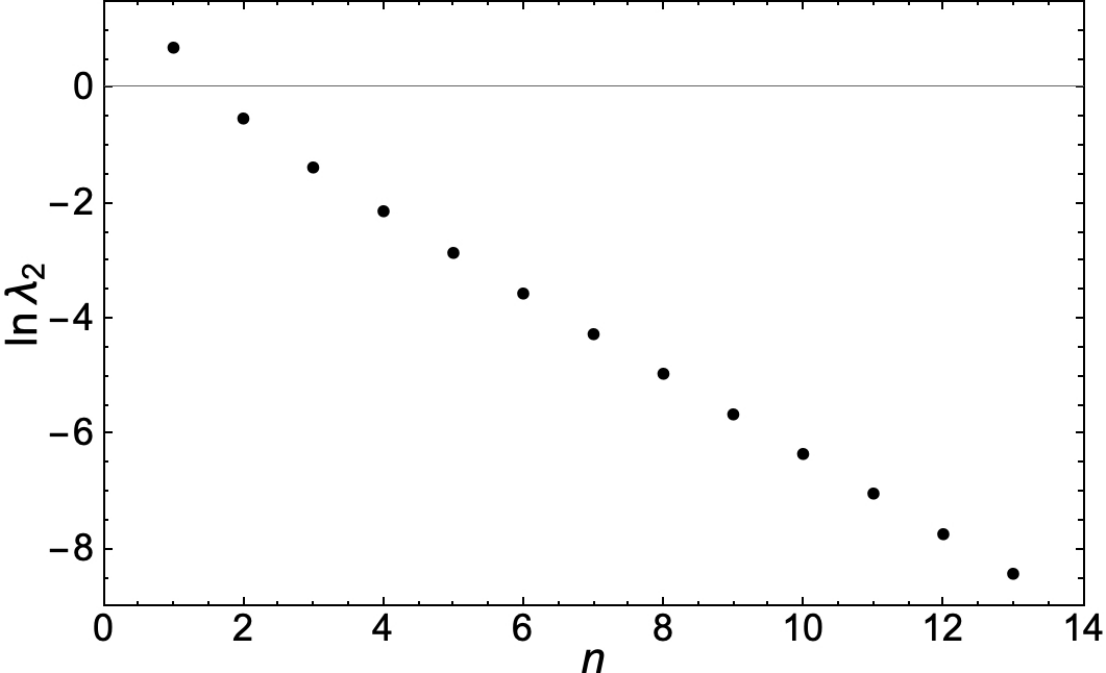}
\caption{
Plot of $\ln \lambda_2$ vs. $n$ showing how the smallest nonzero eigenvalue in the Laplacian spectrum changes with generation $n$ in the DSC$(1)$ trees. 
} 
\label{f10}
\end{center}
\end{figure}

Equations~\eq{Pi-factor}--\eq{pi-rec} enable one to compute the Laplacian matrix spectrum of a DSC$(1)$ tree for any finite $n$ by finding numerically the roots of the equation $\Pi_n(\lambda)=0$. Figure~\ref{f9} show the resulting cumulative densities $g_{\text{cum}}(\lambda,n)$ of the Laplacian matrix spectra for various $n$ up to $n=13$. In this plot, we also indicate the smallest nonzero eigenvalues $\lambda_2(n)$ for various $n$. The smallest eigenvalue in a Laplacian spectrum is zero, $\lambda_1(n) = 0$.  The eigenvalue $\lambda_2$ is also referred to as the Fiedler eigenvalue \cite{fiedler1973algebraic} or the spectral gap. Figure~\ref{f10} shows that $\lambda_2(n)$ decays exponentially with increasing $n$. The large-$n$ exponential asymptotics can be obtained with high precision,  
\be
\lambda_2(n) \cong 1.77525 \  2^{-n}
.
\label{620}
\ee
%%

%%%%%%%%%%%%%%%%%%%%%%%%%%%%%%%%%%%%%%%%%%%%%%%
%%%%%%%%%%%%%%%%%%%%%%%%%%%%%%%%%%%%%%%%%%%%%%%

%%%%%%%%%%%%%%%%%%%%%%%%%%%%%%%%%%%%%%%%%%%%%%%
%%%%%%%%%%%%%%%%%%%%%%%%%%%%%%%%%%%%%%%%%%%%%%%

Figure~\ref{f10} suggests that at small $\lambda$, the spectral density exhibits power-law behavior and hence the spectral dimension $d_s$ is finite. 
We have mentioned that the DSC$(1)$ trees are a special, $d=1$, case of the $d$-pseudofractal deterministic simplicial complexes introduced in \cite{bianconi2020spectral}, where $d_s$ was found only for $d \geq 2$. Let us obtain this spectral dimension for the DSC$(1)$ trees. Notice in figure~\ref{f10} that the eigenvalue $\lambda_2(\tilde{n})$ is present in the Laplacian spectra for all $n \geq \tilde{n}$. Furthermore, for all $\tilde{n} \leq n$, a clear gap exists between this eigenvalue and the closest eigenvalue beneath it. Inspecting the spectra for $n \leq 13$, we notice that the values of $g_{\text{cum}}(\lambda_1(\tilde{n}),n)$ at these special eigenvalues can be written explicitly, 
\be
g_{\text{cum}}(\lambda_2(\tilde{n}),n) = 2^{-\tilde{n}} + 2^{-n}
,  
\label{630}
\ee
where $n \geq \tilde{n}$. Thus, for large $n \gg \tilde{n}$, we have the asymptotic 
\be
\label{631}
g_{\text{cum}}(\lambda_2(\tilde{n}),n) \cong 2^{-\tilde{n}}
,  
\ee
where $\lambda_2(\tilde{n}) \sim 2^{-\tilde{n}}$ for $\tilde{n} \gg 1$ according to equation~\eq{620}.  
Consequently, in the limit $n \to\infty$, the log--log plot of the spectrum (figure~\ref{f10}) lies within an inclined strip of slope $1$, with both width and height equal to $ \ln 2$. This corresponds to the spectral dimension $d_s = 2$, see equation~\eq{435}. 
In general, the bonds on $d_s$ are set by the Hausdorff dimension $d_H$, specifically $2d_H /(d_H + 1) \leq d_s \leq d_H$   
\cite{durhuus2009hausdorff}. Since our simplicial complexes have an infinite Hausdorff dimension, the resulting $d_s=2$ is the minimal possible spectral dimension. 

We can also obtain this value of the spectral dimension by looking directly at how $\lambda_2$ scales with the network size $N_0$. Indeed, equations~\eq{439} and \eq{620} lead to the asymptotic $\lambda_2(N_0) \sim 1/N_0$. Substituting $g_{\text{cum}}(\lambda) \sim \lambda^{d_s/2}$ into the extreme statistics condition  
\be 
%%N_0 \int_0^{\lambda_2(N_0)} \!\! d\lambda\,g(\lambda) = 
N_0 \, g_{\text{cum}}(\lambda_2) \sim 1
,  
\ee
resembling \eq{565}, and assuming that the cumulative spectral density becomes stationary in the $n \to \infty$ limit, cf. \eqref{630}--\eqref{631}, we get  
\be 
\lambda_2(N_0) \sim N_0^{-2/d_s}
.   
\label{650}
\ee
Comparing \eqref{650} with the asymptotic $\lambda_2(N_0) \sim N_0^{-1}$ we again arrive at $d_s=2$. 

Finally, we note that the maximum eigenvalue of the spectrum behaves as $\lambda_{\text{max}} \simeq n$ when $n\gg 1$.  

%%%%%%%%%%%%%%%%%%%%%%%%%%%%%%%%%%%%%%%%%%%%%%%
%%%%%%%%%%%%%%%%%%%%%%%%%%%%%%%%%%%%%%%%%%%%%%%

\begin{figure}[t]
\begin{center}
\includegraphics[width=0.47\textwidth]{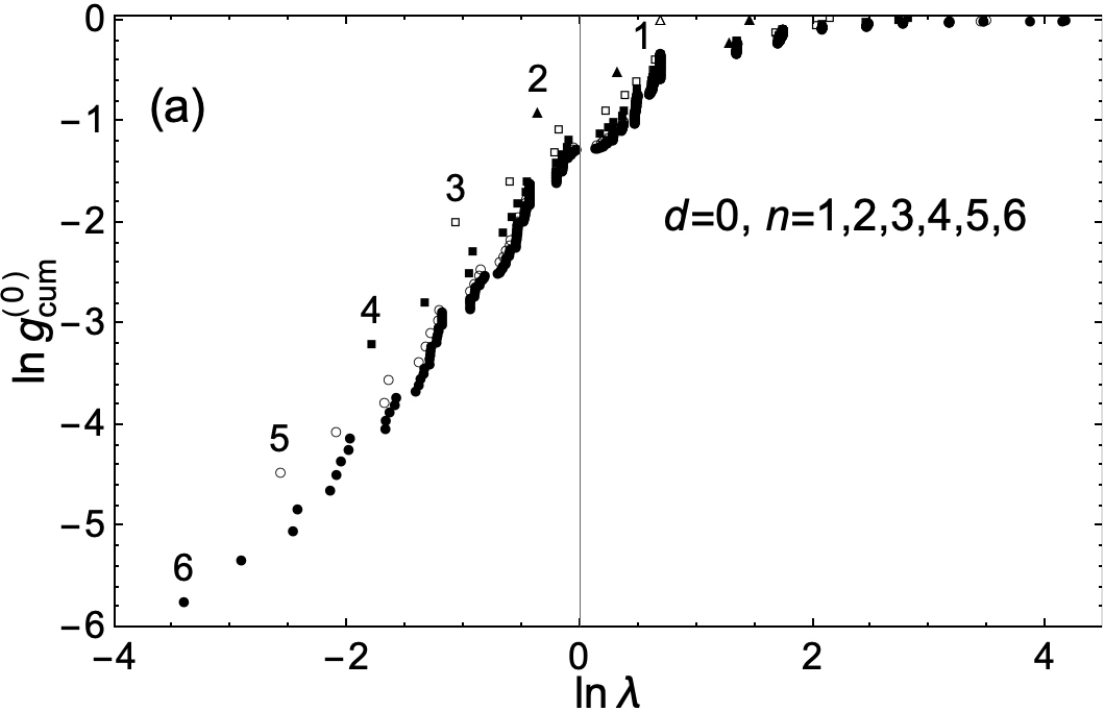}
\\[3pt]
\includegraphics[width=0.47\textwidth]{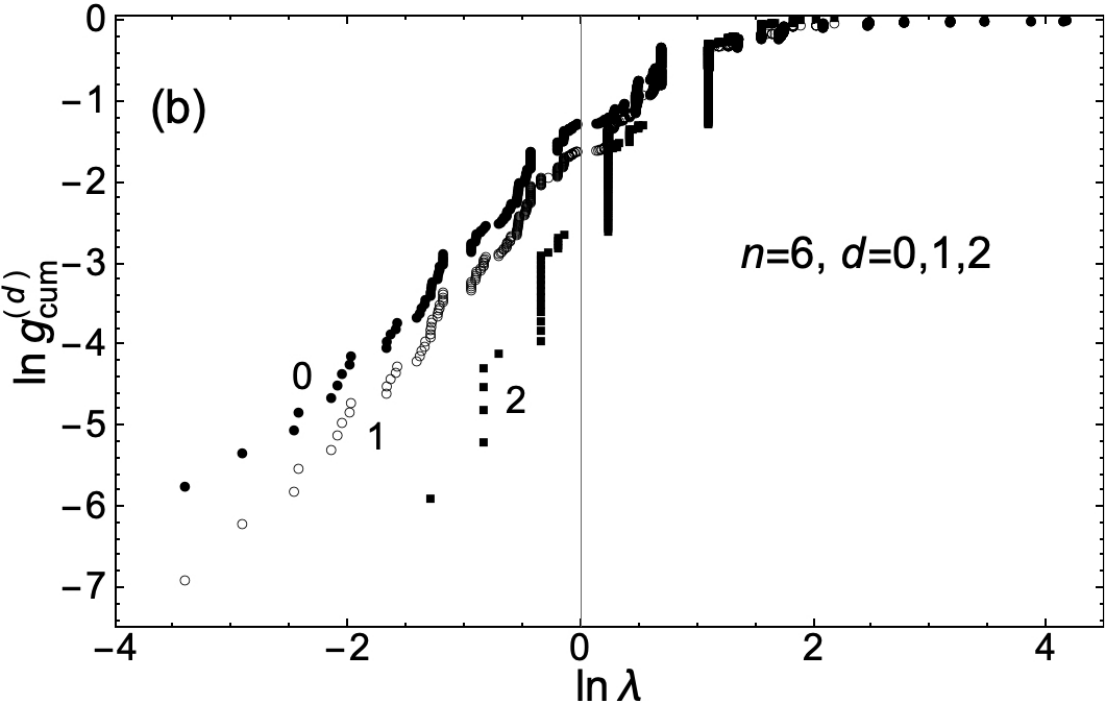}
\caption{The log--log plots of the cumulative spectral densities $g^{(d)}_{\text{cum}}(\lambda,n)$ of Hodge Laplacians for simplicial complexes generated by the DSC$(2)$ model at various values of $n$ and $d$. 
(a)~The log--log plot of the cumulative spectral densities $g^{(0)}_{\text{cum}}(\lambda,n)$ of Laplacians for the $1$-skeletons of the DSC$(2)$ simplicial complexes at different generations $n$. 
Markers represent different values of $n$: 
{\Large \raisebox{-1pt}{$\bullet$}}, {\Large \raisebox{-1pt}{$\circ$}}, 
%%$\Circle[f]$, $\Circle$, 
%%%%%%%%%%%%%%%%%%%%%%%%%%%%%%%%%%%$\fullsquare$,  $\opensquare$, 
{\scriptsize $\blacksquare$}, {\scriptsize $\square$}, 
$\blacktriangle$, and $\protect\rotatebox[origin=c]{90}{$\rhd$}$, 
correspond to $n=6$ to $1$, respectively. 
The labels on the plot mark the smallest nonzero eigenvalues $\lambda_2(n)$  
in the Laplacian spectra. 
(b)~The log--log plots of the cumulative spectral densities $g^{(d)}_{\text{cum}}(\lambda,n)$ of Hodge Laplacians for the DSC$(2)$ simplicial complexes at generation $n=6$ for various $d$. Markers represent different values of $d$: 
{\Large \raisebox{-1pt}{$\bullet$}}, {\Large \raisebox{-1pt}{$\circ$}}, 
%%$\Circle[f]$, $\Circle$, 
{\scriptsize $\blacksquare$}  
correspond to $d=0$, $1$, $2$. 
} 
\label{f11}
\end{center}
\end{figure}

%%%%%%%%%%%%%%%%%%%%%%%%%%%%%%%%%%%%%%%%%%%%%%%
%%%%%%%%%%%%%%%%%%%%%%%%%%%%%%%%%%%%%%%%%%%%%%%

\subsection{Hodge Laplacian spectra of the DSC$(2)$ model}

\begin{figure}[t]
\begin{center}
\includegraphics[width=0.47\textwidth]{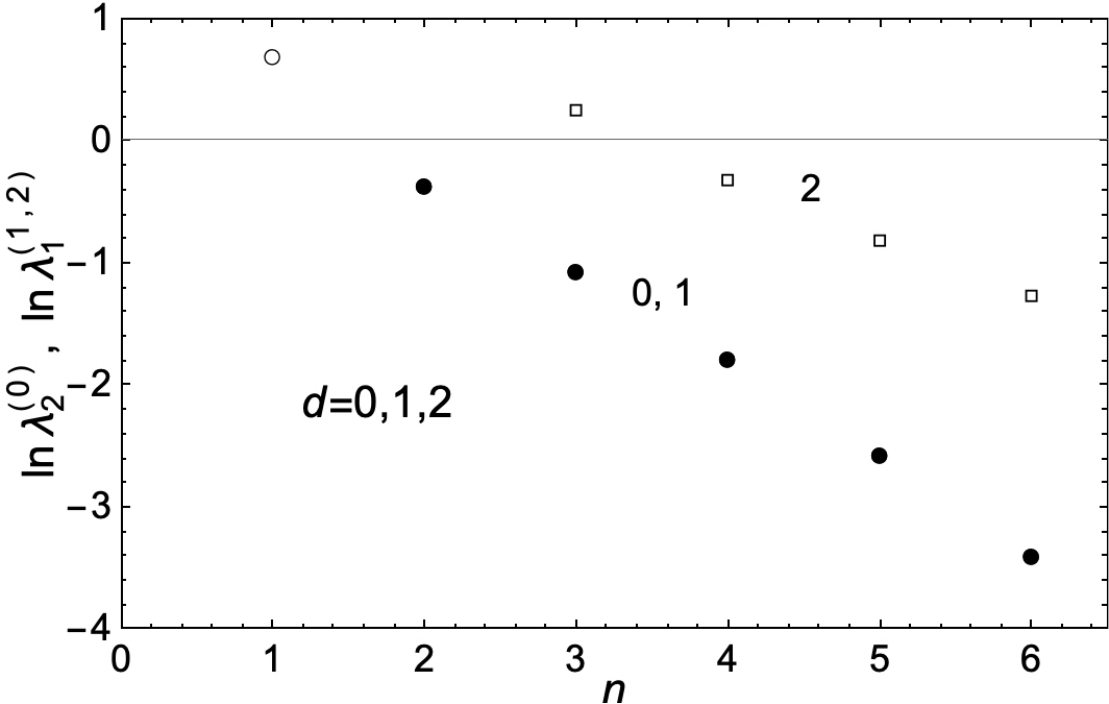}
\caption{
Plot of $\ln \lambda_2^{(0)}$, $\ln \lambda_1^{(1)}$, and $\ln \lambda_1^{(2)}$ vs. $n$ showing how the smallest nonzero eigenvalue in the Hodge Laplacian spectrum changes with generation $n$ in the DSC$(2)$ model. 
Markers represent different values of $d$: 
{\Large \raisebox{-1pt}{$\circ$}}, {\Large \raisebox{-1pt}{$\bullet$}}, 
%%$\Circle$, $\Circle[f]$, 
and 
%%%%%%%%%%%%%%%%%%%%%%%%%%%%%%%%%%%%%%%$\opensquare$ 
{\scriptsize $\square$} 
correspond to $d=0$, $1$, and  $3$, respectively. Notice that for $n \geq 2$, $\lambda_2^{(0)}(n) = \lambda_1^{(1)}(n)$. 
} 
\label{f12}
\end{center}
\end{figure}

Using equations~\eq{410}--\eq{425}, we numerically compute the spectra of the Hodge Laplacians for simplicial complexes generated by the DSC$(2)$ model for generations up to $n = 6$. The log--log plots in figure~\ref{f11} show the cumulative spectral densities of these spectra.  As illustrated in figure~\ref{f11}(a), the cumulative spectral density $g^{(0)}_{\text{cum}}(\lambda)$ of the Laplacian spectra tends toward a stationary limit as $n$ increases. 
Figure~\ref{f11}(b) shows that $g^{(0)}_{\text{cum}}(\lambda)$ becomes increasingly steep at small $\lambda$ as $d$ increases.  
The dependences of the smallest nonzero eigenvalues $\lambda_2^{(0)}$, $\lambda_1^{(1)}$, and $\lambda_1^{(2)}$ on $n$ are close to exponential, see figure~\ref{f12}. The spectra of the 1st and 2nd Hodge Laplacians of these simplicial complexes have no zero eigenvalues. 
By measuring the slope of the dependence $\ln \lambda_2^{(0)}(n)$ in figure~\ref{f12} and applying equations~\eq{N0-2}--\eq{N2-2} and \eq{650}, we estimate the spectral dimension of this model to be $d_s \approx 3.0(4)$.

%%%%%%%%%%%%%%%%%%%%%%%%%%%%%%%%%%%%%%%%%%%%%%%
%%%%%%%%%%%%%%%%%%%%%%%%%%%%%%%%%%%%%%%%%%%%%%%

%%%%%%%%%%%%%%%%%%%%%%%%%%%%%%%%%%%%%%%%%%%%%%%
%%%%%%%%%%%%%%%%%%%%%%%%%%%%%%%%%%%%%%%%%%%%%%%

\section{Conclusions}
\label{s-conclusions}

We investigated the structural properties of simplicial complexes grown via very simple deterministic procedure. In our parameter-free DSC model, at each step, every $d$-simplex of the simplicial complex is transformed into a $(d{+}1)$-simplexl. The diameter of these simplicial complexes, $\mathcal{D}(n)$, grows linearly with the number of evolution steps $n$, while the number of $d$-simplices, $N_d(n)$, scales as $n! (\ln n)^d$, and the total number of simplices, $N(n)$, increases as $n \, n!$. An unexpected observation of Ref.~\cite{Sergey02} is that structural characteristics of deterministic networks are surprisingly close to those of their randomly growing counterparts \cite{dorogovtsev2001size}. Analogously, we anticipate that our results for deterministic simplicial complexes may provide useful insights into the behavior of these characteristics in the corresponding stochastic models. The asymptotics $N_0 \sim n!$ leads to $n \sim \ln N_0 / \ln \ln N_0$ implying the following size dependences in the DSC model: 
\bea 
\  \mathcal{D}(N_0) &\sim& \frac{\ln N_0}{\ln \ln N_0}
,
\\[3pt]
N_d(N_0) &\sim& \frac{N_0(\ln\ln N_0)^d}{d!}
,
\\[3pt]
\, N(N_0) &\sim& \frac{N_0 \ln N_0}{\ln\ln N_0}
. 
\eea
These asymptotics give us a good idea of what to expect from the stochastic version of the DSC model, which is harder to tackle directly.

We have observed power-law $d$-degree distributions in the DSC model. Notably, while the `normal' degree distribution ($0$-degrees) rapidly converges to the infinite-$n$ stationary limiting shape, the infinite-$n$ limit of the upper-degree distributions (that is, for $(d \geq 1)$-degrees) is practically unreachable. Furthermore, we numerically computed the spectra of the Hodge Laplacians for the DSC model and found that the spectral dimension diverges as $n \to \infty$. 

We also analyzed a class of deterministic constrained growth models, the DSC$(m)$ models, parametrized by integer $m\geq 1$. We have shown that both $N_d$ and $N$ grow exponentially with $n$. In the DSC$(1)$ model generating determistic trees, the degree distribution decays exponentially. Interestingly, in the DSC$(2)$ model, the degree distribution follows a power law, while the $1$-degree distribution decays exponentially. 
In the DSC$(3)$ model, the $0$- and $1$-degree distributions are power laws, while the $1$-degree distribution is exponential. 
We have found that in the DSC$(1)$ model, the tail of the cumulative spectral density for the adjacency matrix matches a Gaussian, and the largest eigenvalue scales like $\sqrt{n}$. We have shown that the spectral dimension of the DSC$(1)$ model is $2$ and that the smallest nonzero eigenvalue in the Laplacian spectrum scales as $2^{-n}$. These results are noteworthy, given the surprising scarcity of studies on the spectra of the recursive random trees---the stochastic counterpart of the DSC$(1)$ model. One can get a rough picture of the shape of the adjacency matrix spectra for the recursive random trees from \cite{golinelli2003statistics,bhamidi2012spectra}, but not their asymptotics. Regarding the spectral dimension of random trees, rigorous results have been obtained for various non-growing trees \cite{durhuus2009hausdorff,destri2002spectral}---specifically, it is equal to $2$ when the Hausdorff dimension is infinite---while no such result is available for the recursive random trees. 
Finally, we have obtained the Hodge Laplacian spectra of the DSC$(2)$ model. In particular, unlike the DSC model, the cumulative Laplacian spectral density of the DSC$(2)$ model settles into a stationary power-law shape as $n \to \infty$, and the spectral dimension turns out to be finite.    

Some of our results were obtained by solving recursion relations exactly, others suggested by numerical analyses, and a third set are conjectures that we verified for a sufficiently long sequence of early evolutionary steps. Proving these conjectures is a challenge for future work.

%%%%%%%%%%%%%%%%%%%%%%%%%%%%%%%%%%%%%%%%%%%%%%%\ack

\begin{acknowledgments}

This work was initiated during PLK's visit to the University of Aveiro. PLK is grateful to the Department of Physics of the University of Aveiro for the excellent working conditions. We also want to thank Ginestra Bianconi and Jose Mendes for the discussions. 

\end{acknowledgments}

\appendix

\section{The DSC model in higher generations}
\label{sa1}

Here we present a few structural characteristics of the simplicial complexes generated by the DSC rule in higher generations. The list of simplicies \eqref{decomp} for $n\leq 6$ is given in equations \eqref{N1-6}. For the next three generations 
%%%%%%%%%%%%%%%%%%%%%%%%%%%%%%%%%%%%%%%%%%%%%%%\ack
\begin{widetext}
\begin{equation}
\label{N7-9}
\begin{split}
&\mathcal{N}(7)=[13700, 33375, 30681, 14084,  3535, 490, 35, 1], \\
&\mathcal{N}(8) = [109601,283072,284066,147532,43939,7756,798,44,1],\\
&\mathcal{N}(9) = [986410,2673321,2878284,1646714,558642,117579,15456,1230,54,1]. 
\end{split}
\end{equation}

Equations \eqref{330} give the lists $K^{(d)}(n)$ of distinct $d$-degrees in the simplicial complexes generated by the DSC rule with $n\leq 5$. Looking at these lists and the lists of $K^{(d)}(n)$ with $d\leq n-1$ and $n\leq 8$ one observes that the arrays $K^{(n-p)}(n)$ with fixed $p$ quickly become stationary (i.e., independent on the generation $n$). The stationary values $K^{(n-p)}(n)$ with $n\geq p+1$ and $1\leq p\leq 8$ are 
\bea
&&
\hspace{-72pt} 
K^{(n-1)}(n) = [0,1]
,
\nonumber 
\\[3pt]
&&
\hspace{-72pt} 
K^{(n-2)}(n) = [0,1,2,3]
,
\nonumber 
\\[3pt]
&&
\hspace{-72pt} 
K^{(n-3)}(n) =  [0,1,2,3,6,8]
,
\nonumber 
\\[3pt]
&&
\hspace{-72pt} 
K^{(n-4)}(n) =  [0,1,2,3,4,6,8,11,19,24]
,
\nonumber 
\\[3pt]
&&
\hspace{-72pt} 
K^{(n-5)}(n) = [0,1,2,3,4,5,6,8,11,19,20,24,46,73,89]
,
\nonumber 
\\[3pt]
&&
\hspace{-72pt} 
K^{(n-6)}(n) = [0,1,2,3,4,5,6,8,11,19,20,24,37,46,73,89,117,236,350,415]
,
\nonumber 
\\[3pt]
&&
\hspace{-72pt} 
K^{(n-7)}(n) = [0,1,2,3,4,5,6,7,8,11,19,20,24,37,46,70,73,89,117,236,312,350,415,807,1459,
\nonumber
\\[3pt]
&&
\hspace{-16pt} 
2046,2372]
,
\nonumber
\\[3pt]
&&
\hspace{-72pt} 
K^{(n-8)}(n) = [0,1,2,3,4,5,6,7,8,11,19,20,24,37,46,70,73,89,117,135,236,312,350,415,807,863,
\nonumber
\\[3pt]
&&
\hspace{-16pt}
1459,2046,2372,2878,6352,10527,14115,16072]. 
\label{a10}
\end{eqnarray}
The array $K^{(0)}(p)$ corresponding to $n=p$ lacks zero but otherwise the same as the stationary array $K^{(n-p)}(n)$ with $n>p$; e.g., $K^{(0)}(2)=[1,2,3]$, while $K^{(1)}(3)=K^{(2)}(4)=K^{(3)}(5)=\ldots=[0,1,2,3]$. These results for $K^{(n-p)}(n)$ with $p\leq 8$ and arbitrary $n\geq p$ give all $K^{(d)}(n)$ with $n\leq 8$. We also know all $K^{(d)}(9)$ apart from $K^{(0)}(9)$, all $K^{(d)}(10)$ apart from $K^{(0)}(10)$ and $K^{(1)}(10)$, etc. 

We now turn to the complete upper-degree sequences $Q^{(d)}(n)$ accounting for degeneracies. Inspecting $Q^{(n-1)}(n)$ which we computed  with the help of Mathematica for $n\leq 8$, and by hand for $n\leq 3$, equation \eqref{385}, we arrive at 
\begin{equation}
\label{Qn1}
Q^{(n-1)}(n) = \left[\hspace{1.5pt} \left[0,\tfrac{(n-1)(n+2)}{2}\right], [1,n+1]\hspace{1.5pt}\right]
\end{equation}
which was equivalently stated as \eqref{Dn1n}. Equation \eqref{Qn1} is 
conjecturally valid for all $n\geq 2$. Similarly, we verified 
\begin{equation}
\label{Qn2}
Q^{(n-2)}(n) =  [\hspace{1.5pt}[0, a_n], [1, b_n], [2, c_n], [3, n]\hspace{1.5pt}]
\end{equation}
with 
\begin{equation}
\label{Qn2-abc}
a_n = \frac{(n-2)(n-1)n(3n+11)}{2}\,, \quad
b_n = \frac{n(n^2+n-4)}{2}\,, \quad
c_n =\frac{n(n-1)}{2}
\end{equation}
for $n\leq 8$, and conjecture that \eqref{Qn2}--\eqref{Qn2-abc} hold for all $n\geq 3$. Equations \eqref{Qn2}--\eqref{Qn2-abc} were equivalently stated as \eqref{Dn2n}. 

We now list $Q^{(d)}(n)$ for $4\leq n\leq 8$. Since $Q^{(n-1)}(n)$ and  $Q^{(n-2)}(n)$ can be extracted from \eqref{Qn1}--\eqref{Qn2-abc}, we only list $Q^{(d)}(n)$ with $d\leq n-3$:
\bea
&&
\hspace{-72pt}
Q^{(0)}(4) = [\hspace{1.5pt}[1,16],[2,23],[3,14],[4,1],[6,5],[8,2],[11,1],[19,1],[24,2]\hspace{1.5pt}]
,
\nonumber
\\[3pt]
&&
\hspace{-72pt}
Q^{(1)}(4) = [\hspace{1.5pt}[0, 16], [1, 51], [2, 27], [3, 16], [6, 3], [8, 3]\hspace{1.5pt}]
,
\nonumber
\\[3pt]
&&
\hspace{-72pt}
Q^{(0)}(5) = [\hspace{1.5pt}[1, 65], [2, 116], [3, 81], [4, 14], [5, 1], [6, 23], [8, 5], [11, 9], [19, 5], [20, 1], [24, 2], [46, 1], \nonumber
\\[3pt]
&&
\hspace{-22pt}
[73, 1], [89, 2]\hspace{1.5pt}] 
,
\nonumber
\\[3pt]
&&
\hspace{-72pt}
Q^{(1)}(5) = [\hspace{1.5pt}[0, 65], [1, 248], [2, 195], [3, 107], [4, 5], [6, 27], [8, 12], [11, 4], [19, 3], [24, 3]\hspace{1.5pt}]
,
\nonumber
\\[3pt]
&&
\hspace{-72pt}
Q^{(2)}(5) = [\hspace{1.5pt}[0, 116], [1, 218], [2, 84], [3, 42], [6, 6], [8, 4]\hspace{1.5pt}]
,
\nonumber
\\[3pt]
&&
\hspace{-72pt}
Q^{(0)}(6) =  [\hspace{1.5pt}[1,326],[2,669],[3,535],[4,145],[5,20],[6,117],[8,16],[11,65],[19,23],
\nonumber
\\[3pt]
&&
\hspace{-22pt}
[20,14],[24,5],[37,1],[46,9],[73,5],[89,2],[117,1],[236,1],[350,1],[415,2]\hspace{1.5pt}] 
,
\nonumber
\\[3pt]
&&
\hspace{-72pt}
Q^{(1)}(6) = 
 [\hspace{1.5pt}[0,326],[1,1403],[2,1410],[3,828],[4,100],[5,6],[6,195],[8,51],[11,56],
\nonumber
\\[3pt]
&&
\hspace{-22pt}
[19,27],[20,5],[24,12],[46,4],[73,3],[89,3]\hspace{1.5pt}] 
,
\nonumber
\\[3pt]
&&
\hspace{-72pt}
Q^{(2)}(6) =  [\hspace{1.5pt}[0,669],[1,1526],[2,870],[3,418],[4,15],[6,84],[8,32],[11,10],[19,6],[24,4]\hspace{1.5pt}] 
,
\nonumber
\\[3pt]
&&
\hspace{-72pt}
Q^{(3)}(6) =  [\hspace{1.5pt}[0,470],[1,645],[2,200],[3,85],[6,10],[8,5]\hspace{1.5pt}]
,
\nonumber
\\[3pt]
&&
\hspace{-72pt}
Q^{(0)}(7) =  [\hspace{1.5pt}[1,1957],[2,4429],[3,3960],[4,1415],[5,280],[6,696],[7,1],[8,65],[11,470],
\nonumber
\\[3pt]
&&
\hspace{-22pt}
[19,116],[20,145],[24,16],[37,20],[46,65],[70,1],[73,23],[89,5],[117,14],
\nonumber
\\[3pt]
&&
\hspace{-22pt}
[236,9],[312,1],[350,5],[415,2],[807,1],[1459,1],[2046,1],[2372,2]\hspace{1.5pt}]
,
\nonumber
\\[3pt]
&&
\hspace{-72pt}
Q^{(1)}(7) =  [\hspace{1.5pt}[0,1957],[1,9184],[2,10902],[3,7063],[4,1400],[5,162],[6,1417],[8,248],
\nonumber
\\[3pt]
&&
\hspace{-22pt}
[11,580],[19,195],[20,100],[24,51],[37,6],[46,56],[73,27],[89,12],[117,5],
\nonumber
\\[3pt]
&&
\hspace{-22pt}
[236,4],[350,3],[415,3]\hspace{1.5pt}]
,
\nonumber
\\[3pt]
&&
\hspace{-72pt}
Q^{(2)}(7) =  [\hspace{1.5pt}[0,4429],[1,11571],[2,8490],[3,4326],[4,405],[5,21],[6,870],[8,218],
\nonumber
\\[3pt]
&&
\hspace{-22pt}
[11,200],[19,84],[20,15],[24,32],[46,10],[73,6],[89,4]\hspace{1.5pt}]
,
\nonumber
\\[3pt]
&&
\hspace{-72pt}
Q^{(3)}(7) =  [\hspace{1.5pt}[0,3634],[1,6130],[2,2800],[3,1185],[4,35],[6,200],[8,65],[11,20], [19,10],[24,5]\hspace{1.5pt}]
,
\nonumber
\\[3pt]
&&
\hspace{-72pt}
Q^{(4)}(7) =  [\hspace{1.5pt}[0,1415],[1,1545],[2,405],[3,149],[6,15],[8,6]\hspace{1.5pt}]
,
\nonumber
\\[3pt]
&&
\hspace{-72pt}
Q^{(0)}(8) =  [\hspace{1.5pt}[1,13700],[2,33375],[3,32638],[4,14084],[5,3535],[6,4919],[7,35],[8,327],
\nonumber
\\[3pt]
&&
\hspace{-22pt}
[11,3634],[19,669],[20,1415],[24,65],[37,280],[46,470],[70,27],[73,116],
\nonumber
\\[3pt]
&&
\hspace{-22pt}
[89,16],[117,145],[135,1],[236,65],[312,20],[350,23],[415,5],[807,14],[863,1],
\nonumber
\\[3pt]
&&
\hspace{-22pt}
[1459,9],[2046,5],[2372,2],[2878,1],[6352,1],[10527,1],[14115,1],[16072,2]\hspace{1.5pt}]
,
\nonumber
\\[3pt]
&&
\hspace{-72pt}
Q^{(1)}(8) =  [\hspace{1.5pt}[0,13700],[1,68707],[2,92043],[3,65520],[4,17675],[5,2940],[6,11147],
\nonumber
\\[3pt]
&&
\hspace{-22pt}
[7,8],[8,1403],[11,5660],[19,1410],[20,1400],[24,248],[37,162],[46,580],
\nonumber
\\[3pt]
&&
\hspace{-22pt}
[70,7],[73,195],[89,51],[117,100],[236,56],[312,6],[350,27],[415,12],
\nonumber
\\[3pt]
&&
\hspace{-22pt}
[807,5],[1459,4],[2046,3],[2372,3]\hspace{1.5pt}]
,
\nonumber
\\[3pt]
&&
\hspace{-72pt}
Q^{(2)}(8) =  [\hspace{1.5pt}[0,33375],[1,96472],[2,84504],[3,46921],[4,7350],[5,735],[6,8518],
\nonumber
\\[3pt]
&&
\hspace{-22pt}
[8,1526],[11,2800],[19,870],[20,405],[24,218],[37,21],[46,200],[73,84],
\nonumber
\\[3pt]
&&
\hspace{-22pt}
[89,32],[117,15],[236,10],[350,6],[415,4]\hspace{1.5pt}]
,
\nonumber
\\[3pt]
&&
\hspace{-72pt}
Q^{(3)}(8) =  [\hspace{1.5pt}[0,30681],[1,59970],[2,35350],[3,15930],[4,1225],[5,56],[6,2800],[8,645],
\nonumber
\\[3pt]
&&
\hspace{-22pt}
[11,540],[19,200],[20,35],[24,65],[46,20],[73,10],[89,5]\hspace{1.5pt}]
,
\nonumber
\\[3pt]
&&
\hspace{-72pt}
Q^{(4)}(8) =  [\hspace{1.5pt}[0,14084],[1,19090],[2,7350],[3,2770],[4,70],[6,405],[8,114],[11,35],
\nonumber
\\[3pt]
&&
\hspace{-22pt}
[19,15],[24,6]\hspace{1.5pt}]
,
\nonumber
\\[3pt]
&&
\hspace{-72pt}
Q^{(5)}(8) =  [\hspace{1.5pt}[0,3535],[1,3220],[2,735],[3,238],[6,21],[8,7]\hspace{1.5pt}].
\label{a20}
\eea

\section{Characteristic polynomials of the Laplacian matrix in the DSC$(1)$ model}
\label{sa2}

The characteristic polynomial of the Laplacian matrix of the DSC$(1)$ tree factorizes into the product of $\lambda$ and the polynomials $\pi_i(\lambda)$, see equation~\eqref{Pi-factor}. In this appendix, we present the polynomials $\pi_i(\lambda)$ for $i\leq 7$ directly computed by Mathematica for trees  generated by the DSC$(1)$ model. As in independent check, we derived $\pi_i(\lambda)$ from the recursive relations~\eq{585} and \eq{pi-rec}, confirming their validity. We also derived $\pi_8(\lambda)$ and $\pi_9(\lambda)$. The latter is the polynomial of $\lambda$ of degree $2^8=256$, viz. $\pi_9(\lambda)=\sum_{0\leq j\leq 256}a_9^{(j)}\lambda^j$. Most of the coefficients $a_9^{(j)}$ are huge integers, greatly exceeding the coefficients $a_7^{(j)}$ of the polynomial $\pi_7(\lambda)$ which are also huge, cf. \eqref{pi-7}. The expressions for $\pi_8(\lambda)$ and $\pi_9(\lambda)$ are too lengthy to include here. 

The polynomials $\pi_i(\lambda)$ for $i\leq 6$ read
\bea
\label{pi:1-5}
\pi_1 &=& 2 - \lambda 
,
\label{b10}
\\[3pt]
\pi_2 &=& 2 - 4 \lambda + \lambda^2
,
\label{b20}
\\[3pt]
\pi_3 &=& 2 -12 \lambda + 18 \lambda^2 - 8 \lambda^3 + \lambda^4
,
\label{b30}
\\[3pt]
\pi_4 &=& 2 - 32 \lambda + 168 \lambda^2 - 396 \lambda^3 + 472 \lambda^4 - 296 \lambda^5 
           + 98 \lambda^6 - 16 \lambda^7 +  \lambda^8
,
\label{b40}
\\[3pt]
\pi_5 &=& 2 - 80 \lambda + 1208 \lambda^2 - 9552 \lambda^3 + 45476 \lambda^4
              - 140600 \lambda^5 + 295460 \lambda^6 - 434172 \lambda^7  + 453962 \lambda^8 
\nonumber
\\[3pt]
&&
              - 340944 \lambda^9+ 184428 \lambda^{10} - 71516 \lambda^{11} + 19594 \lambda^{12} 
              - 3684 \lambda^{13} + 450 \lambda^{14} - 32 \lambda^{15} +  \lambda^{16} 
,
\label{pi:1-5}
\\[3pt]
\pi_6 &=&2 - 192 \lambda + 7552 \lambda^2 - 168016 \lambda^3 + 2429056 \lambda^4 - 
   24672704 \lambda^5 + 185087064 \lambda^6   - 1061439008 \lambda^7 
\nonumber
\\[3pt]
&&
 + 4771422936 \lambda^8 -  17131717600 \lambda^9 + 49845993112 \lambda^{10}  - 118859071000 \lambda^{11} + 
   234345967208 \lambda^{12} 
\nonumber
\\[3pt]
&&
  - 384712533224 \lambda^{13} + 528733119200 \lambda^{14} -610901653388 \lambda^{15} + 595213158568 \lambda^{16} 
\nonumber
\\[3pt]
&&
   - 490043745744 \lambda^{17} +  341299300448 \lambda^{18}    - 201108575904 \lambda^{19} + 100168324568 \lambda^{20} 
\nonumber
\\[3pt]
&&
-  42087010696 \lambda^{21} + 14866893448 \lambda^{22}  - 4393269256 \lambda^{23} + 
   1078500576 \lambda^{24}
\nonumber
\\[3pt]
&&
   - 217847304 \lambda^{25} + 35733280 \lambda^{26} - 4673696 \lambda^{27} + 474920 \lambda^{28} - 36064 \lambda^{29} + 1922 \lambda^{30}  - 64 \lambda^{31} + \lambda^{32}
.
\label{pi-6}
\eea
The largest polynomial which we present here is $\pi_7$ of degree $2^6=64$: 
\bea
\pi_7 &=& 2 -  448 \lambda + 43296 \lambda^2 - 2475008 \lambda^3 + 95802352 \lambda^4 - 2710026912 \lambda^5 + 58848566672 \lambda^6
\nonumber
\\[3pt]
&&
    - 1015322599120 \lambda^7 + 14276120264664 \lambda^8 - 166811577840832 \lambda^9 + 1644952088329168 \lambda^{10} 
\nonumber
\\[3pt]
&&
    - 13861672325942416 \lambda^{11} + 100852187226166856 \lambda^{12}  - 639005480497405168 \lambda^{13} 
\nonumber
\\[3pt]
&&
   + 3551759016685178600 \lambda^{14} - 
   17426681983367704176 \lambda^{15} 
\nonumber
\\[3pt]
&&
   + 75885537057453902356 \lambda^{16} -   294652679333249886624 \lambda^{17} 
\nonumber
\\[3pt]
&&
   + 1024345335985165113984 \lambda^{18} - 
   3199835643260413205840 \lambda^{19} 
\nonumber
\\[3pt]
&&
   + 9010097606614633096736 \lambda^{20} - 
   22933328467668082847696 \lambda^{21} 
\nonumber
\\[3pt]
&&
   + 52895197062371681474728 \lambda^{22} - 
   110797548329198863794632 \lambda^{23} 
\nonumber
\\[3pt]
&&
   + 211180285681193265670336 \lambda^{24} - 
   366887826245295683140832  \lambda^{25} 
\nonumber
\\[3pt]
&&
   + 581875029694786789717496  \lambda^{26} - 
   843573067183096269924688  \lambda^{27} 
\nonumber
\\[3pt]
&&
   + 1119229472188587709837192  \lambda^{28} - 
   1360373049232441866421904  \lambda^{29} 
\nonumber
\\[3pt]
&&
   + 1516054243947087449817516  \lambda^{30} - 
   1550259517961348141265516  \lambda^{31} 
\nonumber
\\[3pt]
&&
   + 1455406505514364268533698  \lambda^{32} - 
   1255034375536319939236704  \lambda^{33} 
\nonumber
\\[3pt]
&&
   + 994407288655882984773456  \lambda^{34} - 
   724104605969591122083056  \lambda^{35} 
\nonumber
\\[3pt]
&&
   + 484620861293428811828944  \lambda^{36} - 
   298088604338866354974016  \lambda^{37} 
\nonumber
\\[3pt]
&&
   + 168481042772369799801968  \lambda^{38} - 
   87473693463505454653200  \lambda^{39} 
\nonumber
\\[3pt]
&&
   + 41698241195897687997968  \lambda^{40} - 
   18238748013419543763968  \lambda^{41} 
\nonumber
\\[3pt]
&&
   + 7314030250954793877232  \lambda^{42} - 
   2686380322408351452648 \lambda^{43} 
\nonumber
\\[3pt]
&&
   + 902610279766402908552 \lambda^{44} - 
   277029440260862642776 \lambda^{45} 
\nonumber
\\[3pt]
&&
   + 77535414566037487188 \lambda^{46} - 
   19749102387753234676 \lambda^{47}   
   + 4567107690848518688 \lambda^{48}  
\nonumber
\\[3pt]
&&
 - 956268175503207552 \lambda^{49}     
   + 180698152592272704 \lambda^{50} - 
   30697807426865216 \lambda^{51} 
\nonumber
\\[3pt]
&&
   + 4667515749947904 \lambda^{52} - 
   631780861865648 \lambda^{53} + 75642564296104 \lambda^{54} 
\nonumber
\\[3pt]
&&
   - 7949007809280 \lambda^{55} + 
   726210225040 \lambda^{56} - 56993981992 \lambda^{57} + 3784179976 \lambda^{58} 
\nonumber
\\[3pt]
&&
   - 208316232 \lambda^{59} + 9247506 \lambda^{60} - 317916 \lambda^{61} + 7938 \lambda^{62} - 
   128 \lambda^{63} + \lambda^{64}
.
\label{pi-7}
\eea

\end{widetext}

%%\textcolor{red}{xxxxxxxxxxx}

%%%%%%%%%%%%%%%%%%%%%%%%%%%%%%%%%%%%%%\section*{References}

%%\bibliographystyle{longbibliography}
%%%%%%%%%%%%%%%%%%%%%%%%%%%%%%%%%%%%%%%%%\bibliographystyle{iopart-num}
\bibliography{deterministic}

\begin{thebibliography}{49}%
\makeatletter
\providecommand \@ifxundefined [1]{%
 \@ifx{#1\undefined}
}%
\providecommand \@ifnum [1]{%
 \ifnum #1\expandafter \@firstoftwo
 \else \expandafter \@secondoftwo
 \fi
}%
\providecommand \@ifx [1]{%
 \ifx #1\expandafter \@firstoftwo
 \else \expandafter \@secondoftwo
 \fi
}%
\providecommand \natexlab [1]{#1}%
\providecommand \enquote  [1]{``#1''}%
\providecommand \bibnamefont  [1]{#1}%
\providecommand \bibfnamefont [1]{#1}%
\providecommand \citenamefont [1]{#1}%
\providecommand \href@noop [0]{\@secondoftwo}%
\providecommand \href [0]{\begingroup \@sanitize@url \@href}%
\providecommand \@href[1]{\@@startlink{#1}\@@href}%
\providecommand \@@href[1]{\endgroup#1\@@endlink}%
\providecommand \@sanitize@url [0]{\catcode `\\12\catcode `\$12\catcode
  `\&12\catcode `\#12\catcode `\^12\catcode `\_12\catcode `\%12\relax}%
\providecommand \@@startlink[1]{}%
\providecommand \@@endlink[0]{}%
\providecommand \url  [0]{\begingroup\@sanitize@url \@url }%
\providecommand \@url [1]{\endgroup\@href {#1}{\urlprefix }}%
\providecommand \urlprefix  [0]{URL }%
\providecommand \Eprint [0]{\href }%
\providecommand \doibase [0]{http://dx.doi.org/}%
\providecommand \selectlanguage [0]{\@gobble}%
\providecommand \bibinfo  [0]{\@secondoftwo}%
\providecommand \bibfield  [0]{\@secondoftwo}%
\providecommand \translation [1]{[#1]}%
\providecommand \BibitemOpen [0]{}%
\providecommand \bibitemStop [0]{}%
\providecommand \bibitemNoStop [0]{.\EOS\space}%
\providecommand \EOS [0]{\spacefactor3000\relax}%
\providecommand \BibitemShut  [1]{\csname bibitem#1\endcsname}%
\let\auto@bib@innerbib\@empty
%</preamble>
\bibitem [{\citenamefont {Hatcher}(2002)}]{Hatcher}%
  \BibitemOpen
  \bibfield  {author} {\bibinfo {author} {\bibfnamefont {A.}~\bibnamefont
  {Hatcher}},\ }\href@noop {} {\emph {\bibinfo {title} {Algebraic Topology}}}\
  (\bibinfo  {publisher} {Cambridge University Press},\ \bibinfo {address}
  {Cambridge},\ \bibinfo {year} {2002})\BibitemShut {NoStop}%
\bibitem [{\citenamefont {Pippenger}\ and\ \citenamefont
  {Schleich}(2006)}]{Pippenger}%
  \BibitemOpen
  \bibfield  {author} {\bibinfo {author} {\bibfnamefont {N.}~\bibnamefont
  {Pippenger}}\ and\ \bibinfo {author} {\bibfnamefont {K.}~\bibnamefont
  {Schleich}},\ }\bibfield  {title} {\enquote {\bibinfo {title} {Topological
  characteristics of random triangulated surfaces},}\ }\href {\doibase
  https://doi.org/10.1002/rsa.20080} {\bibfield  {journal} {\bibinfo  {journal}
  {Random Struct. Algorithms}\ }\textbf {\bibinfo {volume} {28}},\ \bibinfo
  {pages} {247--288} (\bibinfo {year} {2006})}\BibitemShut {NoStop}%
\bibitem [{\citenamefont {Linial}\ and\ \citenamefont
  {Peled}(2016)}]{Linial16}%
  \BibitemOpen
  \bibfield  {author} {\bibinfo {author} {\bibfnamefont {N.}~\bibnamefont
  {Linial}}\ and\ \bibinfo {author} {\bibfnamefont {Y.}~\bibnamefont {Peled}},\
  }\bibfield  {title} {\enquote {\bibinfo {title} {On the phase transition in
  random simplicial complexes},}\ }\href {http://www.jstor.org/stable/44072029}
  {\bibfield  {journal} {\bibinfo  {journal} {Ann. Math.}\ }\textbf {\bibinfo
  {volume} {184}},\ \bibinfo {pages} {745--773} (\bibinfo {year}
  {2016})}\BibitemShut {NoStop}%
\bibitem [{\citenamefont {Costa}\ and\ \citenamefont
  {Farber}(2016)}]{Farber16}%
  \BibitemOpen
  \bibfield  {author} {\bibinfo {author} {\bibfnamefont {A.}~\bibnamefont
  {Costa}}\ and\ \bibinfo {author} {\bibfnamefont {M.}~\bibnamefont {Farber}},\
  }\bibfield  {title} {\enquote {\bibinfo {title} {Large random simplicial
  complexes, {I}},}\ }\href {\doibase 10.1142/S179352531650014X} {\bibfield
  {journal} {\bibinfo  {journal} {J. Topol. Anal.}\ }\textbf {\bibinfo {volume}
  {8}},\ \bibinfo {pages} {399--429} (\bibinfo {year} {2016})}\BibitemShut
  {NoStop}%
\bibitem [{\citenamefont {Kahle}(2017)}]{Kahle}%
  \BibitemOpen
  \bibfield  {author} {\bibinfo {author} {\bibfnamefont {M.}~\bibnamefont
  {Kahle}},\ }\bibfield  {title} {\enquote {\bibinfo {title} {Random simplicial
  complexes},}\ }in\ \href {https://doi.org/10.1201/9781315119601} {\emph
  {\bibinfo {booktitle} {Handbook of Discrete and Computational Geometry}}},\
  \bibinfo {editor} {edited by\ \bibinfo {editor} {\bibfnamefont {C.~D.}\
  \bibnamefont {Toth}}, \bibinfo {editor} {\bibfnamefont {J.}~\bibnamefont
  {O'{R}ourke}}, \ and\ \bibinfo {editor} {\bibfnamefont {J.~E.}\ \bibnamefont
  {Goodman}}}\ (\bibinfo  {publisher} {Chapman and Hall/CRC},\ \bibinfo
  {address} {New York},\ \bibinfo {year} {2017})\BibitemShut {NoStop}%
\bibitem [{\citenamefont {Battiston}\ \emph {et~al.}(2020)\citenamefont
  {Battiston}, \citenamefont {Cencetti}, \citenamefont {Iacopini},
  \citenamefont {Latora}, \citenamefont {Lucas}, \citenamefont {Patania},
  \citenamefont {Young},\ and\ \citenamefont {Petri}}]{Petri20}%
  \BibitemOpen
  \bibfield  {author} {\bibinfo {author} {\bibfnamefont {Federico}\
  \bibnamefont {Battiston}}, \bibinfo {author} {\bibfnamefont {Giulia}\
  \bibnamefont {Cencetti}}, \bibinfo {author} {\bibfnamefont {Iacopo}\
  \bibnamefont {Iacopini}}, \bibinfo {author} {\bibfnamefont {Vito}\
  \bibnamefont {Latora}}, \bibinfo {author} {\bibfnamefont {Maxime}\
  \bibnamefont {Lucas}}, \bibinfo {author} {\bibfnamefont {Alice}\ \bibnamefont
  {Patania}}, \bibinfo {author} {\bibfnamefont {Jean-Gabriel}\ \bibnamefont
  {Young}}, \ and\ \bibinfo {author} {\bibfnamefont {Giovanni}\ \bibnamefont
  {Petri}},\ }\bibfield  {title} {\enquote {\bibinfo {title} {Networks beyond
  pairwise interactions: Structure and dynamics},}\ }\href {\doibase
  https://doi.org/10.1016/j.physrep.2020.05.004} {\bibfield  {journal}
  {\bibinfo  {journal} {Phys. Reports}\ }\textbf {\bibinfo {volume} {874}},\
  \bibinfo {pages} {1--92} (\bibinfo {year} {2020})},\ \bibinfo {note}
  {networks beyond pairwise interactions: Structure and dynamics}\BibitemShut
  {NoStop}%
\bibitem [{\citenamefont {Bianconi}(2021)}]{Ginestra21}%
  \BibitemOpen
  \bibfield  {author} {\bibinfo {author} {\bibfnamefont {G.}~\bibnamefont
  {Bianconi}},\ }\href {https://doi.org/10.1017/9781108770996} {\emph {\bibinfo
  {title} {Higher-Order Networks}}}\ (\bibinfo  {publisher} {Cambridge
  University Press},\ \bibinfo {address} {Cambridge, UK},\ \bibinfo {year}
  {2021})\BibitemShut {NoStop}%
\bibitem [{\citenamefont {Bobrowski}\ and\ \citenamefont
  {Krioukov}(2022)}]{bobrowski2022random}%
  \BibitemOpen
  \bibfield  {author} {\bibinfo {author} {\bibfnamefont {O.}~\bibnamefont
  {Bobrowski}}\ and\ \bibinfo {author} {\bibfnamefont {D.}~\bibnamefont
  {Krioukov}},\ }\bibfield  {title} {\enquote {\bibinfo {title} {{Random
  simplicial complexes: Models and phenomena}},}\ }in\ \href
  {https://doi.org/10.1007/978-3-030-91374-8_2} {\emph {\bibinfo {booktitle}
  {Higher-Order Systems. Understanding Complex Systems}}},\ \bibinfo {editor}
  {edited by\ \bibinfo {editor} {\bibfnamefont {F.}~\bibnamefont {Battiston}}\
  and\ \bibinfo {editor} {\bibfnamefont {G.}~\bibnamefont {Petri}}}\ (\bibinfo
  {publisher} {Springer},\ \bibinfo {address} {Cham, CH},\ \bibinfo {year}
  {2022})\ pp.\ \bibinfo {pages} {59--96}\BibitemShut {NoStop}%
\bibitem [{\citenamefont {Bianconi}\ \emph {et~al.}(2015)\citenamefont
  {Bianconi}, \citenamefont {Rahmede},\ and\ \citenamefont
  {Wu}}]{bianconi2015complex}%
  \BibitemOpen
  \bibfield  {author} {\bibinfo {author} {\bibfnamefont {G.}~\bibnamefont
  {Bianconi}}, \bibinfo {author} {\bibfnamefont {C.}~\bibnamefont {Rahmede}}, \
  and\ \bibinfo {author} {\bibfnamefont {Z.}~\bibnamefont {Wu}},\ }\bibfield
  {title} {\enquote {\bibinfo {title} {Complex quantum network geometries:
  Evolution and phase transitions},}\ }\href {\doibase
  10.1103/PhysRevE.92.022815} {\bibfield  {journal} {\bibinfo  {journal} {Phys.
  Rev. E}\ }\textbf {\bibinfo {volume} {92}},\ \bibinfo {pages} {022815}
  (\bibinfo {year} {2015})}\BibitemShut {NoStop}%
\bibitem [{\citenamefont {Bianconi}\ and\ \citenamefont
  {Rahmede}(2017)}]{bianconi2017emergent}%
  \BibitemOpen
  \bibfield  {author} {\bibinfo {author} {\bibfnamefont {G.}~\bibnamefont
  {Bianconi}}\ and\ \bibinfo {author} {\bibfnamefont {C.}~\bibnamefont
  {Rahmede}},\ }\bibfield  {title} {\enquote {\bibinfo {title} {Emergent
  hyperbolic network geometry},}\ }\href {https://doi.org/10.1038/srep41974}
  {\bibfield  {journal} {\bibinfo  {journal} {Sci. Rep.}\ }\textbf {\bibinfo
  {volume} {7}},\ \bibinfo {pages} {41974} (\bibinfo {year}
  {2017})}\BibitemShut {NoStop}%
\bibitem [{\citenamefont {da~Silva}\ \emph {et~al.}(2018)\citenamefont
  {da~Silva}, \citenamefont {Bianconi}, \citenamefont {da~Costa}, \citenamefont
  {Dorogovtsev},\ and\ \citenamefont {Mendes}}]{da2018complex}%
  \BibitemOpen
  \bibfield  {author} {\bibinfo {author} {\bibfnamefont {D.~C.}\ \bibnamefont
  {da~Silva}}, \bibinfo {author} {\bibfnamefont {G.}~\bibnamefont {Bianconi}},
  \bibinfo {author} {\bibfnamefont {R.~A.}\ \bibnamefont {da~Costa}}, \bibinfo
  {author} {\bibfnamefont {S.~N.}\ \bibnamefont {Dorogovtsev}}, \ and\ \bibinfo
  {author} {\bibfnamefont {J.~F.~F.}\ \bibnamefont {Mendes}},\ }\bibfield
  {title} {\enquote {\bibinfo {title} {Complex network view of evolving
  manifolds},}\ }\href {\doibase 10.1103/PhysRevE.97.032316} {\bibfield
  {journal} {\bibinfo  {journal} {Phys. Rev. E}\ }\textbf {\bibinfo {volume}
  {97}},\ \bibinfo {pages} {032316} (\bibinfo {year} {2018})}\BibitemShut
  {NoStop}%
\bibitem [{\citenamefont {Xie}\ \emph {et~al.}(2023)\citenamefont {Xie},
  \citenamefont {Wang}, \citenamefont {Xu}, \citenamefont {Zhu}, \citenamefont
  {Li},\ and\ \citenamefont {Zhang}}]{xie2023combinatorial}%
  \BibitemOpen
  \bibfield  {author} {\bibinfo {author} {\bibfnamefont {Z.}~\bibnamefont
  {Xie}}, \bibinfo {author} {\bibfnamefont {Y.}~\bibnamefont {Wang}}, \bibinfo
  {author} {\bibfnamefont {W.}~\bibnamefont {Xu}}, \bibinfo {author}
  {\bibfnamefont {L.}~\bibnamefont {Zhu}}, \bibinfo {author} {\bibfnamefont
  {W.}~\bibnamefont {Li}}, \ and\ \bibinfo {author} {\bibfnamefont
  {Z.}~\bibnamefont {Zhang}},\ }\bibfield  {title} {\enquote {\bibinfo {title}
  {Combinatorial properties for a class of simplicial complexes extended from
  pseudo-fractal scale-free web},}\ }\href
  {https://doi.org/10.1142/S0218348X23500226} {\bibfield  {journal} {\bibinfo
  {journal} {Fractals}\ }\textbf {\bibinfo {volume} {31}},\ \bibinfo {pages}
  {2350022} (\bibinfo {year} {2023})}\BibitemShut {NoStop}%
\bibitem [{\citenamefont {Barab{\'a}si}\ \emph {et~al.}(2001)\citenamefont
  {Barab{\'a}si}, \citenamefont {Ravasz},\ and\ \citenamefont
  {Vicsek}}]{barabasi2001deterministic}%
  \BibitemOpen
  \bibfield  {author} {\bibinfo {author} {\bibfnamefont {A.-L.}\ \bibnamefont
  {Barab{\'a}si}}, \bibinfo {author} {\bibfnamefont {E.}~\bibnamefont
  {Ravasz}}, \ and\ \bibinfo {author} {\bibfnamefont {T.}~\bibnamefont
  {Vicsek}},\ }\bibfield  {title} {\enquote {\bibinfo {title} {Deterministic
  scale-free networks},}\ }\href {\doibase
  https://doi.org/10.1016/S0378-4371(01)00369-7} {\bibfield  {journal}
  {\bibinfo  {journal} {Physica A}\ }\textbf {\bibinfo {volume} {299}},\
  \bibinfo {pages} {559} (\bibinfo {year} {2001})}\BibitemShut {NoStop}%
\bibitem [{\citenamefont {Dorogovtsev}\ \emph {et~al.}(2002)\citenamefont
  {Dorogovtsev}, \citenamefont {Goltsev},\ and\ \citenamefont
  {Mendes}}]{Sergey02}%
  \BibitemOpen
  \bibfield  {author} {\bibinfo {author} {\bibfnamefont {S.~N.}\ \bibnamefont
  {Dorogovtsev}}, \bibinfo {author} {\bibfnamefont {A.~V.}\ \bibnamefont
  {Goltsev}}, \ and\ \bibinfo {author} {\bibfnamefont {J.~F.~F.}\ \bibnamefont
  {Mendes}},\ }\bibfield  {title} {\enquote {\bibinfo {title} {Pseudofractal
  scale-free web},}\ }\href {\doibase 10.1103/PhysRevE.65.066122} {\bibfield
  {journal} {\bibinfo  {journal} {Phys. Rev. E}\ }\textbf {\bibinfo {volume}
  {65}},\ \bibinfo {pages} {066122} (\bibinfo {year} {2002})}\BibitemShut
  {NoStop}%
\bibitem [{\citenamefont {Andrade~Jr}\ \emph {et~al.}(2005)\citenamefont
  {Andrade~Jr}, \citenamefont {Herrmann}, \citenamefont {Andrade},\ and\
  \citenamefont {Da~Silva}}]{andrade2005apollonian}%
  \BibitemOpen
  \bibfield  {author} {\bibinfo {author} {\bibfnamefont {J.~S.}\ \bibnamefont
  {Andrade~Jr}}, \bibinfo {author} {\bibfnamefont {H.~J.}\ \bibnamefont
  {Herrmann}}, \bibinfo {author} {\bibfnamefont {R.~F.~S.}\ \bibnamefont
  {Andrade}}, \ and\ \bibinfo {author} {\bibfnamefont {L.~R.}\ \bibnamefont
  {Da~Silva}},\ }\bibfield  {title} {\enquote {\bibinfo {title} {Apollonian
  networks: Simultaneously scale-free, small world, euclidean, space filling,
  and with matching graphs},}\ }\href {\doibase 10.1103/PhysRevLett.94.018702}
  {\bibfield  {journal} {\bibinfo  {journal} {Phys. Rev. Lett.}\ }\textbf
  {\bibinfo {volume} {94}},\ \bibinfo {pages} {018702} (\bibinfo {year}
  {2005})}\BibitemShut {NoStop}%
\bibitem [{\citenamefont {Hwang}\ \emph {et~al.}(2010)\citenamefont {Hwang},
  \citenamefont {Yun}, \citenamefont {Lee}, \citenamefont {Kahng},\ and\
  \citenamefont {Kim}}]{hwang2010spectral}%
  \BibitemOpen
  \bibfield  {author} {\bibinfo {author} {\bibfnamefont {S.}~\bibnamefont
  {Hwang}}, \bibinfo {author} {\bibfnamefont {C.-K.}\ \bibnamefont {Yun}},
  \bibinfo {author} {\bibfnamefont {D.-S.}\ \bibnamefont {Lee}}, \bibinfo
  {author} {\bibfnamefont {B.}~\bibnamefont {Kahng}}, \ and\ \bibinfo {author}
  {\bibfnamefont {D.}~\bibnamefont {Kim}},\ }\bibfield  {title} {\enquote
  {\bibinfo {title} {Spectral dimensions of hierarchical scale-free networks
  with weighted shortcuts},}\ }\href {\doibase 10.1103/PhysRevE.82.056110}
  {\bibfield  {journal} {\bibinfo  {journal} {Phys. Rev. E}\ }\textbf {\bibinfo
  {volume} {82}},\ \bibinfo {pages} {056110} (\bibinfo {year}
  {2010})}\BibitemShut {NoStop}%
\bibitem [{\citenamefont {Hasegawa}\ and\ \citenamefont
  {Nemoto}(2013)}]{hasegawa2013hierarchical}%
  \BibitemOpen
  \bibfield  {author} {\bibinfo {author} {\bibfnamefont {T.}~\bibnamefont
  {Hasegawa}}\ and\ \bibinfo {author} {\bibfnamefont {K.}~\bibnamefont
  {Nemoto}},\ }\bibfield  {title} {\enquote {\bibinfo {title} {Hierarchical
  scale-free network is fragile against random failure},}\ }\href {\doibase
  10.1103/PhysRevE.88.062807} {\bibfield  {journal} {\bibinfo  {journal} {Phys.
  Rev. E}\ }\textbf {\bibinfo {volume} {88}},\ \bibinfo {pages} {062807}
  (\bibinfo {year} {2013})}\BibitemShut {NoStop}%
\bibitem [{\citenamefont {Dorogovtsev}\ and\ \citenamefont
  {Mendes}(2022)}]{dorogovtsev2022the}%
  \BibitemOpen
  \bibfield  {author} {\bibinfo {author} {\bibfnamefont {S.~N.}\ \bibnamefont
  {Dorogovtsev}}\ and\ \bibinfo {author} {\bibfnamefont {J.~F.~F.}\
  \bibnamefont {Mendes}},\ }\href
  {https://doi.org/10.1093/oso/9780199695119.001.0001} {\emph {\bibinfo {title}
  {The Nature of Complex Networks}}}\ (\bibinfo  {publisher} {Oxford University
  Press},\ \bibinfo {address} {Oxford},\ \bibinfo {year} {2022})\BibitemShut
  {NoStop}%
\bibitem [{\citenamefont {Dorogovtsev}\ \emph {et~al.}(2001)\citenamefont
  {Dorogovtsev}, \citenamefont {Mendes},\ and\ \citenamefont
  {Samukhin}}]{dorogovtsev2001size}%
  \BibitemOpen
  \bibfield  {author} {\bibinfo {author} {\bibfnamefont {S.~N.}\ \bibnamefont
  {Dorogovtsev}}, \bibinfo {author} {\bibfnamefont {J.~F.~F.}\ \bibnamefont
  {Mendes}}, \ and\ \bibinfo {author} {\bibfnamefont {A.~N.}\ \bibnamefont
  {Samukhin}},\ }\bibfield  {title} {\enquote {\bibinfo {title} {Size-dependent
  degree distribution of a scale-free growing network},}\ }\href {\doibase
  10.1103/PhysRevE.63.062101} {\bibfield  {journal} {\bibinfo  {journal} {Phys.
  Rev. E}\ }\textbf {\bibinfo {volume} {63}},\ \bibinfo {pages} {062101}
  (\bibinfo {year} {2001})}\BibitemShut {NoStop}%
\bibitem [{\citenamefont {Bianconi}\ and\ \citenamefont
  {Ziff}(2018)}]{bianconi2018topological}%
  \BibitemOpen
  \bibfield  {author} {\bibinfo {author} {\bibfnamefont {G.}~\bibnamefont
  {Bianconi}}\ and\ \bibinfo {author} {\bibfnamefont {R.~M.}\ \bibnamefont
  {Ziff}},\ }\bibfield  {title} {\enquote {\bibinfo {title} {Topological
  percolation on hyperbolic simplicial complexes},}\ }\href {\doibase
  10.1103/PhysRevE.98.052308} {\bibfield  {journal} {\bibinfo  {journal} {Phys.
  Rev. E}\ }\textbf {\bibinfo {volume} {98}},\ \bibinfo {pages} {052308}
  (\bibinfo {year} {2018})}\BibitemShut {NoStop}%
\bibitem [{\citenamefont {Bianconi}\ and\ \citenamefont
  {Dorogovtsev}(2020)}]{bianconi2020spectral}%
  \BibitemOpen
  \bibfield  {author} {\bibinfo {author} {\bibfnamefont {G.}~\bibnamefont
  {Bianconi}}\ and\ \bibinfo {author} {\bibfnamefont {S.~N.}\ \bibnamefont
  {Dorogovtsev}},\ }\bibfield  {title} {\enquote {\bibinfo {title} {{The
  spectral dimension of simplicial complexes: A renormalization group
  theory}},}\ }\href {\doibase 10.1088/1742-5468/ab5d0e} {\bibfield  {journal}
  {\bibinfo  {journal} {J. Stat. Mech.}\ }\textbf {\bibinfo {volume} {2020}},\
  \bibinfo {pages} {014005} (\bibinfo {year} {2020})}\BibitemShut {NoStop}%
\bibitem [{\citenamefont {Reitz}\ and\ \citenamefont
  {Bianconi}(2020)}]{reitz2020higher}%
  \BibitemOpen
  \bibfield  {author} {\bibinfo {author} {\bibfnamefont {M.}~\bibnamefont
  {Reitz}}\ and\ \bibinfo {author} {\bibfnamefont {G.}~\bibnamefont
  {Bianconi}},\ }\bibfield  {title} {\enquote {\bibinfo {title} {{The
  higher-order spectrum of simplicial complexes: A renormalization group
  approach}},}\ }\href {\doibase 10.1088/1751-8121/ab9338} {\bibfield
  {journal} {\bibinfo  {journal} {J. Phys. A: Math. Theor.}\ }\textbf {\bibinfo
  {volume} {53}},\ \bibinfo {pages} {295001} (\bibinfo {year}
  {2020})}\BibitemShut {NoStop}%
\bibitem [{\citenamefont {Hodge}(1934)}]{hodge1934dirichlet}%
  \BibitemOpen
  \bibfield  {author} {\bibinfo {author} {\bibfnamefont {W.~V.~D.}\
  \bibnamefont {Hodge}},\ }\bibfield  {title} {\enquote {\bibinfo {title} {A
  {D}irichlet problem for harmonic functionals, with applications to analytic
  varities},}\ }\href@noop {} {\bibfield  {journal} {\bibinfo  {journal} {Proc.
  Lond. Math. Soc.}\ }\textbf {\bibinfo {volume} {2}},\ \bibinfo {pages}
  {257--303} (\bibinfo {year} {1934})}\BibitemShut {NoStop}%
\bibitem [{\citenamefont {Dodziuk}(1976)}]{Dodziuk-hodge}%
  \BibitemOpen
  \bibfield  {author} {\bibinfo {author} {\bibfnamefont {J.}~\bibnamefont
  {Dodziuk}},\ }\bibfield  {title} {\enquote {\bibinfo {title}
  {Finite-difference approach to the {H}odge theory of harmonic forms},}\
  }\href {https://doi.org/10.2307/2373615} {\bibfield  {journal} {\bibinfo
  {journal} {Amer. J. Math.}\ }\textbf {\bibinfo {volume} {98}},\ \bibinfo
  {pages} {79--104} (\bibinfo {year} {1976})}\BibitemShut {NoStop}%
\bibitem [{\citenamefont {Lim}(2020)}]{lim2020hodge}%
  \BibitemOpen
  \bibfield  {author} {\bibinfo {author} {\bibfnamefont {L.-H.}\ \bibnamefont
  {Lim}},\ }\bibfield  {title} {\enquote {\bibinfo {title} {{Hodge Laplacians
  on graphs}},}\ }\href {\doibase 10.1137/18M1223101} {\bibfield  {journal}
  {\bibinfo  {journal} {SIAM Rev.}\ }\textbf {\bibinfo {volume} {62}},\
  \bibinfo {pages} {685--715} (\bibinfo {year} {2020})}\BibitemShut {NoStop}%
\bibitem [{\citenamefont {Schaub}\ \emph {et~al.}(2020)\citenamefont {Schaub},
  \citenamefont {Benson}, \citenamefont {Horn}, \citenamefont {Lippner},\ and\
  \citenamefont {Jadbabaie}}]{schaub2020random}%
  \BibitemOpen
  \bibfield  {author} {\bibinfo {author} {\bibfnamefont {M.~T.}\ \bibnamefont
  {Schaub}}, \bibinfo {author} {\bibfnamefont {A.~R.}\ \bibnamefont {Benson}},
  \bibinfo {author} {\bibfnamefont {P.}~\bibnamefont {Horn}}, \bibinfo {author}
  {\bibfnamefont {G.}~\bibnamefont {Lippner}}, \ and\ \bibinfo {author}
  {\bibfnamefont {A.}~\bibnamefont {Jadbabaie}},\ }\bibfield  {title} {\enquote
  {\bibinfo {title} {{Random walks on simplicial complexes and the normalized
  Hodge $1$-Laplacian}},}\ }\href {\doibase 10.1137/18M1201019} {\bibfield
  {journal} {\bibinfo  {journal} {SIAM Rev.}\ }\textbf {\bibinfo {volume}
  {62}},\ \bibinfo {pages} {353--391} (\bibinfo {year} {2020})}\BibitemShut
  {NoStop}%
\bibitem [{\citenamefont {Brouwer}\ and\ \citenamefont
  {Haemers}(2011)}]{brouwer2011spectra}%
  \BibitemOpen
  \bibfield  {author} {\bibinfo {author} {\bibfnamefont {A.~E.}\ \bibnamefont
  {Brouwer}}\ and\ \bibinfo {author} {\bibfnamefont {W.~H.}\ \bibnamefont
  {Haemers}},\ }\href {https://doi.org/10.1007/978-1-4614-1939-6} {\emph
  {\bibinfo {title} {Spectra of Graphs}}}\ (\bibinfo  {publisher} {Springer},\
  \bibinfo {address} {New York, NY},\ \bibinfo {year} {2011})\BibitemShut
  {NoStop}%
\bibitem [{\citenamefont {Frankl}(1984)}]{Frankl}%
  \BibitemOpen
  \bibfield  {author} {\bibinfo {author} {\bibfnamefont {P.}~\bibnamefont
  {Frankl}},\ }\bibfield  {title} {\enquote {\bibinfo {title} {A new short
  proof for the {K}ruskal-{K}atona theorem},}\ }\href {\doibase
  https://doi.org/10.1016/0012-365X(84)90193-6} {\bibfield  {journal} {\bibinfo
   {journal} {Discrete Math.}\ }\textbf {\bibinfo {volume} {48}},\ \bibinfo
  {pages} {327} (\bibinfo {year} {1984})}\BibitemShut {NoStop}%
\bibitem [{\citenamefont {Sloane}(2025)}]{Sloane}%
  \BibitemOpen
  \bibfield  {author} {\bibinfo {author} {\bibfnamefont {N.~J.~A.}\
  \bibnamefont {Sloane}},\ }\href@noop {} {\enquote {\bibinfo {title} {{\em The
  Online Encyclopedia of Integer Sequences (OEIS)}},}\ }\bibinfo {howpublished}
  {\url{https://oeis.org/}} (\bibinfo {year} {2025})\BibitemShut {NoStop}%
\bibitem [{\citenamefont {Flajolet}\ and\ \citenamefont
  {Sedgewick}(2009)}]{Flajolet}%
  \BibitemOpen
  \bibfield  {author} {\bibinfo {author} {\bibfnamefont {P.}~\bibnamefont
  {Flajolet}}\ and\ \bibinfo {author} {\bibfnamefont {R.}~\bibnamefont
  {Sedgewick}},\ }\href {https://doi.org/10.1017/CBO9780511801655} {\emph
  {\bibinfo {title} {Analytic Combinatorics}}}\ (\bibinfo  {publisher}
  {Cambridge University Press},\ \bibinfo {address} {Cambridge},\ \bibinfo
  {year} {2009})\BibitemShut {NoStop}%
\bibitem [{\citenamefont {Evgrafov}(2020)}]{evgrafov2020asymptotic}%
  \BibitemOpen
  \bibfield  {author} {\bibinfo {author} {\bibfnamefont {M.~A.}\ \bibnamefont
  {Evgrafov}},\ }\href@noop {} {\emph {\bibinfo {title} {Asymptotic Estimates
  and Entire Functions}}}\ (\bibinfo  {publisher} {Dover Publications},\
  \bibinfo {address} {New York},\ \bibinfo {year} {2020})\BibitemShut {NoStop}%
\bibitem [{\citenamefont {Jung}\ \emph {et~al.}(2002)\citenamefont {Jung},
  \citenamefont {Kim},\ and\ \citenamefont {Kahng}}]{jung2002geometric}%
  \BibitemOpen
  \bibfield  {author} {\bibinfo {author} {\bibfnamefont {S.}~\bibnamefont
  {Jung}}, \bibinfo {author} {\bibfnamefont {S.}~\bibnamefont {Kim}}, \ and\
  \bibinfo {author} {\bibfnamefont {B.}~\bibnamefont {Kahng}},\ }\bibfield
  {title} {\enquote {\bibinfo {title} {Geometric fractal growth model for
  scale-free networks},}\ }\href {\doibase 10.1103/PhysRevE.65.056101}
  {\bibfield  {journal} {\bibinfo  {journal} {Phys. Rev. E}\ }\textbf {\bibinfo
  {volume} {65}},\ \bibinfo {pages} {056101} (\bibinfo {year}
  {2002})}\BibitemShut {NoStop}%
\bibitem [{\citenamefont {Dorogovtsev}\ \emph {et~al.}(2006)\citenamefont
  {Dorogovtsev}, \citenamefont {Mendes},\ and\ \citenamefont
  {Oliveira}}]{dorogovtsev2006degree}%
  \BibitemOpen
  \bibfield  {author} {\bibinfo {author} {\bibfnamefont {S.~N.}\ \bibnamefont
  {Dorogovtsev}}, \bibinfo {author} {\bibfnamefont {J.~F.~F.}\ \bibnamefont
  {Mendes}}, \ and\ \bibinfo {author} {\bibfnamefont {J.~G.}\ \bibnamefont
  {Oliveira}},\ }\bibfield  {title} {\enquote {\bibinfo {title}
  {Degree-dependent intervertex separation in complex networks},}\ }\href
  {\doibase 10.1103/PhysRevE.73.056122} {\bibfield  {journal} {\bibinfo
  {journal} {Phys. Rev. E}\ }\textbf {\bibinfo {volume} {73}},\ \bibinfo
  {pages} {056122} (\bibinfo {year} {2006})}\BibitemShut {NoStop}%
\bibitem [{\citenamefont {Qi}\ \emph {et~al.}(2009)\citenamefont {Qi},
  \citenamefont {Zhang}, \citenamefont {Ding}, \citenamefont {Zhou},\ and\
  \citenamefont {Guan}}]{qi2009structural}%
  \BibitemOpen
  \bibfield  {author} {\bibinfo {author} {\bibfnamefont {Y.}~\bibnamefont
  {Qi}}, \bibinfo {author} {\bibfnamefont {Z.}~\bibnamefont {Zhang}}, \bibinfo
  {author} {\bibfnamefont {B.}~\bibnamefont {Ding}}, \bibinfo {author}
  {\bibfnamefont {S.}~\bibnamefont {Zhou}}, \ and\ \bibinfo {author}
  {\bibfnamefont {J.}~\bibnamefont {Guan}},\ }\bibfield  {title} {\enquote
  {\bibinfo {title} {{Structural and spectral properties of a family of
  deterministic recursive trees: Rigorous solutions}},}\ }\href {\doibase
  10.1088/1751-8113/42/16/165103} {\bibfield  {journal} {\bibinfo  {journal}
  {J. Phys. A: Math. Theor.}\ }\textbf {\bibinfo {volume} {42}},\ \bibinfo
  {pages} {165103} (\bibinfo {year} {2009})}\BibitemShut {NoStop}%
\bibitem [{\citenamefont {Pittel}(1994)}]{Pittel94}%
  \BibitemOpen
  \bibfield  {author} {\bibinfo {author} {\bibfnamefont {B.}~\bibnamefont
  {Pittel}},\ }\bibfield  {title} {\enquote {\bibinfo {title} {Note on the
  heights of random recursive trees and random $m-$ary search trees},}\ }\href
  {https://doi.org/10.1002/rsa.3240050207} {\bibfield  {journal} {\bibinfo
  {journal} {Random Struct. Algorithms}\ }\textbf {\bibinfo {volume} {5}},\
  \bibinfo {pages} {337--347} (\bibinfo {year} {1994})}\BibitemShut {NoStop}%
\bibitem [{\citenamefont {Krapivsky}\ and\ \citenamefont
  {Redner}(2001)}]{KR01}%
  \BibitemOpen
  \bibfield  {author} {\bibinfo {author} {\bibfnamefont {P.~L.}\ \bibnamefont
  {Krapivsky}}\ and\ \bibinfo {author} {\bibfnamefont {S.}~\bibnamefont
  {Redner}},\ }\bibfield  {title} {\enquote {\bibinfo {title} {Organization of
  growing random networks},}\ }\href {\doibase 10.1103/PhysRevE.63.066123}
  {\bibfield  {journal} {\bibinfo  {journal} {Phys. Rev. E}\ }\textbf {\bibinfo
  {volume} {63}},\ \bibinfo {pages} {066123} (\bibinfo {year}
  {2001})}\BibitemShut {NoStop}%
\bibitem [{\citenamefont {Krapivsky}\ and\ \citenamefont
  {Redner}(2002)}]{KR02}%
  \BibitemOpen
  \bibfield  {author} {\bibinfo {author} {\bibfnamefont {P.~L.}\ \bibnamefont
  {Krapivsky}}\ and\ \bibinfo {author} {\bibfnamefont {S.}~\bibnamefont
  {Redner}},\ }\bibfield  {title} {\enquote {\bibinfo {title} {Statistics of
  changes in lead node in connectivity-driven networks},}\ }\href {\doibase
  10.1103/PhysRevLett.89.258703} {\bibfield  {journal} {\bibinfo  {journal}
  {Phys. Rev. Lett.}\ }\textbf {\bibinfo {volume} {89}},\ \bibinfo {pages}
  {258703} (\bibinfo {year} {2002})}\BibitemShut {NoStop}%
\bibitem [{\citenamefont {Janson}(2005)}]{Janson05}%
  \BibitemOpen
  \bibfield  {author} {\bibinfo {author} {\bibfnamefont {S.}~\bibnamefont
  {Janson}},\ }\bibfield  {title} {\enquote {\bibinfo {title} {Asymptotic
  degree distribution in random recursive trees},}\ }\href
  {https://doi.org/10.1002/rsa.20046} {\bibfield  {journal} {\bibinfo
  {journal} {Random Struct. Algorithms}\ }\textbf {\bibinfo {volume} {26}},\
  \bibinfo {pages} {69--83} (\bibinfo {year} {2005})}\BibitemShut {NoStop}%
\bibitem [{\citenamefont {Drmota}(2009)}]{Drmota}%
  \BibitemOpen
  \bibfield  {author} {\bibinfo {author} {\bibfnamefont {M.}~\bibnamefont
  {Drmota}},\ }\href {https://doi.org/10.1007/978-3-211-75357-6} {\emph
  {\bibinfo {title} {Random Trees}}}\ (\bibinfo  {publisher} {Springer},\
  \bibinfo {address} {Vienna},\ \bibinfo {year} {2009})\BibitemShut {NoStop}%
\bibitem [{\citenamefont {Frieze}\ and\ \citenamefont
  {Karo\'{n}ski}(2016)}]{Frieze}%
  \BibitemOpen
  \bibfield  {author} {\bibinfo {author} {\bibfnamefont {A.}~\bibnamefont
  {Frieze}}\ and\ \bibinfo {author} {\bibfnamefont {M.}~\bibnamefont
  {Karo\'{n}ski}},\ }\href {https://doi.org/10.1017/CBO9781316339831} {\emph
  {\bibinfo {title} {Introduction to Random Graphs}}}\ (\bibinfo  {publisher}
  {Cambridge University Press},\ \bibinfo {address} {Cambridge},\ \bibinfo
  {year} {2016})\BibitemShut {NoStop}%
\bibitem [{\citenamefont {Holmgren}\ and\ \citenamefont
  {Janson}(2015)}]{Janson15}%
  \BibitemOpen
  \bibfield  {author} {\bibinfo {author} {\bibfnamefont {C.}~\bibnamefont
  {Holmgren}}\ and\ \bibinfo {author} {\bibfnamefont {S.}~\bibnamefont
  {Janson}},\ }\bibfield  {title} {\enquote {\bibinfo {title} {Limit laws for
  functions of fringe trees for binary search trees and random recursive
  trees},}\ }\href {https://doi.org/10.1214/EJP.v20-3627} {\bibfield  {journal}
  {\bibinfo  {journal} {Electron. J. Probab.}\ }\textbf {\bibinfo {volume}
  {20}},\ \bibinfo {pages} {1--51} (\bibinfo {year} {2015})}\BibitemShut
  {NoStop}%
\bibitem [{\citenamefont {Janson}(2019)}]{Janson19}%
  \BibitemOpen
  \bibfield  {author} {\bibinfo {author} {\bibfnamefont {S.}~\bibnamefont
  {Janson}},\ }\bibfield  {title} {\enquote {\bibinfo {title} {Random recursive
  trees and preferential attachment trees are random split trees},}\ }\href
  {https://doi.org/10.1017/S0963548318000226} {\bibfield  {journal} {\bibinfo
  {journal} {Comb. Probab. Comput.}\ }\textbf {\bibinfo {volume} {28}},\
  \bibinfo {pages} {81--99} (\bibinfo {year} {2019})}\BibitemShut {NoStop}%
\bibitem [{\citenamefont {Godsil}\ and\ \citenamefont
  {Gutman}(1981)}]{godsil1981theory}%
  \BibitemOpen
  \bibfield  {author} {\bibinfo {author} {\bibfnamefont {C.~D.}\ \bibnamefont
  {Godsil}}\ and\ \bibinfo {author} {\bibfnamefont {I.}~\bibnamefont
  {Gutman}},\ }\bibfield  {title} {\enquote {\bibinfo {title} {On the theory of
  the matching polynomial},}\ }\href {https://doi.org/10.1002/jgt.3190050203}
  {\bibfield  {journal} {\bibinfo  {journal} {J. Graph Theory}\ }\textbf
  {\bibinfo {volume} {5}},\ \bibinfo {pages} {137--144} (\bibinfo {year}
  {1981})}\BibitemShut {NoStop}%
\bibitem [{\citenamefont {Diestel}(2017)}]{Diestel}%
  \BibitemOpen
  \bibfield  {author} {\bibinfo {author} {\bibfnamefont {R.}~\bibnamefont
  {Diestel}},\ }\href {https://doi.org/10.1007/978-3-662-53622-3} {\emph
  {\bibinfo {title} {Graph Theory}}}\ (\bibinfo  {publisher} {Springer},\
  \bibinfo {address} {Heidelberg},\ \bibinfo {year} {2017})\BibitemShut
  {NoStop}%
\bibitem [{\citenamefont {Fiedler}(1973)}]{fiedler1973algebraic}%
  \BibitemOpen
  \bibfield  {author} {\bibinfo {author} {\bibfnamefont {M.}~\bibnamefont
  {Fiedler}},\ }\bibfield  {title} {\enquote {\bibinfo {title} {Algebraic
  connectivity of graphs},}\ }\href {\doibase 10.21136/CMJ.1973.101168}
  {\bibfield  {journal} {\bibinfo  {journal} {Czech. Math. J.}\ }\textbf
  {\bibinfo {volume} {23}},\ \bibinfo {pages} {298--305} (\bibinfo {year}
  {1973})}\BibitemShut {NoStop}%
\bibitem [{\citenamefont {Durhuus}(2009)}]{durhuus2009hausdorff}%
  \BibitemOpen
  \bibfield  {author} {\bibinfo {author} {\bibfnamefont {B.}~\bibnamefont
  {Durhuus}},\ }\bibfield  {title} {\enquote {\bibinfo {title} {Hausdorff and
  spectral dimension of infinite random graphs},}\ }\href@noop {} {\bibfield
  {journal} {\bibinfo  {journal} {Acta Phys. Polon. B}\ }\textbf {\bibinfo
  {volume} {40}},\ \bibinfo {pages} {3509} (\bibinfo {year}
  {2009})}\BibitemShut {NoStop}%
\bibitem [{\citenamefont {Golinelli}(2003)}]{golinelli2003statistics}%
  \BibitemOpen
  \bibfield  {author} {\bibinfo {author} {\bibfnamefont {O.}~\bibnamefont
  {Golinelli}},\ }\bibfield  {title} {\enquote {\bibinfo {title} {Statistics of
  delta peaks in the spectral density of large random trees},}\ }\href
  {https://arxiv.org/abs/cond-mat/0301437} {\bibfield  {journal} {\bibinfo
  {journal} {arXiv:cond-mat/0301437}\ } (\bibinfo {year} {2003})}\BibitemShut
  {NoStop}%
\bibitem [{\citenamefont {Bhamidi}\ \emph {et~al.}(2012)\citenamefont
  {Bhamidi}, \citenamefont {Evans},\ and\ \citenamefont
  {Sen}}]{bhamidi2012spectra}%
  \BibitemOpen
  \bibfield  {author} {\bibinfo {author} {\bibfnamefont {S.}~\bibnamefont
  {Bhamidi}}, \bibinfo {author} {\bibfnamefont {S.~N.}\ \bibnamefont {Evans}},
  \ and\ \bibinfo {author} {\bibfnamefont {A.}~\bibnamefont {Sen}},\ }\bibfield
   {title} {\enquote {\bibinfo {title} {Spectra of large random trees},}\
  }\href {https://doi.org/10.1007/s10959-011-0360-9} {\bibfield  {journal}
  {\bibinfo  {journal} {J. Theor. Probab.}\ }\textbf {\bibinfo {volume} {25}},\
  \bibinfo {pages} {613} (\bibinfo {year} {2012})}\BibitemShut {NoStop}%
\bibitem [{\citenamefont {Destri}\ and\ \citenamefont
  {Donetti}(2002)}]{destri2002spectral}%
  \BibitemOpen
  \bibfield  {author} {\bibinfo {author} {\bibfnamefont {C.}~\bibnamefont
  {Destri}}\ and\ \bibinfo {author} {\bibfnamefont {L.}~\bibnamefont
  {Donetti}},\ }\bibfield  {title} {\enquote {\bibinfo {title} {The spectral
  dimension of random trees},}\ }\href {\doibase 10.1088/0305-4470/35/45/301}
  {\bibfield  {journal} {\bibinfo  {journal} {J. Phys. A: Math. Gen.}\ }\textbf
  {\bibinfo {volume} {35}},\ \bibinfo {pages} {9499} (\bibinfo {year}
  {2002})}\BibitemShut {NoStop}%
\end{thebibliography}%

\end{document}